\documentclass[11pt,reqno]{amsart}

\usepackage{amsmath, amsfonts, amsthm, amssymb, color, comment}

\newtheorem{Theorem}{Theorem}[section]
\newtheorem{Lemma}{Lemma}[section]

\theoremstyle{definition}

\newtheorem{Remark}{Remark}[section]

\numberwithin{equation}{section}

\allowdisplaybreaks

\renewcommand{\u}{{\bf u}}

\newcommand{\R}{{\mathbb R}}
\newcommand{\C}{{\mathbb C}}
\newcommand{\Dv}{{\rm div}}

\newcommand{\V}{{\mathcal V}}
\newcommand{\D}{{\mathcal D}}

\def\f{\frac}
\renewcommand{\O}{\Omega}
\renewcommand{\o}{\omega}

\def\hf1{^\f{1}{1-\xi^2}}

\def\be{\begin{equation}}
\def\en{\end{equation}}
\def\bs{\begin{split}}
\def\es{\end{split}}
\def\rb{\begin{Remark}}
\def\re{\end{Remark}}

\newcommand{\F}{{\mathtt F}}

\newcommand{\M}{{\bf M}}

\newcommand{\e}{\varepsilon}
\newcommand{\p}{\partial}

\newcommand{\nc}{\textrm{nc}}
\newcommand{\n}{\textrm{n}}

\newcommand{\vertiii}[1]{{\left\vert\kern-0.25ex\left\vert\kern-0.25ex\left\vert #1
    \right\vert\kern-0.25ex\right\vert\kern-0.25ex\right\vert}}
\newcommand{\ti}{i}

\author[R. M. Chen]{Robin Ming Chen}
\address{Department of Mathematics, University of Pittsburgh,
                           Pittsburgh, PA 15260.}
\email{mingchen@pitt.edu}

\author[J. Hu]{Jilong Hu}
\address{Department of Mathematics, University of Pittsburgh,
                           Pittsburgh, PA 15260.}
\email{jih62@pitt.edu}

\author[D. Wang]{Dehua Wang}
\address{Department of Mathematics, University of Pittsburgh,
                           Pittsburgh, PA 15260.}
\email{dwang@math.pitt.edu}

\title
[Stability of Vortex Sheets in Elastodynamics]
{Linear Stability of Compressible Vortex Sheets in 2D Elastodynamics: Variable Coefficients}


\date{\today}

\begin{document}

\begin{abstract}
The linear stability with variable coefficients of the vortex sheets for the two-dimensional compressible elastic flows is studied. 
As in our earlier work \cite{chen2015linear} on the linear stability with constant coefficients, 
the problem has a free boundary which is characteristic, and also the Kreiss-Lopatinskii condition  is not uniformly satisfied. 
In addition, the roots of the Lopatinskii determinant of the para-linearized system may coincide with the poles of the system. Such a new collapsing phenomenon causes serious difficulties when applying the bicharacteristic extension method of \cite{coulombel2002weak,coulombel2004weakly,coulombel2004stability}. Motivated by our method introduced in the constant-coefficient case \cite{chen2015linear}, we perform an upper triangularization to the para-linearized system to separate the outgoing mode into a closed form where the outgoing mode only appears at the leading order. This procedure results in a gain of regularity for the outgoing mode,  which allows us to overcome the loss of regularity of the characteristic components at the poles and hence to close all the energy estimates. We find that, analogous to the constant-coefficient case, elasticity generates notable stabilization effects, and there are additional stable subsonic regions compared with the isentropic Euler flows. Moreover, since our method does not rely on the construction of the bicharacterisic curves, it can also be applied to other fluid models such as the non-isentropic Euler equations and the MHD equations. 
\end{abstract}

\maketitle

\tableofcontents

\section{Introduction}
 
In this paper, we continue our study of   stability for vortex sheets in inviscid compressible elastic fluid flows initiated in \cite{chen2015linear}. In particular, we consider the vortex sheets in two-dimensional elastic fluids which can be described by the following system \cite{dafermos2010hyperbolic, gurtin1981introduction, joseph1990fluid}:
\begin{equation}\label{cons form}
\begin{cases}
\rho_{t}+\Dv(\rho\u)=0, \\
(\rho\u)_{t}+\Dv(\rho\u\otimes\u)+\nabla p
-\Dv(\rho\F\F^\top)=0,\\
(\rho\F_j)_t+\Dv(\rho\F_j\otimes\u-\u\otimes\rho\F_j)=0,
\end{cases}
\end{equation}
where $\rho$ is the density, $\u = (v,u)^\top \in\R^2$ is the velocity, $\F_j$  ($j=1,2$)  is the $j$th column of  the deformation gradient $\F=(F_{ij})\in\M^{2\times2}$, and $p$ is the pressure with $p=p(\rho)$ a smooth strictly increasing function on $(0,+\infty)$. To obtain this divergence form, we took the advantage of the following intrinsic property \cite{hu2011global}:
\begin{equation}\label{div free F}
\Dv(\rho\F_j)=0 \quad \text{for } j=1,2,
\end{equation}
which holds for all $t>0$, if it is satisfied initially. For a brief review of the history of the study of vortex sheets please refer to \cite{chen2015linear} and the references therein. 

The purpose of this paper is to establish the linear stability of compressible vortex sheets with variable coefficients for \eqref{cons form}.
In our previous paper \cite{chen2015linear}, we linearized the system \eqref{cons form} around a class of constant background states, which leads to a linear system with constant coefficients. For this linear system with constant coefficients, we obtained a constant-coefficient ODE system by the Fourier transform. Then by carefully examining the spectrum of the ODEs with respect to the initial conditions, we proved the linear stability of the system. In this paper, we take a step forward to linearize \eqref{cons form} around some {\it non-constant} background states near the constant states studied in \cite{chen2015linear}. Specifically, we follow the idea of Francheteau-M\'etivier \cite{francheteau2000existence} to straighten the wave front by a change of variable that satisfies the eikonal equations (cf. \eqref{variableeikonal}),
which results in a linear system with variable coefficients. Such a reformulation is different from those in Trakhinin \cite{trakhinin2005existencevairable, trakhinin2005existencelinearized}. Our goal is to analyze the spectrum and establish the stability for the new linearized system. With this important step and result in hand, we are able to finally prove in a subsequent paper \cite{chen2018nonlinear} the nonlinear stability and hence the local existence of  the vortex sheets for \eqref{cons form} in the two-dimensional case. 
We remark that the study of the linear and nonlinear stability of the  vortex sheets in the three-dimensional compressible elastic fluids is even more challenging due to the more complicated structure of the system, and will appear in a forthcoming paper.

As addressed in \cite{chen2015linear}, a common difficulty in the compressible vortex sheet problem is that the system has a free boundary and the free boundary is characteristic. This suggests that one would expect the loss of control on the trace of characteristic part of the solution. Moreover, the Kreiss-Lopatinskii condition in this problem is not uniformly satisfied. For the case of constant-coefficient system, the use of Fourier transform allows us to perform an upper triangularization (\cite{chen2015linear}) of the system pointwise in the frequency space. This way we only need to estimate the incoming modes, which can be derived directly from the Lopantiskii determinant. Due to the degeneracy induced by the non-uniformity of the Kreiss-Lopatinskii condition, the stability result exhibits loss of derivatives.

When dealing with the variable-coefficient system, however, a direct Fourier transform would lead to a non-local ODE system, making the analysis much more complicated. Instead, a common practice is to use para-differential calculus, which can be thought of as a generalization of the exact calculus of Fourier multipliers; see, for example \cite{Serre2000}. 
Unlike the classical Fourier multiplier method, the para-linearization is performed in the physical space, and hence it prohibits one to treat the system pointwise in the frequency space. A natural approach is to construct a weight function to characterize the degeneracy of the Kreiss-Lopatinskii condition on the boundary in the frequency space.

More specifically, in order to obtain the linear stability result, one needs to combine the differential equations in the interior domain with the boundary conditions. This is straightforward in the constant-coefficient case. But there is an extra difficulty in the variable-coefficient case, namely when performing the energy estimates to the para-linearized system, the commutator generated from the symmetrization is of the leading order, preventing one from closing the energy estimates. Coulombel \cite{coulombel2002weak,coulombel2004weakly} and Coulombel-Secchi \cite{coulombel2004stability} introduced a very intelligent  way 
to extend the weight function in the frequency space along the bicharacteristic curves emanating from the boundary. Such an extension provides crucial cancellation effect at the leading order commutator terms in the energy estimates. Based on such a technique they successfully proved the stability of shock waves \cite{coulombel2002weak, coulombel2004weakly} and isentropic compressible vortex sheets \cite{coulombel2004stability}. This idea has also been adapted in the study of linear stability of other models, such as two-phase flows \cite{ruanrectilinear}, steady isentropic flows \cite{wang2013stability},
nonisentropic flows \cite{morando2008two}, relativistic vortex sheets \cite{chen2017relativistic}, and so on. However, for this method to work, one requires that the ODE system be diagonalizable along the bicharacteristic curves which stem from the points near the roots of the Lopatinskii determinant. 

On the other hand, the characteristic boundary introduces an algebraic system for the characteristic part of the unknowns which may fail to be uniquely solvable at a set of boundary points, namely the poles.  In addition, the above mentioned ODE system becomes non-diagonalizable near those poles. Therefore it is unclear whether the above bicharacteristic extension method can work when the roots of the Lopatinskii determinant 
coincide with the poles of the system, which indeed happens in our case due to the elastic effect. This collapsing phenomenon is also observed in the three-dimensional compressible steady flows \cite{wang2015stability}, where the authors managed to prove the linear stability of contact discontinuity without exploring the cancellation effect as in \cite{coulombel2002weak, coulombel2004weakly,coulombel2004stability}, at the price of closing the estimates in a rougher regularity class.

We propose here a different approach which does not rely on the cancellation effect as well, but can provide sharper estimates for para-linearized system in a finer regularity space. Motivated by our upper triangularization method in the constant-coefficient case \cite{chen2015linear}, we can also upper triangularize our para-linearized system in terms of the symbols. Such a procedure would separate the {\it outgoing mode} into a closed form  where the outgoing mode only appears at the leading order. Roughly speaking, the upper-triangular system that we obtain by the para-differential calculus formally reads as
\begin{equation}\label{formal system}
I_2\frac{d}{dx_2}\begin{pmatrix}
W^{c}\\
W^{in}\\
W^{out}
\end{pmatrix}=
\begin{pmatrix}
\mathbb{F}  & 0 & *\\
*  & \mathbb{G} & 0\\
0 & 0  & \mathbb{H}\\
\end{pmatrix}\begin{pmatrix}
W^{c}\\
W^{in}\\
W^{out}
\end{pmatrix}+
\text{ lower order  terms},
\end{equation}
where $W^{in}$ is the incoming mode, $W^{out}$ is the outgoing mode, $W^{c}$ denotes the characteristic components and $I_2$ is a diagonal matrix with diagonal entries being 1 at the last two diagonal places and 0 elsewhere. Here the matrix on the right-hand side is understood as symbols of para-differential operation.
Specifically, the crucial point is that the symbol $\mathbb{H}$ is a scalar valued symbol with positive real part. This way reading off the equation for the outgoing mode $W^{out}$:
\begin{equation*}
\frac{d}{dx_2}W^{out}=\mathbb{H}W^{out}+ \text{ lower order terms},
\end{equation*}
one may bound $W^{out}$ in terms of the lower order terms. Moreover, performing this separation of modes allows us to extract a better property of $\mathbb{H}$, namely,  $\mathbb{H}$ is in fact homogeneous of degree one near the roots of the Lopatinskii determinant and the poles.  Hence we can obtain an improved regularity estimate for $W^{out}$. As a result, when plugging such an estimate into the system for $W^{c}$, since the symbol $\mathbb{F}$ is coercive in the sense that $\langle f, \mathbb{F} f \rangle \gtrsim \|f\|_{L^2}^2$, one can bound $W^{c}$ by the lower order terms as well. Combining the boundary conditions, the weight function and the control of $W^c$ will yield the desired estimates for $W^{in}$. This gain of regularity of the outgoing mode is  analogous to the constant-coefficient case where the closed system for the outgoing modes does not even have the lower order part, implying that the outgoing modes indeed all vanish. 

We note that our new argument does not require the explicit solvability of $W^c$. Indeed the symbol $\mathbb{F}$ is not invertible in the sense of para-differential calculus. Therefore the method can be applied to treat the poles.  
We also want to point out that higher-order estimates for $W^c$ can be obtained in the region away from the poles, whereas near the poles we do not in general expect this gain of regularity of $W^c$. This is again due to the fact that the symbol $\mathbb{F}$ is not invertible. It is the improved regularity of $W^{out}$ that helps us overcome this loss of higher regularity of $W^c$ at the poles.

Another advantage of our method is that the extension of the weight function does not rely on the construction of the bicharacterisic curves. Instead, we explore the relation between the  weight function and the background states. In particular we observe that the expression of the weight function on the boundary only involves the values of the background states, not the position of the points. This gives us the flexibility to extend the weight function by extending that relation between the  weight function and the background states into the interior of the domain; cf. Section \ref{subsec_poles and roots}, \ref{subsec_roots}. On the other hand, the poles of the system are also determined by an equation whose coefficients only depend on the background states. In this paper we will only consider the case that when the roots of the Lopantiskii determinant  coincide with the poles on the boundary, they share the same relation in terms of the background states. In this case we see that the roots of the weight function propagate along the poles. Therefore we are able to avoid the discrepancy between the bicharacteristic extension and the distribution of the poles. Note that our way of extending the weight function can also be applied to classical Euler and MHD cases where the roots of the Lopantiskii determinant  do not coincide with poles of the system.

The rest of the paper is organized as follows. In Section \ref{pre}, we review the nonlinear system of vortex sheets, linearize the system around a class of non-constant vortex sheets and introduce the necessary notations in the proof and our main theorem. In Section \ref{reductionofthesystem}, we reduce the problem  only for estimating some boundary component of the unknowns and use para-differential calculus to obtain an ODE system from which we expect to deduce those estimates. In Section \ref{microlocal}, we micro-localize the frequency space and identify the possible cases that we need to deal with. In Section \ref{Estimateineachcase}, for each case, we deduce the expected energy estimates. Finally in Section \ref{ProofofMain}, we combine all the energy estimates to obtain the main result.

\section{Preliminaries and Main Result}\label{pre}

In this section, we first review the nonlinear system which governs the dynamics of the vortex sheets in elastic fluids. Then we will introduce a class of perturbed states around the constant vortex sheets and linearize the nonlinear system around these states. Finally we collect some notations we will use in the main theorem and the proof. 

We assume that the vortex sheet is a free surface which is the graph $\Gamma$ of a function $x_2=\varphi(t,x_1)$ in $\R^3$, and denote the states on the two sides of the surface to be
\begin{equation*}
  U(t,x_1,x_2)=\begin{cases}
               U^+(t,x_1,x_2),\quad\quad \text{when $x_2>\varphi(t,x_1)$},\\
               U^-(t,x_1,x_2),\quad\quad \text{when $x_2<\varphi(t,x_1)$},
            \end{cases}
\end{equation*}
where $U^\pm = (\rho^\pm,\u^\pm,\F^\pm)$. Since this is a free boundary problem, as in the setup for the linear stability of constant coefficients \cite{chen2015linear}, we consider a class of functions $\Phi(t, x_1, x_2)$ such that $\inf \partial_2 \Phi > 0$ and $\Phi(t, x_1, 0) = \varphi(t, x_1)$. Then we define $\Phi^\pm(t,x_1,x_2) := \Phi(t,x_1, \pm x_2)$ for $x_2 \ge 0$. Now we can straighten the domain by introducing the functions
\begin{equation}\label{change}
U^+_{\#}(t, x_1, x_2) := U^+(t, x_1, \Phi^+(t,x_1,x_2)), \quad U^-_{\#}(t, x_1, x_2) := U^-(t, x_1, \Phi^-(t,x_1,x_2)).
\end{equation}
The new spatial domain is $\R^2_+:=\left\{(x_1,x_2)\in\R^2:\; x_2>0\right\},$
and the new domain $\Omega$ of our differential equations and its boundary $\o$ are defined respectively as:   
\begin{equation}\label{defn_domain}
\O:=\left\{(t,x_1,x_2)\in\R^3:\; x_2>0\right\},\quad \o:=\left\{(t,x_1,x_2)\in\O:\; x_2=0\right\}.
\end{equation}
In the following, we will only work on the fixed domain $\O$, and for convenience, we will drop the $\#$ subscript.

After the change of unknowns \eqref{change}, we obtain
\begin{equation}\label{differentialequation}
\begin{split}
\mathbb{L}(U^\pm,\Phi^\pm):=\p_tU^\pm&+A_1(U^\pm)\p_{1}U^\pm\\
&+\frac{1}{\p_{2}\Phi^\pm}\left[A_2(U^\pm)-\p_t\Phi^\pm I-\p_{1}\Phi^\pm A_3(U^\pm)\right]\p_{2}U^\pm=0,
\end{split}
\end{equation}
for $x_2>0$, where
\begin{align*}
A_1(U) =A_3(U)=
 \begin{pmatrix}
  v & \rho & 0 & 0 & 0& 0 & 0 \\
  \frac{p'}{\rho} & v & 0 & -F_{11} & 0 & -F_{12} & 0 \\
  0 & 0 & v & 0 & -F_{11} & 0 & -F_{12} \\
  0 & -F_{11} & 0 & v & 0 & 0 & 0 \\
  0 & 0 & -F_{11} & 0 & v & 0 & 0 \\
  0 & -F_{12} & 0 & 0 & 0 & v & 0 \\
  0 & 0 & -F_{12} & 0 & 0 & 0 & v \\
 \end{pmatrix}
\end{align*}
and
\begin{align*}
A_2(U) =
 \begin{pmatrix}
  u & 0 & \rho & 0 & 0& 0 & 0 \\
  0 & u & 0 & -F_{21} & 0 & -F_{22} & 0 \\
  \frac{p'}{\rho} & 0 & u & 0 & -F_{21} & 0 & -F_{22} \\
  0 & -F_{21} & 0 & u & 0 & 0 & 0 \\
  0 & 0 & -F_{21} & 0 & u & 0 & 0 \\
  0 & -F_{22} & 0 & 0 & 0 & u & 0 \\
  0 & 0 & -F_{22} & 0 & 0 & 0 & u \\
 \end{pmatrix}.
\end{align*}
Similarly as in the previous paper \cite{chen2015linear}, from the Rankine-Hugoniot conditions, we can obtain the boundary condition at $x_2=0$:
\begin{align}\label{boundarycondition}
\mathbb{B}(U|_{x_2=0},\varphi):=\begin{cases}
(v^+-v^-)\p_1\varphi-(u^+-u^-)=0, \\
\p_t\varphi+v^+\p_1\varphi-u^+=0, \\
(F_{11}^+-F_{11}^-)\p_1\varphi-(F_{21}^+-F_{21}^-)=0, \\
F_{11}^+\p_1\varphi-F_{21}^+=0, \\
(F_{12}^+-F_{12}^-)\p_1\varphi-(F_{22}^+-F_{22}^-)=0, \\
F_{12}^+\p_1\varphi-F_{22}^+=0, \\
\rho^+-\rho^-=0,
\end{cases}
\end{align}
where  we recall that $\Phi^\pm=\varphi$ at $x_2=0$.

\begin{Remark}
We note that  the above Rankine-Hugoniot conditions lead to $\text{det} F=0$,
which is the degenerate case corresponding to the vortex sheets in
elastodynamics. This degeneracy is due to the mathematical formulation
of the singular vortex sheet structure in elastic fluids. To the best
of the authors' knowledge the associated physical interpretation of
this degeneracy is still unclear.  On the other hand, the mathematical
methods we develop in this work can be applied to treat many other
physically relevant problems of vortex sheets arising from gas
dynamics, MHD, relativistic fluids, and so on. In fact, as mentioned
in the introduction, we can handle more  technical cases such as the
classical Euler and MHD systems where the roots of the Lopantiskii
determinant do not coincide with poles of the system.
\end{Remark}

 For the linear stability of the variable coefficient case, we consider the following background states:
\begin{equation}\label{variablestates}
\begin{split}
&U^{r,l}=\begin{pmatrix}
\rho^{r,l} \\
v^{r,l} \\
u^{r,l} \\
F_{11}^{r,l} \\
F_{21}^{r,l} \\
F_{12}^{r,l} \\
F_{22}^{r,l}
\end{pmatrix}:= \bar{U}^{r,l}+\dot{U}^{r,l}=
\begin{pmatrix}
\bar{\rho} \\
\pm\bar{v} \\
0 \\
\pm\bar{F}_{11} \\
0 \\
\pm\bar{F}_{12} \\
0
\end{pmatrix}+\begin{pmatrix}
\dot{\rho}^{r,l} \\
\dot{v}^{r,l} \\
\dot{u}^{r,l} \\
\dot{F}_{11}^{r,l} \\
\dot{F}_{21}^{r,l} \\
\dot{F}_{12}^{r,l} \\
\dot{F}_{22}^{r,l}
\end{pmatrix}, \\
&\Phi^{r,l}(t,x_1,x_2):=\pm x_2+\dot{\Phi}^{r,l},
\end{split}
\end{equation}
where $U^{r,l}$ and $\Phi^{r,l}$ are states and changes of variables on each side of the vortex sheet respectively;  $\bar{\rho}>0$, $\bar{v}$, $\bar{F}_{11}$ and $\bar{F}_{12}$ are constants;  $\dot{U}^{r,l}$ and $\dot{\Phi}^{r,l}$ are functions which represent the perturbation around the constant states.

\begin{Remark}
When the perturbations $\dot{U}^{r,l}$ and $\dot{\Phi}^{r,l}$ are zero, the background states reduce to the rectilinear vortex sheets, and the linear stability around this kind of states has been discussed in the previous paper \cite{chen2015linear}.
\end{Remark}
Moreover, we require the perturbation of the background states to satisfy
\begin{align}\label{variablesmallness}
&\dot{U}^{r,l}\in W^{2,\infty}(\O), \dot{\Phi}^{r,l}\in W^{3,\infty}(\O),\quad\|(\dot{U}^r,\dot{U}^l)\|_{W^{2,\infty}(\O)}+\|(\dot{\Phi}^r,\dot{\Phi}^l)\|_{W^{3,\infty}(\O)}\leq K,
\end{align}
where $K$ is a suitable positive constant, $\dot{U}^{r,l}$ and $\dot{\Phi}^{r,l}$ have compact support in the domain. 
We also need the perturbed states \eqref{variablestates} to satisfy the Rankine-Hugoniot conditions: 
 \begin{align}\label{variableboundary}
 \begin{cases}
(v^r-v^l)\p_1\varphi-(u^r-u^l)=0, \\
\p_t\varphi+v^r\p_1\varphi-u^r=0, \\
(F_{11}^r-F_{11}^l)\p_1\varphi-(F_{21}^r-F_{21}^l)=0, \\
F_{11}^r\p_1\varphi-F_{21}^r=0, \\
(F_{12}^r-F_{12}^l)\p_1\varphi-(F_{22}^r-F_{22}^l)=0, \\
F_{12}^r\p_1\varphi-F_{22}^r=0, \\
\rho^r-\rho^l=0,
\end{cases}
\end{align}
on $x_2=0$, where $\varphi=\Phi^{r}|_{x_2=0}=\Phi^{l}|_{x_2=0}$. To keep the constant rank of boundary matrix in the whole domain and have a simpler formulation, inspired by \cite{coulombel2004stability, francheteau2000existence}, we assume that the following conditions on the perturbed states \eqref{variablestates} hold:
\begin{equation}\label{variableeikonal}
\begin{split}
\begin{cases}
\p_t\Phi^{r,l}+v^{r,l}\p_1\Phi^{r,l}-u^{r,l}=0,\\
F_{11}^{r,l}\p_1\Phi^{r,l}-F_{21}^{r,l}=0,\\
F_{12}^{r,l}\p_1\Phi^{r,l}-F_{22}^{r,l}=0,
\end{cases}
\end{split}
\end{equation}
and 
\begin{align}\label{variablenondegeneracy}
\p_2\Phi^r\geq \kappa_0 \text{ and } \p_2\Phi^l\leq-\kappa_0,
\end{align}
for all $(t,x)\in \O$ and  some positive constant $\kappa_0$.

Now we linearize the differential equation \eqref{differentialequation} around the states \eqref{variablestates} and denote by $(V^\pm, \Psi^\pm)$ the perturbation  of the states $(U^{r,l}, \Phi^{r,l})$. Then the linearized equations are 
\begin{equation*}
\begin{split}
&\p_t V^\pm+ A_1(U^{r,l})\p_1V^\pm+\frac{1}{\p_2\Phi^{r,l}}(A_2(U^{r,l})-\p_t\Phi^{r,l}-\p_1\Phi^{r,l}A_1(U^{r,l}))\p_2V^\pm\\
&+[dA_1(U^{r,l})V^\pm]\p_1U^{r,l}-\frac{\p_2\Psi^{\pm}}{(\p_2\Phi^{r,l})^2}(A_2(U^{r,l})-\p_t\Phi^{r,l}-\p_1\Phi^{r,l}A_1(U^{r,l}))\p_2U^{r,l}\\
&+\frac{1}{\p_2\Phi^{r,l}}[dA_2(U^{r,l})V^\pm-\p_t\Psi^{\pm}-\p_1\Psi^{\pm}A_1(U^{r,l})-\p_1\Phi^{r,l}dA_1(U^{r,l})V^\pm]\p_2U^{r,l}=f,
\end{split}
\end{equation*} 
for $x_2>0$. We define the first order linear operator
\begin{equation*}
\begin{split}
&L(U^{r,l},\nabla\Phi^{r,l})V^{\pm}:=\\
&\quad\quad\p_tV^{\pm}+A_1(U^{r,l})\p_1V^{\pm}+\frac{1}{\p_2\Phi^{r,l}}(A_2(U^{r,l})-\p_t\Phi^{r,l}-\p_1\Phi^{r,l}A_1(U^{r,l}))\p_2V^{\pm},
\end{split}
\end{equation*}
and introduce the Alinhac's `good unknown' \cite{alinhac1989existence}:
\begin{equation*}
\dot{V}^\pm=\begin{pmatrix}
\dot{\rho}^\pm \\
\dot{v}^\pm \\
\dot{u}^\pm \\
\dot{F}_{11}^\pm \\
\dot{F}_{21}^\pm \\
\dot{F}_{12}^\pm \\
\dot{F}_{22}^\pm
\end{pmatrix}:=V^{\pm}-\frac{\Psi^{\pm}}{\p_2\Phi^{r,l}}\p_2U^{r,l}.
\end{equation*}
Then we can rewrite the above equations as
\begin{equation*}
L(U^{r,l},\nabla\Phi^{r,l})\dot{V}^\pm+C(U^{r,l},\nabla U^{r,l}, \nabla\Phi^{r,l})\dot{V}^\pm+\frac{\Psi^{\pm}}{\p_2\Phi^{r,l}}\p_2[L(U^{r,l},\nabla\Phi^{r,l})U^{r,l}]=f^{r,l},
\end{equation*}
where
\begin{equation*}
\begin{split}
&C(U^{r,l},\nabla U^{r,l}, \nabla\Phi^{r,l})\dot{V}^\pm:=\\
&\qquad[dA_1(U^{r,l})\dot{V}^\pm]\p_1U^{r,l}+\frac{1}{\p_2\Phi^{r,l}}[dA_2(U^{r,l})\dot{V}^\pm-\p_1\Phi^{r,l}dA_1(U^{r,l})\dot{V}^\pm]\p_2U^{r,l}.
\end{split}
\end{equation*}
From the same argument as in \cite{chen2015linear}, we can neglect the zeroth order terms of $\Psi^{\pm}$ and only consider the following differential equations:
\begin{align}
L'_{r,l}\dot{V}^\pm:=L(U^{r,l},\nabla\Phi^{r,l})\dot{V}^\pm+C(U^{r,l},\nabla U^{r,l}, \nabla\Phi^{r,l})\dot{V}^\pm=f^{r,l}.
\end{align} 
Since $U^{r,l}$ are in $W^{2,\infty}(\O)$, we know that the coefficients in $L(U^{r,l},\nabla\Phi^{r,l})$ are in $W^{2,\infty}(\O)$ and the coefficients in $C(U^{r,l},\nabla U^{r,l},\nabla\Phi^{r,l})$ are in $W^{1,\infty}(\O)$. 

Now we linearize the boundary conditions around the same perturbed states to obtain
\begin{equation*}
\begin{split}
&(v^r-v^l)\p_1\psi+(v^+-v^-)\p_1\varphi-(u^+-u^-)=g_1,\\
&\p_t\psi+v^r\p_1\psi+v^+\p_1\varphi-u^+=g_2,\\
&\rho^+-\rho^-=g_3,\\
&(F_{11}^r-F_{11}^l)\p_1\psi+(F_{11}^+-F_{11}^-)\p_1\varphi-(F_{21}^+-F_{21}^-)=g_4,\\
&F_{11}^r\p_1\psi+F_{11}^+\p_1\varphi-F_{21}^+=g_5,\\
&(F_{12}^r-F_{12}^l)\p_1\psi+(F_{12}^+-F_{12}^-)\p_1\varphi-(F_{22}^+-F_{22}^-)=g_6,\\
&F_{12}^r\p_1\psi+F_{12}^+\p_1\varphi-F_{22}^+=g_7,
\end{split}
\end{equation*}
at $x_2=0$, where $\psi=\Psi^+|_{x_2=0}=\Psi^-|_{x_2=0}$. We denote the above equations in the following form:
\begin{equation*}
\underline{b}\nabla\psi+\underline{M}V|_{x_2=0}=g,
\end{equation*}
where $V=(V^+,V^-)^\top$, $\nabla\psi=(\p_t\psi,\p_1\psi)^\top$, 
$g=(g_1,g_2,g_3,g_4,g_5,g_6,g_7)^\top$,
\begin{equation*}
\underline{b}(t,x_1):=\begin{pmatrix}
0 & (v^r-v^l)|_{x_2=0}\\
1 & v^r|_{x_2=0}\\
0 & 0\\
0 & (F_{11}^r-F_{11}^l)|_{x_2=0}\\
0 & F_{11}^r|_{x_2=0}\\
0 & (F_{12}^r-F_{12}^l)|_{x_2=0}\\
0 & F_{12}^r|_{x_2=0}
\end{pmatrix},
\end{equation*}
and
\begin{equation*}
\begin{split}
&\underline{M}(t,x_1):=\\
&\left(
\begin{array}{*{20}c}
 0 & \p_1\varphi & -1 & 0 & 0 & 0 & 0 & 0 & -\p_1\varphi & 1 & 0 & 0 & 0 & 0\\
 0 & \p_1\varphi & -1 & 0 & 0 & 0 & 0 & 0 & 0 & 0 & 0 & 0 & 0 & 0\\
 1 & 0 & 0 & 0 & 0 & 0 & 0 & -1 & 0 & 0 & 0 & 0 & 0 & 0\\
 0 & 0 & 0 & \p_1\varphi & -1 & 0 & 0 & 0 & 0 & 0 & -\p_1\varphi & 1 & 0 & 0\\
 0 & 0 & 0 & \p_1\varphi & -1 & 0 & 0 & 0 & 0 & 0 & 0 & 0 & 0 & 0 \\
 0 & 0 & 0 & 0 & 0 & \p_1\varphi & -1 & 0 & 0 & 0 & 0 & 0 & -\p_1\varphi & 1 \\
 0 & 0 & 0 & 0 & 0 & \p_1\varphi & -1 & 0 & 0 & 0 & 0 & 0 & 0 & 0
 \end{array}
\right).
\end{split}
\end{equation*}
Also, by using the Alinhac's `good unknown', we obtain
\begin{align}
B'(\dot{V}, \psi):=\underline{b}\nabla\psi+\underline{M}\begin{pmatrix}
\p_2 U^r/\p_2\Phi^r\\
\p_2 U^l/\p_2\Phi^l
\end{pmatrix}\psi+
\underline{M}\dot{V}|_{x_2=0}=g.
\end{align}
Hence our problem is reduced to  the following linear one:
\begin{align}\label{variablesystem}
\begin{cases}
L'_{r,l}\dot{V}^\pm =f^{r,l},      \quad x_2>0,\\
B'(\dot{V}, \psi)=g,               \quad x_2=0.
\end{cases}
\end{align}
We notice that the boundary condition does not involove the tangential components of $\dot{V}|_{x_2=0}$. Thus we define the components of $\dot{V}|_{x_2=0}$ which are involved in the boundary condition by $\dot{V}^\n|_{x_2=0}$ as  
\begin{align*}
&\dot{V}^\textnormal{n}|_{x_2=0}=\left(\dot{\rho}^+,\dot{u}^+-\dot{v}^+\p_1\Phi^r,\dot{F}_{21}^+-\dot{F}_{11}^+\p_1\Phi^r,\dot{F}_{22}^+-\dot{F}_{12}^+\p_1\Phi^r,\right.\\
&\qquad\qquad\qquad\qquad\qquad\qquad\qquad\left.\dot{\rho}^-,\dot{u}^--\dot{v}^-\p_1\Phi^l,\dot{F}_{21}^--\dot{F}_{11}^-\p_1\Phi^l,\dot{F}_{22}^--\dot{F}_{12}^-\p_1\Phi^l \right)^\top|_{x_2=0}.
\end{align*}

To describe our main result, we introduce the following functional spaces. First we define the Sobolev space with the weight on the time variable as in \cite{coulombel2004stability},
\begin{align*}
& H^s_{\gamma}(\o):=\{u(t,x_1)\in\D'(\R^2)\;:\;e^{-\gamma t}u(t,x_1)\in H^s(\o)\},
\end{align*}
for $s\in \mathbb{N},\ \gamma\geq 1$, equipped with the norms
\begin{equation*}
\|u\|_{H^s_{\gamma}(\o)}:=\|e^{-\gamma t}u\|_{H^s(\o)},
\end{equation*}
where $\o$ is given in \eqref{defn_domain}. 
We notice that
the above norm has the following equivalent form   by the Fourier transform:
\begin{equation*}
\|\tilde{u}\|_{s,\gamma}^2=\frac{1}{(2\pi)^2}\int_{\R^2}(\gamma^2+|\xi|^2)^s|\hat{\tilde{u}}(\xi)|^2d\xi,\quad \forall\; u\in H^s(\o),
\end{equation*}
where $\hat{\tilde{u}}(\xi)$ is the Fourier transform of $\tilde{u}$ and $\tilde{u}:=e^{-\gamma t}u$. We denote the above statement by $\|u\|_{H^s_{\gamma}(\o)}\simeq\|\tilde{u}\|_{s,\gamma}$, where $\simeq$ is the equivalence relation of norms. Now we define the space $L^2(\R_+;H^s_\gamma(\o))$, equipped with the norm
\begin{equation*}
\vertiii{v}^2_{L^2(H^s_\gamma(\o))}=\int_0^{+\infty}\|v(\cdot,x_2)\|^2_{H^s_\gamma(\o)}dx_2.
\end{equation*}
Again we have 
\begin{equation*}
\vertiii{v}_{L^2(H^s_\gamma(\o))}^2\simeq\vertiii{\tilde{v}}_{s,\gamma}^2:=\int_0^{+\infty}\|\tilde{v}(\cdot,x_2)\|^2_{s,\gamma}dx_2.
\end{equation*}
Note that $\|\cdot\|_{0,\gamma}$ is actually the usual norm on $L^2(\o)$ and $\vertiii{\cdot}_{0,\gamma}$ is the usual norm on $L^2(\O)$. 
Our main theorem reads as follows.
\begin{Theorem}\label{variablemain}
Suppose that the particular solution defined by \eqref{variablestates} satisfies one of the following two conditions
\begin{align*}
&(i) \qquad\bar{v}^2> 2c(\bar{\rho})^2+\bar{F}_{11}^2+\bar{F}_{12}^2, \ \text{ or }\\
&(ii) \qquad\bar{v}^2< \bar{F}_{11}^2+\bar{F}_{12}^2 \quad\text{ but }\\
&\qquad\qquad\bar{v}^2\neq\frac{\bar{F}_{11}^2+\bar{F}_{12}^2}{4},\quad \bar{v}^2\neq\frac{\left(\sqrt{\bar{F}_{11}^2+\bar{F}_{12}^2+c(\bar{\rho})^2}-\sqrt{\bar{F}_{11}^2+\bar{F}_{12}^2}\right)^2}{4}, \\
&\qquad\qquad\bar{v}^2\neq\frac{\bar{F}_{11}^2+\bar{F}_{12}^2+c(\bar{\rho})^2}{4},\quad \bar{v}^2\neq\frac{\left(\bar{F}_{11}^2+\bar{F}_{12}^2\right)\left(2c(\bar{\rho})^2+\bar{F}_{11}^2+\bar{F}_{12}^2\right)}{4\left(\bar{F}_{11}^2+\bar{F}_{12}^2+c(\bar{\rho})^2\right)};
\end{align*}
and moreover,  the perturbations $\dot{U}^{r,l}$ and $\dot{\Phi}^{r,l}$ have compact support,  and $K$ in \eqref{variablesmallness} is small enough. 
Then there are two constants $C_0$ and $\gamma_0$ which are determined by the particular solution, such that for all $\dot{V}$ and $\psi$ and all $\gamma\geq\gamma_0$ the following estimate holds:
\begin{equation*}
\gamma\vertiii{\dot{V}}_{L^2(H^0_\gamma)}^2+\left\|\dot{V}^\textnormal{n}|_{x_2=0}\right\|_{L^2_\gamma(\R^2)}^2+\left\|\psi\right\|_{H^1_\gamma(\R^2)}^2\leq C_0\left(\frac{1}{\gamma^3}\vertiii{L'\dot{V}}_{L^2(H^1_\gamma)}^2+\frac{1}{\gamma^2}\left\|B'(\dot{V},\psi)\right\|_{H^1_\gamma(\R^2)}^2\right),
\end{equation*}
where $L'\dot{V}=(L'_r\dot{V}^+,L'_l\dot{V}^-)$. 
\end{Theorem}
\begin{Remark}
Note that in comparison with the Euler case (see \cite{coulombel2004weakly}), the stabilization effect from elasticity can be inferred from (ii) in the above theorem. Such a phenomenon has been observed in the constant-coefficient case \cite{chen2015linear}.
\end{Remark}
\begin{Remark}
Compared with the stability theorem in the constant-coefficient case (see \cite[Theorem 2.1 (1)]{chen2015linear}), here we impose more restrictions on the non-constant background states in part (ii) of Theorem \ref{variablemain}. This is because when applying the para-differential calculus to the variable-coefficient problem, one needs to separate certain special frequencies in order to establish a uniform extension of the weight function. See the end of Section \ref{microlocal} for a detailed discussion. Such a separation of frequencies is not necessary in the constant-coefficient case since a Fourier transform method can be performed pointwise in the frequency space. 
\end{Remark}

\section{Reduction of the System}\label{reductionofthesystem}

In this section, we are going to transform the system \eqref{variablesystem} into an ordinary differential system. This requires a linear transformation on the unknown variables and a para-linearization on the equations of the transformed variables. 

For the system \eqref{variablesystem}, it is not difficult to find a symmetrizer 
\begin{align*}
S^{r,l}:=\text{diag}\{p'(\rho^{r,l})/\rho^{r,l},\rho^{r,l},\rho^{r,l},1,1,1,1\}.
\end{align*}
We multiply $(S^{r,l}\dot{V}^\pm)^\top$  to the interior equations of \eqref{variablesystem} and integrate by parts to obtain the following lemma:
\begin{Lemma}\label{L2ofdotV}
There are two positive constants $C$ and $\gamma_1\geq 1$ such that for any $\gamma\geq \gamma_1$ the following holds:
\begin{align*}
\gamma\vertiii{\dot{V}^\pm}_{L^2(H^0_\gamma)}^2\leq C\left(\frac{1}{\gamma}\vertiii{L'_{r,l}\dot{V}^\pm}_{L^2_\gamma(\R^3_+)}^2+\left\|\dot{V}^\n|_{x_2=0}\right\|_{L_\gamma^2(\R^2)}^2\right).
\end{align*}

\end{Lemma}

With the help of this lemma, we only need to estimate $\dot{V}^\textnormal{n}|_{x_2=0}$ and $\psi$ through the system \eqref{variablesystem}. To achieve this we will investigate this system carefully 
using para-linearization. Before we do that, we need to transform the linear operator $L(U^{r,l},\nabla\Phi^{r,l})$ to make the boundary matrix
\begin{align*}
\tilde{A}^{r,l}_2:=\frac{1}{\p_2\Phi^{r,l}}(A_2(U^{r,l})-\p_t\Phi^{r,l}-\p_1\Phi^{r,l}A_1(U^{r,l}))
\end{align*}
reduced to a constant diagonal matrix. This transformation can always be expected by observing that $\tilde{A}^{r,l}_2$ has a constant rank in the closure of the domain $\R^3_+$. In particular, we consider the following transformation matrix
\begin{align*}
&T(U^{r,l},\nabla\Phi^{r,l}):=\\
&\begin{pmatrix}
0 &  \langle\p_1\Phi^{r,l}\rangle & \langle\p_1\Phi^{r,l}\rangle & 0 & 0 & 0 & 0 \\
1 &  -\frac{c(\rho^{r,l})}{\rho^{r,l}}\p_1\Phi^{r,l} & \frac{c(\rho^{r,l})}{\rho^{r,l}}\p_1\Phi^{r,l} & 0 & 0 & 0 & 0\\
\p_1\Phi^{r,l} & \frac{c(\rho^{r,l})}{\rho^{r,l}} & -\frac{c(\rho^{r,l})}{\rho^{r,l}} & 0 & 0 & 0 & 0\\ 
0 & 0 & 0 & 1 & -\p_1\Phi^{r,l} & 0 & 0\\
0 & 0 & 0 & \p_1\Phi^{r,l} & 1 & 0 & 0\\
0 & 0 & 0 & 0 & 0 & 1 & -\p_1\Phi^{r,l}\\
0 & 0 & 0 & 0 & 0 & \p_1\Phi^{r,l} & 1
\end{pmatrix},
\end{align*}
where $\langle\p_1\Phi^{r,l}\rangle=\sqrt{1+{\p_1\Phi^{r,l}}^2}$. This diagonalizes the boundary matrix $\tilde{A}^{r,l}$ as
\begin{align*}
T^{-1}(U^{r,l},\nabla\Phi^{r,l})\tilde{A}^{r,l}_2T(U^{r,l},\nabla\Phi^{r,l})=\text{diag}\{0,\frac{c(\rho^{r,l})\langle\p_1\Phi^{r,l}\rangle}{\p_2\Phi^{r,l}},-\frac{c(\rho^{r,l})\langle\p_1\Phi^{r,l}\rangle}{\p_2\Phi^{r,l}} , 0,0,0,0\}.
\end{align*}
Once multiplying the above by the following matrix:
\begin{align*}
A_0^{r,l}=\text{diag}\{1,\frac{\p_2\Phi^{r,l}}{c(\rho^{r,l})\langle\p_1\Phi^{r,l}\rangle},-\frac{\p_2\Phi^{r,l}}{c(\rho^{r,l})\langle\p_1\Phi^{r,l}\rangle} , 1,1,1,1\},
\end{align*}
 we can make the boundary matrix to be a constant matrix,  in particular, in the following form:
 \begin{align*}
 I_2=\text{diag}\{0,1,1,0,0,0,0\}
 \end{align*}
With the above transformation matrix, the interior equations of $\eqref{variablesystem}$ for the new unknowns $W^\pm:=T^{-1}(U^{r,l},\nabla\Phi^{r,l})\dot{V}^\pm$ are
\begin{align}\label{Tsystem}
A_0^{r,l}\p_tW^{\pm}+A_1^{r,l}\p_1W^\pm+I_2\p_2W^\pm+A_0^{r,l}C^{r,l}W^\pm=F^{r,l},
\end{align}
where
\begin{align*}
A_1^{r,l}=&A_0^{r,l}T^{-1}(U^{r,l},\nabla\Phi^{r,l})A_1(U^{r,l})T(U^{r,l},\nabla\Phi^{r,l}),\\
C^{r,l}=&T^{-1}(U^{r,l},\nabla\Phi^{r,l})\p_tT(U^{r,l},\nabla\Phi^{r,l})+T^{-1}(U^{r,l},\nabla\Phi^{r,l})A_1(U^{r,l})\p_1T(U^{r,l},\nabla\Phi^{r,l})\\
&+T^{-1}(U^{r,l},\nabla\Phi^{r,l})C(U^{r,l},\nabla U^{r,l},\nabla\Phi^{r,l})T(U^{r,l},\nabla\Phi^{r,l})\\
&+T^{-1}(U^{r,l},\nabla\Phi^{r,l})\tilde{A}_2^{r,l}\p_2T(U^{r,l},\nabla\Phi^{r,l}),\\
F^{r,l}=&A_0^{r,l}T^{-1}(U^{r,l},\nabla\Phi^{r,l})f^{r,l}.
\end{align*}

To simplify the argument in obtaining the estimates in the weighted Sobolev norm,  we consider the weighted unknown $\widetilde{W}^\pm=e^{-\gamma t}W^\pm$ and rewrite  \eqref{Tsystem} as
\begin{align*}
\mathcal{L}_{r,l}^\gamma\widetilde{W}^\pm:=\gamma A_0^{r,l}\widetilde{W}^\pm+A_0^{r,l}\p_t\widetilde{W}^\pm+A_1^{r,l}\p_1\widetilde{W}^\pm+I_2\p_2\widetilde{W}^\pm+A_0^{r,l}C^{r,l}\widetilde{W}^\pm=e^{-\gamma t}F^{r,l},
\end{align*}
where $A_j^{r,l}\in W^{2,\infty}(\O)$ and $C^{r,l}\in W^{1,\infty}(\O)$. 
Similarly, applying  the transformation $T$, we rewrite the boundary condition of $\eqref{variablesystem}$ as
\begin{align}\label{Tboundary}
\underline{b}\nabla\psi+\underline{M}\begin{pmatrix}
\p_2 U^r/\p_2\Phi^r\\
\p_2 U^l/\p_2\Phi^l
\end{pmatrix}\psi+
\underline{M}
\begin{pmatrix}
T(U^{r},\nabla\Phi^{r}) & 0\\
0 & T(U^{l},\nabla\Phi^{l})
\end{pmatrix}
W|_{x_2=0}=g.
\end{align} 
In terms of $\widetilde{W}=e^{-\gamma t}W$ and $\widetilde{\psi}=e^{-\gamma t}\psi$, we obtain
\begin{align*}
\mathcal{B}^\gamma(\widetilde{W}|_{x_2=0},\widetilde{\psi}):=\gamma b_0\widetilde{\psi}+\underline{b}\nabla\widetilde{\psi}+\underline{M}\begin{pmatrix}
\p_2 U^r/\p_2\Phi^r\\
\p_2 U^l/\p_2\Phi^l
\end{pmatrix}\widetilde{\psi}+
\underline{M}\begin{pmatrix}
T^r & 0\\
0 & T^l
\end{pmatrix}
\widetilde{W}|_{x_2=0}=e^{-\gamma t}g,
\end{align*}
where $b_0$ is the first column of $\underline{b}$. Obviously, we have $\underline{b}$, $\underline{M}$ and $T$  all in $W^{2,\infty}(\O)$, and
\begin{align*}
\check{b}:=\underline{M}\begin{pmatrix}
\p_2 U^r/\p_2\Phi^r\\
\p_2 U^l/\p_2\Phi^l
\end{pmatrix}\in W^{1,\infty}(\O).
\end{align*}
As we pointed out just before Theorem \ref{variablemain}, the boundary condition $\mathcal{B}^\gamma(\widetilde{W}|_{x_2=0},\widetilde{\psi})$ only involves part of components of $\widetilde{W}$. By examining the matrix coefficient in front of $\widetilde{W}$ in $\mathcal{B}^\gamma(\widetilde{W}|_{x_2=0},\widetilde{\psi})$, we denote the involved components by $\widetilde{W}^{\n}$ and obtain
\begin{align}
\widetilde{W}^{\n}=(\widetilde{W}_2,\widetilde{W}_3,\widetilde{W}_5,\widetilde{W}_7,\widetilde{W}_9,\widetilde{W}_{10},\widetilde{W}_{12},\widetilde{W}_{14}),
\end{align}
which are actually the normal components of $\widetilde{W}$ on the vortex sheets. Thus we rewrite the boundary condition as  
\begin{align*}
\mathcal{B}^\gamma(\widetilde{W}^\n|_{x_2=0},\widetilde{\psi})=e^{-\gamma t}g.
\end{align*}
Then by denoting $\widetilde{F}^{r,l}=e^{-\gamma t}F^{r,l}$ and $\widetilde{g}=e^{-\gamma t}g$, we rewrite the system \eqref{variablesystem} as
\begin{equation}\label{variableweightedsystem}
\begin{cases}
\mathcal{L}_{r,l}^\gamma \widetilde{W}^{\pm}=\widetilde{F}^{r,l},\\
\mathcal{B}^\gamma(\widetilde{W}^\n|_{x_2=0},\widetilde{\psi})=\widetilde{g}.
\end{cases}
\end{equation}

Similarly to the constant coefficient case, the key idea in proving Theorem \ref{variablemain} is to transform the variable coefficient linear system \eqref{variableweightedsystem} into an ODE. Instead of using Fourier transformation as in the constant-coefficient case, we perform the para-linearization to \eqref{variableweightedsystem}. Thus in the rest of this section, we derive the para-linearized system of \eqref{variableweightedsystem} and estimate the errors in replacing the system \eqref{variableweightedsystem} by its para-linearization. On para-differential calculus, we refer to the appendix of \cite{coulombel2004stability} for a list of properties, and \cite{bony1981calcul, meyer1981remarques} for a more detailed theory.

First we recall the frequency space defined in the constant coefficient case \cite{chen2015linear}:
\begin{align*}
\Pi:=\{(\tau,\eta):\;\tau=\gamma+i\delta\in\C,\ \eta\in\R,\ |\tau|^2+\eta^2\neq0,\ \Re\tau\geq 0\},
\end{align*}
where $\delta$ and $\eta$ are the Fourier variables with respect to $t$ and $x_1$ respectively. Since we need to utilize the homogeneity structure of the system, many of our arguments will be developed on the following unit hemisphere in the frequency space $\Pi$:
\begin{align*}
\Sigma=\{(\tau,\eta):\;|\tau|^2+\eta^2=1\text{ and } \Re\tau\geq 0\},
\end{align*}
and then will be extended to the whole frequency space $\Pi$.

To simplify notations, we will drop the tilde in the system for the rest of this section. In the boundary conditions, we denote 
\begin{align*}
b_0:=\begin{pmatrix}
0\\
1\\
0\\
0\\
0\\
0\\
0
\end{pmatrix}, \quad 
b_1(t, x_1):=\begin{pmatrix}
v^r-v^l\\
v^r\\
0\\
F_{11}^r-F_{11}^l\\
F_{11}^r\\
F_{12}^r-F_{12}^l\\
F_{12}^r
\end{pmatrix}(t,x_1,0),\quad \text{ and }
b=\tau b_0+\ti\eta b_1.
\end{align*}
From the theory of para-linearization \cite[Theorem \textbf{B.9.}]{coulombel2004stability}, we have
\begin{align*}
&\gamma b_0 \psi+b_0\p_t\psi=T_{\tau b_0}^\gamma \psi,\\
&\|b_1\p_1\psi-T_{\ti\eta b_1}^\gamma\psi\|_{1,\gamma}\leq C\|b_1\|_{W^{2,\infty}(\R^2)}\|\psi\|_0\leq\frac{C}{\gamma}\|\psi\|_{1,\gamma},\\
&\|\check{b}\psi-T_{\check{b}}^\gamma\psi\|_{1,\gamma}\leq C\|\check{b}\|_{W^{1,\infty}(\R^2)}\|\psi\|_0\leq\frac{C}{\gamma}\|\psi\|_{1,\gamma},\\
&\|T^\gamma_{\check{b}}\psi\|_{1,\gamma}\leq C\|\check{b}\|_{L^{\infty}(\R^2)}\|\psi\|_{1,\gamma}\leq C\|\psi\|_{1,\gamma},
\end{align*}
for some positive constant $C$. Then we consider the coefficients of $W^\n$,
\begin{align*}
&\underline{M}\text{diag}\{T^r,T^l\}W=:\mathbf{M}W^\n=\\
&{\small \left(\begin{array}{*{20}c}
& -\frac{c^r}{\rho^r}\langle\p_1\varphi\rangle^2 & \frac{c^r}{\rho^r}\langle\p_1\varphi\rangle^2 & 0 & 0 & \frac{c^l}{\rho^l}\langle\p_1\varphi\rangle^2 & -\frac{c^l}{\rho^l}\langle\p_1\varphi\rangle^2 & 0 & 0\\
& -\frac{c^r}{\rho^r}\langle\p_1\varphi\rangle^2 & \frac{c^r}{\rho^r}\langle\p_1\varphi\rangle^2 & 0 & 0 & 0 & 0 & 0 & 0\\
& \langle\p_1\varphi\rangle & \langle\p_1\varphi\rangle & 0 & 0 & -\langle\p_1\varphi\rangle & -\langle\p_1\varphi\rangle & 0 & 0 \\
 & 0 & 0 & -\langle\p_1\varphi\rangle^2 & 0 & 0 & 0 & \langle\p_1\varphi\rangle^2 & 0\\
 & 0 & 0 & -\langle\p_1\varphi\rangle^2 & 0 & 0 & 0 & 0 & 0\\
 & 0 & 0 & 0 & -\langle\p_1\varphi\rangle^2 & 0 & 0 & 0 & \langle\p_1\varphi\rangle^2 \\
 & 0 & 0  & 0 & -\langle\p_1\varphi\rangle^2 & 0 & 0 & 0 & 0  \\
\end{array}
\right)W^\n,}
\end{align*}
and
\begin{align*}
\|\mathbf{M}W^\n|_{x_2=0}-T^\gamma_{\mathbf{M}}W^\n|_{x_2=0}\|_{1,\gamma}\leq\frac{C}{\gamma}\|\mathbf{M}\|_{W^{2,\infty}(\R^2)}\|W^{\n}|_{x_2=0}\|_0\leq\frac{C}{\gamma}\|W^{\n}|_{x_2=0}\|_{0},
\end{align*}
where $C$ is some positive constant. Adding all the estimates above, we have
\begin{align}\label{paraerrorboundary}
\|\mathcal{B}^\gamma(W^\n|_{x_2=0},\psi)-T^\gamma_{b}\psi-T^\gamma_{\mathbf{M}}W^\n|_{x_2=0}\|_{1,\gamma}\leq C\left(\|\psi\|_{1,\gamma}+\frac{1}{\gamma}\|W^{\n}|_{x_2=0}\|_0\right).
\end{align}
Next we consider the interior differential equations, 
\begin{align*}
\vertiii{\gamma A_0^r W^+-T^{\gamma}_{\gamma A_0^r}W^+}_{1,\gamma}^2 & =\int_0^{+\infty}\gamma^2\|A_0^rW^+(\cdot,x_2)-T^{\gamma}_{\gamma A_0^r}W^+(\cdot,x_2)\|^2_{1,\gamma}dx_2\\
& \leq C \|A_0^r\|^2_{W^{2,\infty}(\O)}\vertiii{W^+}^2_0\leq C\vertiii{W^+}^2_0,\\
\vertiii{A_0^r\p_tW^+-T^\gamma_{\ti \delta A_0^r}W^+}_{1,\gamma} & \leq C\vertiii{W^+}_0,\\
\vertiii{A_1^r\p_1W^+-T^\gamma_{\ti \eta A_1^r}W^+}_{1,\gamma} & \leq C\vertiii{W^+}_0,\\
\vertiii{A_0^rC^rW^+-T^\gamma_{ A_0^rC^r}W^+}_{1,\gamma} & \leq C\vertiii{W^+}_0.\\
\end{align*}
Again we have
\begin{align}\label{paraerrorinterior}
\vertiii{\mathcal{L}_{r,l}^\gamma W^{\pm}-T^{\gamma}_{\tau A_0^{r,l}+\ti\eta A_1^{r,l}+ A_0^{r,l}C^{r,l}} W^\pm-I_2\p_2 W^\pm}_{1,\gamma}\leq C\vertiii{W^\pm}_0.
\end{align}

The estimates \eqref{paraerrorboundary} and \eqref{paraerrorinterior} guarantee that if we can estimate the terms on the right-hand side, we can just use the para-linearized system in our  proof below. We will discuss this in details in the following sections. Now we derive the specific expression of the para-linearized system. First we need to eliminate the wave front $\psi$ by observing that 
\begin{align*}
|b(t,x_1,\delta,\eta,\gamma)|^2\geq c(\gamma^2+\delta^2+\eta^2).
\end{align*}
Then, by the G\r{a}rding's inequality \cite[Theorem \textbf{B.7.}]{coulombel2004stability}, we have
\begin{align*}
\Re\langle T^\gamma_{b^*b}\psi,\psi\rangle_{L^2(\R^2)}\geq\frac{c}{2}\|\psi\|^2_{1,\gamma},
\end{align*}
for all $\gamma\geq\gamma_0$, where $\gamma_0$ only depends on $K_0$. By the basic properties of para-differential operators, we have $T^\gamma_{b^*b}=(T^\gamma_{b})^*T^\gamma_{b}+R_{1}^\gamma$ where $R_{1}^{\gamma}$ is an operator of degree $1$, which implies
\begin{align*}
\|\psi\|_{1,\gamma}\leq C\|T^\gamma_{b}\psi\|_0,
\end{align*}
for all $\gamma\geq\gamma_0$. By taking \eqref{paraerrorboundary} into account, we have
\begin{equation}\label{wavefront}
\begin{split}
\|\psi\|_{1,\gamma}& \leq C\left(\|T^\gamma_{b}\psi+T^\gamma_{\mathbf{M}}W^\n|_{x_2=0}\|_{0}+\|T^\gamma_{\mathbf{M}}W^\n|_{x_2=0}\|_{0}\right)\\
& \leq C\left(\frac{1}{\gamma}\|T^\gamma_{b}\psi+T^\gamma_{\mathbf{M}}W^\n|_{x_2=0}\|_{1,\gamma}+\|W^{\n}|_{x_2=0}\|_{0}\right)\\
& \leq C\left(\frac{1}{\gamma}\|\mathcal{B}^\gamma(W^\n,\psi)-T^\gamma_{b}\psi-T^\gamma_{\mathbf{M}}W^\n|_{x_2=0}\|_{1,\gamma} +\frac{1}{\gamma}\|\mathcal{B}^\gamma(W^\n,\psi)\|_{1,\gamma}+\|W^{\n}|_{x_2=0}\|_{0}\right)\\
& \leq C\left(\frac{1}{\gamma}\|\mathcal{B}^\gamma(W^\n,\psi)\|_{1,\gamma}+\|W^{\n}|_{x_2=0}\|_{0}\right).
\end{split}
\end{equation}
Thus we want to identify the part of the boundary condition where the wave front $\psi$ is not involved. Let us consider the following matrix:
\begin{align*}
&\Pi(t,x_1,\delta,\eta,\gamma):=\\
&\begin{pmatrix}
0 & 0 & 1 & 0 & 0 & 0 & 0\\
\tau+\ti v^r \eta & -\ti\eta(v^r-v^l) & 0 & 0 & 0 & 0 & 0\\
-(F_{11}^r-F_{11}^l) & 0 & 0 & v^r-v^l & 0 & 0 & 0\\
-F_{11}^r & 0 & 0 & 0 & v^r-v^l & 0 & 0\\
-(F_{12}^r-F_{12}^l) & 0 & 0 & 0 & 0 & v^r-v^l & 0\\
-F_{12}^r & 0 & 0 & 0 & 0 & 0 & v^r-v^l
\end{pmatrix},
\end{align*}
for $(\tau,\eta)\in\Sigma$. Note that it is homogeneous of degree $0$ with respect to $(\tau,\eta)$. In particular $\Pi b=0$, and 
\begin{align*}
&\Pi\mathbf{M}=\\
&{\Small \left(\begin{array}{*{20}c}
\langle\p_1\varphi\rangle & \langle\p_1\varphi\rangle & 0 & 0 & -\langle\p_1\varphi\rangle & -\langle\p_1\varphi\rangle & 0 & 0\\
-\frac{c^r\langle\p_1\varphi\rangle^2}{\rho^r}(\tau+\ti v^l\eta) & \frac{c^r\langle\p_1\varphi\rangle^2}{\rho^r}(\tau+\ti v^l\eta) & 0 & 0 & \frac{c^l\langle\p_1\varphi\rangle^2}{\rho^l}(\tau+\ti v^r\eta) & -\frac{c^l\langle\p_1\varphi\rangle^2}{\rho^l}(\tau+\ti v^r\eta) & 0 & 0\\
* & * & -a & 0 & * & * & a & 0\\
* & * & -a & 0 & * & * & 0 & 0\\
* & * & 0 & -a & * & * & 0 & a\\
* & * & 0 & -a & * &* & 0 & 0\\
\end{array}
\right),}
\end{align*}
and is homogeneous of degree $0$ with respect to $(\tau,\eta)$, where $a=\langle\p_1\varphi\rangle^2(v^r-v^l)$ and $*$ are some nonzero elements whose exact expressions are not important. By the assumptions in Theorem \ref{variablemain}, we can easily verify that $a\neq0$. Thus the last four rows in $\Pi\mathbf{M}$ are linearly independent, and hence we can expect to estimate $(W_5,W_7,W_{12}, W_{14})^\top|_{x_2=0}$ in terms of $W^\nc|_{x_2=0}:=(W_2, W_3,W_9,W_{10})^\top|_{x_2=0}$. 

We denote the last four rows in $\Pi$ by $\Pi_4$. Then by the exact expression of $\Pi\mathbf{M}$, we have
\begin{align*}
T^\gamma_{\Pi_4\mathbf{M}}W^\n|_{x_2=0}=T^\gamma_{A}W^\nc|_{x_2=0}+T^\gamma_{B}(W_5,W_7,W_{12}, W_{14})^\top|_{x_2=0},
\end{align*}
where $B$ is an invertible matrix in the whole domain and is homogeneous of degree $0$. 
It is easy to check that
\begin{align*}
|B(t,x_1)|^2\geq c,
\end{align*}
for some constant $c>0$. By the G\r{a}rding's inequality, we have
\begin{align*}
\Re\langle T^\gamma_{B^*B} (W_5,W_7,W_{12}, W_{14})^\top,(W_5,W_7,W_{12}, W_{14})^\top\rangle|_{x_2=0} \geq \frac{c}{2} \|(W_5,W_7,W_{12}, W_{14})^\top|_{x_2=0}\|_0^2,
\end{align*}
for all $\gamma\geq\gamma_0$. This implies
\begin{align*}
 \|(W_5,W_7,W_{12}, W_{14})^\top|_{x_2=0}\|_0&\leq C\|T^\gamma_B (W_5,W_7,W_{12}, W_{14})^\top|_{x_2=0}\|_0\\
 &\leq C(\|T^\gamma_{\Pi_4\mathbf{M}}W^\n|_{x_2=0}\|_0+\|T^\gamma_{A}W^\nc|_{x_2=0}\|_0).
\end{align*}
Similarly,  for the wave front $\psi$, we have
\begin{align*}
&\|T^\gamma_{\Pi_4\mathbf{M}}W^\n|_{x_2=0}\|_0 \leq \|T^\gamma_{\Pi_4}\mathcal{B}^\gamma(W^\n,\psi)-T^\gamma_{\Pi_4}T^\gamma_{\mathbf{M}}W^\n|_{x_2=0}-T^\gamma_{\Pi_4}T^\gamma_{b}\psi\|_0\\
&\qquad\qquad\qquad\qquad\qquad\qquad\qquad+\|T^\gamma_{\Pi_4}\mathcal{B}^\gamma(W^\n,\psi)\|_0+\|W^\n|_{x_2=0}\|_{-1,\gamma}+\|\psi\|_0\\
&\leq\frac{1}{\gamma} \|\mathcal{B}^\gamma(W^\n,\psi)-T^\gamma_{\mathbf{M}}W^\n|_{x_2=0}-T^\gamma_{b}\psi\|_{1,\gamma}+\frac{1}{\gamma}\|\mathcal{B}^\gamma(W^\n,\psi)\|_{1,\gamma}+\|W^\n|_{x_2=0}\|_{-1,\gamma}+\|\psi\|_0\\
&\leq \frac{1}{\gamma}\|\psi\|_{1,\gamma}+\frac{1}{\gamma}\|\mathcal{B}^\gamma(W^\n,\psi)\|_{1,\gamma}+\frac{1}{\gamma}\|W^\n|_{x_2=0}\|_0.
\end{align*}
By taking $\gamma$ large enough, we have
\begin{align}\label{normalelasticity}
 \|(W_5,W_7,W_{12}, W_{14})^\top|_{x_2=0}\|_0\leq C\left( \frac{1}{\gamma}\|\psi\|_{1,\gamma}+\frac{1}{\gamma}\|\mathcal{B}^\gamma(W^\n,\psi)\|_{1,\gamma}+\|W^\nc|_{x_2=0}\|_0\right).
\end{align}
Combining \eqref{wavefront} and \eqref{normalelasticity} together, we have
\begin{align*}
\|(W_5,W_7,W_{12}, W_{14})^\top|_{x_2=0}\|_0+\|\psi\|_{1,\gamma}\leq C(\frac{1}{\gamma}\|\mathcal{B}^\gamma(W^\n,\psi)\|_{1,\gamma}+\|W^\nc|_{x_2=0}\|_0),
\end{align*}
which suggests that we can obtain the estimate of $(W_5,W_7,W_{12}, W_{14})^\top|_{x_2=0}$ and $\psi$ through the source terms and $W^\nc|_{x_2=0}$. To obtain the estimate of $W^\nc|_{x_2=0}$, we need to utilize  the other part of the boundary conditions which we denote by
\begin{align*}
T^\gamma_{\beta}W^{\nc}|_{x_2=0}=\tilde{G},
\end{align*}
where 
\begin{align*}
\beta=\begin{pmatrix}
\langle\p_1\varphi\rangle & \langle\p_1\varphi\rangle & -\langle\p_1\varphi\rangle & -\langle\p_1\varphi\rangle\\
-\frac{c^r\langle\p_1\varphi\rangle^2}{\rho^r}(\tau+\ti v^l\eta) & \frac{c^r\langle\p_1\varphi\rangle^2}{\rho^r}(\tau+\ti v^l\eta) & \frac{c^l\langle\p_1\varphi\rangle^2}{\rho^l}(\tau+\ti v^r\eta) & -\frac{c^l\langle\p_1\varphi\rangle^2}{\rho^l}(\tau+\ti v^r\eta)
\end{pmatrix},
\end{align*}
and is homogeneous of degree $0$ with respect to $(\tau,\eta)$. The above matrix $\beta$ is from the first two rows of $\Pi\mathbf{M}$ and in the symbol class $\Gamma_2^0$. 

Finally, with the help of the above estimates, we can define the para-linearized system as
\begin{equation}\label{paralinearizationsystem}
\begin{cases}
T^\gamma_{\tau A_0^r+\ti\eta A_1^r+A_0^rC^r}W^++I_2\p_2W^+=\tilde{F}^+,\\
T^\gamma_{\tau A_0^l+\ti\eta A_1^l+A_0^lC^l}W^-+I_2\p_2W^-=\tilde{F}^-,\\
T^\gamma_{\beta}W^{\nc}|_{x_2=0}=\tilde{G}.
\end{cases}
\end{equation}
Moreover, our aim is to derive the following estimate:
\begin{align}\label{noncharacteristicpart}
\|W^{\nc}|_{x_2=0}\|_0^2\leq C_0\left(\frac{1}{\gamma^3}\vertiii{\tilde{F}}_{1,\gamma}^2+\frac{1}{\gamma^2}\|\tilde{G}\|_{1,\gamma}^2\right),
\end{align}
where $\tilde{F}=(\tilde{F}^+,\tilde{F}^-)^\top$. In fact,  \eqref{wavefront}, \eqref{normalelasticity} and \eqref{noncharacteristicpart}  together imply Theorem \ref{variablemain}, which can be shown as follows.

First, by the above discussion, we have
\begin{align*}
&\vertiii{\mathcal{L}_{r,l}^\gamma W^{\pm}-T^{\gamma}_{\tau A_0^{r,l}+\ti\eta A_1^{r,l}+ A_0^{r,l}C^{r,l}} W^\pm-I_2\p_2 W^\pm}_{1,\gamma}\leq C\vertiii{W^\pm}_0,\\
&\|\mathcal{B}^\gamma(W^\n,\psi)-T^\gamma_{b}\psi-T^\gamma_{\mathbf{M}}W^\n|_{x_2=0}\|_{1,\gamma}\leq C\left(\|\psi\|_{1,\gamma}+\frac{1}{\gamma}\|W^{\n}|_{x_2=0}\|_0\right).
\end{align*}
Then
\begin{align*} 
&\|\tilde{G}\|_{1,\gamma}=\|T^\gamma_{\beta}W^{\nc}|_{x_2=0}\|_{1,\gamma}=\|T^\gamma_{\Pi_2\mathbf{M}}W^\nc|_{x_2=0}\|_{1,\gamma}=\|T^\gamma_{\Pi_2b}\psi+T^\gamma_{\Pi_2\mathbf{M}}W^\nc|_{x_2=0}\|_{1,\gamma},
\end{align*}
where $\Pi_2$ is the first two rows of $\Pi$. Hence we have
\begin{align*}
\|\tilde{G}\|_{1,\gamma} & \leq\|T^\gamma_{\Pi_2}(T^\gamma_{b}\psi+T^\gamma_{\mathbf{M}}W^\nc|_{x_2=0})\|_{1,\gamma}+\|\psi\|_{1,\gamma}+\|W^\nc|_{x_2=0}\|_0\\
&\leq\|\mathcal{B}^\gamma(W^\n,\psi)-T^\gamma_{b}\psi-T^\gamma_{\mathbf{M}}W^\n|_{x_2=0}\|_{1,\gamma}+\|\mathcal{B}^\gamma(W^\n,\psi)\|_{1,\gamma}+\|\psi\|_{1,\gamma}+\|W^\nc|_{x_2=0}\|_0\\
&\leq \|\psi\|_{1,\gamma}+\frac{1}{\gamma}\|W^\n|_{x_2=0}\|_0+\|W^\nc|_{x_2=0}\|_0+\|\mathcal{B}^\gamma(W^\n,\psi)\|_{1,\gamma}.
\end{align*}
Moreover,
\begin{align*}
&\vertiii{\tilde{F}^\pm}_{1,\gamma}=\vertiii{\mathcal{L}_{r,l}^\gamma W^{\pm}}_{1,\gamma}+\vertiii{\tilde{F}^\pm-\mathcal{L}_{r,l}^\gamma W^{\pm}}_{1,\gamma}\leq \vertiii{\mathcal{L}_{r,l}^\gamma W^{\pm}}_{1,\gamma} + C\vertiii{W^\pm}_0.
\end{align*}
Thus \eqref{noncharacteristicpart} implies
\begin{align*}
\|W^{\nc}|_{x_2=0}\|_0^2\leq C_0 & \left(\frac{1}{\gamma^3}\vertiii{\mathcal{L}_{r}^\gamma W^{+}}_{1,\gamma}^2+\frac{1}{\gamma^3}\vertiii{W^+}_0^2+\frac{1}{\gamma^3}\vertiii{\mathcal{L}_{l}^\gamma W^{-}}_{1,\gamma}^2+\frac{1}{\gamma^3}\vertiii{W^-}_0^2\right.\\
&\left.+\frac{1}{\gamma^2}\|\mathcal{B}^\gamma(W^\n,\psi)\|_{1,\gamma}^2 +\frac{1}{\gamma^4} \|W^\n|_{x_2=0}\|_0^2+\frac{1}{\gamma^2}\|\psi\|_{1,\gamma}^2\right).
\end{align*}
Combining the above with \eqref{wavefront} and \eqref{normalelasticity}, we have
\begin{align*}
&\|W^\n|_{x_2=0}\|_0^2+\|\psi\|_{1,\gamma}^2\leq C  \left(\frac{1}{\gamma^3}\vertiii{\mathcal{L}_{r}^\gamma W^{+}}_{1,\gamma}^2+\frac{1}{\gamma^3}\vertiii{\mathcal{L}_{l}^\gamma W^{-}}_{1,\gamma}^2+\frac{1}{\gamma^2}\|\mathcal{B}^\gamma(W^\n,\psi)\|_{1,\gamma}^2 \right.\\
&\hspace{1.6in} \left.+\frac{1}{\gamma^3}\vertiii{W^+}_0^2+\frac{1}{\gamma^3}\vertiii{W^-}_0^2+\frac{1}{\gamma^4} \|W^\n|_{x_2=0}\|_0^2+\frac{1}{\gamma^2}\|\psi\|_{1,\gamma}^2\right)\\
&\ \ \leq C\left(\frac{1}{\gamma^3}\vertiii{\mathcal{L}_{r}^\gamma W^{+}}_{1,\gamma}^2+\frac{1}{\gamma^3}\vertiii{\mathcal{L}_{l}^\gamma W^{-}}_{1,\gamma}^2+\frac{1}{\gamma^2}\|\mathcal{B}^\gamma(W^\n,\psi)\|_{1,\gamma}^2+\frac{1}{\gamma^3}\vertiii{W^+}_0^2+\frac{1}{\gamma^3}\vertiii{W^-}_0^2\right).
\end{align*}
Then by Lemma \ref{L2ofdotV}, we obtain the Theorem \ref{variablemain}. 

Therefore the only thing left to prove is the estimate \eqref{noncharacteristicpart} from the para-linearized system \eqref{paralinearizationsystem}.

\section{Microlocalization}\label{microlocal}

From the previous section, we obtain a system of ODE in \eqref{paralinearizationsystem}. We note that the para-linearization is an analogue of the Fourier transform in the constant coefficient case. To derive the estimate \eqref{noncharacteristicpart}, we need to investigate at each frequency point $(\tau,\eta)\in\Pi$ if the boundary conditions in \eqref{paralinearizationsystem} provide enough information on the incoming mode of the ODE. As in the constant coefficient case, we want to estimate the rate of vanishing for the Lopatinskii determinant near its roots. 

In this section, we  analyze and classify the situations which we will need to deal with in the system \eqref{paralinearizationsystem} for each frequency point $(\tau,\eta)\in\Pi$. Specifically, we need to identify those frequencies where the standard Kreiss symmetrizer method can not be applied as in \cite{coulombel2002weak, coulombel2004stability}. More precisely,  these frequency points are points where the system \eqref{paralinearizationsystem} can not be reduced to that  involving only the non-characteristic part of the unknown $W^\nc$ and where the Lopatinskii determinant vanishes. 
For simplicity in the argument, we only focus on the case on the unit hemisphere $\Sigma=\{(\tau,\eta):\;|\tau|^2+\eta^2=1\text{ and } \Re\tau\geq 0\}$ in the frequency space $\Pi$. Then we can extend the results to the whole frequency space $\Pi$ by their appropriate homogeneities.

\subsection{Poles} In this part, we identify the frequencies where the system \eqref{paralinearizationsystem} can not be reduced to a  problem  involving only the non-characteristic part of the unknown $W^\nc$. As in the constant coefficient case \cite{coulombel2004stability}, we call these points the poles of the system. For this purpose, we just formally consider the first order symbols $\tau A_0^r+\ti\eta A_1^r$ and $\tau A_0^l+\ti\eta A_1^l$ and the boundary symbol $\beta$ in \eqref{paralinearizationsystem}. In particular, we focus on the following differential system:
\begin{equation}\label{formaldifferentialequation}
\begin{cases}
&(\tau A_0^r+\ti\eta A_1^r)W^++I_2\p_2W^+=0,\\
&(\tau A_0^l+\ti\eta A_1^l)W^-+I_2\p_2W^-=0,\\
&\beta W^{\nc}|_{x_2=0}=0,
\end{cases}
\end{equation}
where $A_0^{r,l}$, $A_1^{r,l}$, $I_2$ and $\beta$ are the matrices we specified in the previous section. We note that $I_2$ is a diagonal matrix with only two elements in diagonal being $1$ and others being $0$. So there are two differential equations and five algebraic equations for $W^+$ and $W^-$ respectively. This suggests that our system is characteristic. However, for most of the points $(\tau,\eta)$ in the frequency space, we can use the algebraic equations to reduce the system to a system only involving differential equations, which is a non-characteristic system.  To find out those points, we consider the algebraic equations for $W^+$: 
\begin{align*}
&\begin{pmatrix}
\tau+\ti v^r \eta & \frac{\ti\eta {c^r}^2}{\langle\p_1\Phi^r\rangle\rho^r} & \frac{\ti\eta {c^r}^2}{\langle\p_1\Phi^r\rangle\rho^r} & -\ti\eta F_{11}^r & 0 & -\ti\eta F_{12}^r & 0\\
-\ti F_{11}^r\eta & 0 & 0 & \tau+\ti v^r\eta & 0 & 0 & 0\\
0 & -\frac{\ti c^r F_{11}^r\eta}{\rho^r} & \frac{\ti c^r F_{11}^r\eta}{\rho^r} & 0 & \tau+\ti v^r\eta & 0 & 0\\
-\ti F_{12}^r\eta & 0 & 0 & 0 & 0 & \tau+\ti v^r\eta & 0\\
0 & -\frac{\ti c^r F_{12}^r\eta}{\rho^r} & \frac{\ti c^r F_{12}^r\eta}{\rho^r} & 0 & 0 & 0 & \tau+\ti v^r\eta
\end{pmatrix}W^+=0.
\end{align*}
The algebraic equation for $W^-$ can be treated similarly.
It is obvious that, as in the constant coefficient case, if $(\tau+\ti v^r \eta)((\tau+\ti v^r \eta)^2+\left({F_{11}^r}^2+{F_{12}^r}^2\right)$ $\eta^2)\neq 0$, we can uniquely determine $W_1$, $W_4$, $W_5$, $W_6$, $W_7$ from $W_2$, $W_3$. Then by using the differential equations in \eqref{formaldifferentialequation}, we can obtain the differential equations only involving $W_2$ and $W_3$, which are 
\begin{align*}
\mathbb{A}^r\begin{pmatrix}
W_2\\
W_3
\end{pmatrix}+\p_2\begin{pmatrix}
W_2\\
W_3
\end{pmatrix}=0,
\end{align*} 
where
\begin{align*}
\mathbb{A}^r:=\begin{pmatrix}
\mu^r & -m^r\\
m^r & -\mu^r
\end{pmatrix}+\ti \eta \frac{\p_1\Phi^r\p_2\Phi^r}{\langle\p_1\Phi^r\rangle^2}\begin{pmatrix}
1 & 0\\
0 & 1
\end{pmatrix},
\end{align*}
and
\begin{align*}
&\mu^r=-\frac{\p_2\Phi^r(\tau+\ti v^r\eta)}{c^r\langle\p_1\Phi^r\rangle}-\frac{\p_2\Phi^r({F_{11}^r}^2+{F_{12}^r}^2)\eta^2}{2\langle\p_1\Phi^r\rangle c^r(\tau+\ti v^r\eta)}-\frac{\p_2\Phi^r c^r(\tau+\ti v^r\eta)\eta^2}{2\langle\p_1\Phi^r\rangle^3\left((\tau+\ti v^r\eta)^2+({F_{11}^r}^2+{F_{12}^r}^2)\eta^2\right)},\\
&m^r=-\frac{\p_2\Phi^r({F_{11}^r}^2+{F_{12}^r}^2)\eta^2}{2\langle\p_1\Phi^r\rangle c^r(\tau+\ti v^r\eta)}+\frac{\p_2\Phi^r c^r(\tau+\ti v^r\eta)\eta^2}{2\langle\p_1\Phi^r\rangle^3\left((\tau+\ti v^r\eta)^2+({F_{11}^r}^2+{F_{12}^r}^2)\eta^2\right)},\\
\end{align*}
Combining the above discussions on the equations for $W^-$, we see that the points in the frequency space where we can not reduce the system to the non-characteristic form  are the points where 
\begin{equation}\label{poledefn}
(\tau+\ti v^{r,l} \eta)\left((\tau+\ti v^{r,l}\eta)^2+\left({F_{11}^{r,l}}^2+{F_{12}^{r,l}}^2\right)\eta^2\right)= 0. 
\end{equation}
We call this kind of points the poles of the system. We define the sets of poles as 
\begin{align*}
\Upsilon_p:=\left\{(t,x_1,x_2,\tau,\eta)\in \R^3_+\times\Sigma: (\tau+\ti v^{r,l} \eta)\left((\tau+\ti v^{r,l} \eta)^2+\left({F_{11}^{r,l}}^2+{F_{12}^{r,l}}^2\right)\eta^2\right)= 0\right\}.
\end{align*}
Here for each fixed $(t,x_1,x_2)\in\O$, there are six kinds of frequencies in $\Upsilon_p$. Each one corresponds to the root of one factor in the above definition.
For the root of each factor in the definition, fixing $(t,x_1)\in\R^2$, we can view it as a curve in  the frequency space parametrized by $x_2$   originating from the boundary $x_2=0$ and propagate into the interior of the domain $x_2>0$.

\subsection{Roots of the Lopatinskii determinant}
Now we identify the frequencies where the Lopatinskii determinant vanishes. First we need to derive the Lopatinskii determinant. This requires us to find out the eigenvectors of $\mathbb{A}^{r,l}$ corresponding to the eigenvalues with negative real parts. By a direct computation,  we denote the eigenvalue of $\mathbb{A}^r$ with a negative real part by $\o^r+\ti\frac{\p_1\Phi^r\p_2\Phi^r}{\langle\p_1\Phi^r\rangle^2}$ where
\begin{align*}
({\o^r})^2=({\mu^r})^2-({m^r})^2=\frac{({\p_2\Phi^r})^2}{({c^r})^2\langle\p_1\Phi^r\rangle^4}\left[\langle\p_1\Phi^r\rangle^2\left((\tau+\ti v^r\eta)^2+({F_{11}^r}^2+{F_{12}^r}^2)\eta^2\right)+({c^r})^2\eta^2\right].
\end{align*}
The corresponding eigenvector is 
\begin{align}
E^r=\begin{pmatrix}
-\alpha^r(\mu^r+\o^r)\\
-\alpha^rm^r
\end{pmatrix},
\end{align}
where $\alpha^r=(\tau+\ti v^r\eta)\left((\tau+\ti v^r\eta)^2+({F_{11}^r}^2+{F_{12}^r}^2)\eta^2\right)$. It is easy to check that $E^r$ is homogeneous of degree $0$ with respect to $(\tau,\eta)$. 

The case for $W^-$ is similar. The eigenvalue with a negative real part is $\o^l+\ti\frac{\p_1\Phi^l\p_2\Phi^l}{\langle\p_1\Phi^l\rangle^2}$, where
\begin{align*}
({\o^l})^2=({\mu^l})^2-({m^l})^2=\frac{({\p_2\Phi^l})^2}{({c^l})^2\langle\p_1\Phi^l\rangle^4}\left[\langle\p_1\Phi^l\rangle^2\left((\tau+\ti v^l\eta)^2+({F_{11}^l}^2+{F_{12}^l}^2)\eta^2\right)+({c^l})^2\eta^2\right].
\end{align*}
The corresponding eigenvector is 
\begin{align}
E^l=\begin{pmatrix}
\alpha^lm^l\\
\alpha^l(\mu^l-\o^l)
\end{pmatrix},
\end{align}
where $\alpha^l=(\tau+\ti v^l\eta)\left((\tau+\ti v^l\eta)^2+({F_{11}^l}^2+{F_{12}^l}^2)\eta^2\right)$. Similarly $E^l$ is homogeneous of degree $0$ with respect to $(\tau,\eta)$, 

Note that by a similar computation in the constant coefficient cases, the above eigenvalues and eigenvectors are all well-defined and smooth on the whole space $\R^3\times\Pi$. Hence the Lopantiskii determinant is well-defined for all the points in the frequency space, which is 
\begin{equation}\label{variableLopatinskii}
{\small 
\begin{split}
\left.\text{det}\left(\beta\begin{pmatrix}
E^r & 0\\
0 & E^l
\end{pmatrix}\right)\right|_{x_2=0}&=\frac{c^4a_1^2}{\rho}(\tau+\ti v^r\eta)(\tau +\ti v^l \eta)\left(\frac{a_1^4}{a_2^ra_2^l}\o^r\o^l+\eta^2\right)\left(\frac{\o^r}{a_2^r}-\frac{\o^l}{a_2^l}\right)\\
&\times\left(\frac{a_2^r}{a_1 c}\left((\tau+\ti v^r\eta)^2+({F_{11}^r}^2+{F_{12}^r}^2)\eta^2\right)-(\tau+\ti v^r \eta)\o^r\right)\\
&\times\left(\frac{a_2^l}{a_1 c}\left((\tau+\ti v^l\eta)^2+({F_{11}^l}^2+{F_{12}^l}^2)\eta^2\right)+(\tau+\ti v^l \eta)\o^l\right),
\end{split}}
\end{equation}
and is homogeneous of degree $0$ with respect to $(\tau,\eta)$, where $a_1=\langle\p_1\varphi\rangle$, $a_2^{r,l}=\p_2\Phi^{r,l}|_{x_2=0}$, $c=c^r|_{x_2=0}=c^l|_{x_2=0}$. The last equality for $c$ is from the fact that $c^{r,l}$ only depends on the density $\rho^{r,l}$ which is continuous at $x_2=0$ by \eqref{variableboundary}. 

By analyzing each factor in the Lopatinskii determinant we see that the last two factors in \eqref{variableLopatinskii} are never equal to zeros, and the first two factors corresponds to the roots $\tau=-\ti v^{r,l} \eta$ respectively. Thus we only need to discuss the third and forth factors.

We note that all the coefficients in the factors of the Lopatinskii determinant are coninuous with respect to the background state $U|_{x_2=0}:=(U^r|_{x_2=0},U^l|_{x_2=0})$ and $\Phi|_{x_2=0}:=(\Phi^r|_{x_2=0},\Phi^l|_{x_2=0})$, and those factors reduce to the corresponding factors in the constant coefficient case, if the perturbation in \eqref{variablestates} is zero. By assuming $K$  sufficient small and using the continuity argument, we conclude that the number of the roots in the third and forth factors in \eqref{variableLopatinskii} is the same as the number of roots in the corresponding factors in the constant coefficient case.
Hence there are two roots of $\frac{a_1^4}{a_2^ra_2^l}\o^r\o^l+\eta^2$ which we denote by $\tau=\ti V_1\eta$ and $\tau=\ti V_2 \eta$, and there is one root of $\frac{\o^r}{a_2^r}-\frac{\o^l}{a_2^l}$  denoted by $\tau=\ti V_3\eta$, if $\bar{v}^2>2c(\bar{\rho})^2+\bar{F}_{11}^2+\bar{F}_{12}^2$. 

We can verify that $V_1(U|_{x_2=0},\nabla\Phi|_{x_2=0})$, $V_2(U|_{x_2=0},\nabla\Phi|_{x_2=0})$ and $V_3(U|_{x_2=0},\nabla\Phi|_{x_2=0})$ are all real and continuously depend on the background state $U|_{x_2=0}$ and $\nabla\Phi|_{x_2=0}$. Furthermore, all the roots above are simple if they do not coincide. So the set of roots of the Lopatinskii determinant can be represented by the following set on the boundary of the domain:
\begin{align} \label{rootsofL}
\Upsilon_r^0:=\{(t,x_1,\tau,\eta)\in \R^2\times \Sigma : \Re\tau=0 \text{ and } \sigma=0 \},
\end{align}
where 
\begin{equation*}
\sigma = \left\{\begin{array}{ll}
(\delta+v^r|_{x_2=0}\eta)(\delta+v^l|_{x_2=0}\eta)(\delta-V_1\eta)(\delta-V_2\eta)(\delta-V_3\eta), & \text{if } \bar{v}^2>2c(\bar{\rho})^2+\bar{F}_{11}^2+\bar{F}_{12}^2, \\
(\delta+v^r|_{x_2=0}\eta)(\delta+v^l|_{x_2=0}\eta)(\delta-V_1\eta)(\delta-V_2\eta), & \text{if } \bar{v}^2\leq \bar{F}_{11}^2+\bar{F}_{12}^2
\end{array}\right.
\end{equation*}
on $\Sigma$. Here the set $\Upsilon_r^0$ and function $\sigma$ can be naturally extended into the interior of physical domain where $x_2>0$. More precisely, we can define the coefficients in $\sigma$ for $x_2>0$ by the continuous dependence of $V_1$, $V_2$ and $V_3$ on the background state $U$ and $\nabla\Phi$. We denote the extended set by 
\begin{align*}
\Upsilon_r:=\{(t,x_1,x_2,\tau,\eta)\in \R^3_+\times\Sigma: \Re\tau=0 \text{ and } \sigma=0 \}.
\end{align*}
Just like $\Upsilon_p$, for fixed $(t,x_1)\in\R^2$, the root of each factor of $\sigma$ can be considered as a curve in the frequency space $\Sigma$ parametrized by $x_2$, which originates from the boundary $x_2=0$ and propagates into the interior of the domain $x_2>0$. 

\begin{Remark}\label{rk_evalue}
In fact, there is one more type of points we need to pay attention to due to the restraint of para-differential calculus. Those are the points where the eigenvalues $\o^r$ and $\o^l$ vanish. When we deal with these points in deriving the energy estimates, we can not perform the upper triangularization method for poles and roots of the Lopatinskii determinant, because the upper triangularization method will introduce $\o^r$ and $\o^l$ as the symbols but $\o^r$ and $\o^l$ do not satisfy the requirements (\cite[Definition \textbf{B.1.}]{coulombel2004stability}) of the symbol class near the points where $\o^r=0$ and $\o^l=0$ respectively. Hence we require that the roots of the eigenvalues $\o^r$ and $\o^l$ do not coincide with poles and roots of the Lopatinskii determinant. Then we can perform the classical Kreiss symmetrizer method to deal with the zeros of the eigenvalues, which is the same way as in dealing with the points other than the poles and roots of the Lopatinskii determinant. 
\end{Remark}

From the above remark and the definition of $\Upsilon_p$, $\Upsilon_r$, it is important to consider the following frequency sets:
\begin{align*}
\text{(1)}\quad&\Upsilon_p^{(1)} = \Upsilon_r^{(1)} := \{(t, x_1, x_2, \tau, \eta):\ \tau = -\ti v^{r,l}\eta \},\\
\text{(2)}\quad &\Upsilon_r^{(2)} : = \left\{(t, x_1, x_2, \tau, \eta):\ \begin{array}{l} \tau = \ti V_1\eta,\quad \tau = \ti V_2\eta \text{ and }\\
 \tau = \ti V_3\eta \quad \text{if } \bar{v}^2>2c(\bar{\rho})^2+\bar{F}_{11}^2+\bar{F}_{12}^2 \end{array}\right\},\\
\text{(3)}\quad&\Upsilon_p^{(2)} := \left\{(t, x_1, x_2, \tau, \eta):\ \tau = - \ti \left(v^{r,l}\pm\sqrt{{F_{11}^{r,l}}^2+{F_{12}^{r,l}}^2}\right)\eta \right\},
\end{align*}
where (1) corresponds to the zeros to the first factor in the definition of $\Upsilon_p$, which are also the zeros of the first two factors of $\sigma$ in the definition of $\Upsilon_r$; (2)  corresponds to the zeros to the last three or two factors of $\sigma$ in $\Upsilon_r$; and (3) corresponds to the zeros to the second factor in $\Upsilon_p$. Notice that $\Upsilon_{p} = \Upsilon_p^{(1)} \cup \Upsilon_p^{(2)}$ and $\Upsilon_{r} = \Upsilon_r^{(1)} \cup \Upsilon_r^{(2)}$.

From Remark \ref{rk_evalue}  we need to require that $\Upsilon_\o : = \{ \o^{r,l} = 0 \}$ does not intersect with $\Upsilon_p$ and $\Upsilon_r$, i.e. $\Upsilon_\o \cap \left( \Upsilon_p \cup \Upsilon_r \right) = \emptyset$.  Note that in general any two curves in the classes $\Upsilon_p$ and $\Upsilon_r$ may also intersect with each other, which causes extreme difficulties in deriving the energy estimates through the symbolic cut-off functions in the frequency space, unless the two intersecting curves are identical (that corresponds to the case where the curves are in $\Upsilon_p^{(1)} = \Upsilon_r^{(1)}$). Hence we need further conditions on the variable background states \eqref{variablestates} which guarantee that any two curves in $\Upsilon_p$ and $\Upsilon_r$ do not intersect with each other in the whole domain $\R^3_+\times\Sigma$, unless they are identical. By the continuous dependence of these three types of curves with respect to the background states, it suffices to assume that all the poles, the roots of the Lopatinskii determinant and the zeros of eigenvalues $\o^{r,l}$ for the constant background states $\bar{U}$ in \eqref{variablestates} do not coincide with each other,  and $K$ in \eqref{variablesmallness} is small enough.

Tracking carefully what the frequency sets $\Upsilon_p, \Upsilon_r$ and $\Upsilon_\o$ are at the constant background states $\bar U$ (see \cite[Section 5]{chen2015linear}), we know that in the supersonic regime $\bar{v}^2> 2c(\bar{\rho})^2+\bar{F}_{11}^2+\bar{F}_{12}^2$, any two curves from these three sets do not intersect. In the subsonic regime $\bar{v}^2< \bar{F}_{11}^2+\bar{F}_{12}^2$, on the other hand,  we have 
\begin{align*}
& \Upsilon_p^{(2)} \cap \Upsilon_r^{(2)} = \Upsilon_r^{(2)} \cap \Upsilon_\o = \emptyset,\\
& \bar{v}^2\neq\frac{\bar{F}_{11}^2+\bar{F}_{12}^2}{4} \quad \Rightarrow \quad \Upsilon_p^{(1)} \cap \Upsilon_p^{(2)} = \Upsilon_r^{(1)} \cap \Upsilon_p^{(2)} = \emptyset,\\
& \bar{v}^2\neq\frac{\left(\sqrt{\bar{F}_{11}^2+\bar{F}_{12}^2+c(\bar{\rho})^2}-\sqrt{\bar{F}_{11}^2+\bar{F}_{12}^2}\right)^2}{4} \quad \Rightarrow \quad \Upsilon_p^{(2)} \cap \Upsilon_\o = \emptyset,\\
& \bar{v}^2\neq\frac{\bar{F}_{11}^2+\bar{F}_{12}^2+c(\bar{\rho})^2}{4} \quad \Rightarrow \quad \Upsilon_p^{(1)} \cap \Upsilon_\o = \Upsilon_r^{(1)} \cap \Upsilon_\o = \emptyset,  \\
& \bar{v}^2\neq\frac{\left(\bar{F}_{11}^2+\bar{F}_{12}^2\right)\left(2c(\bar{\rho})^2+\bar{F}_{11}^2+\bar{F}_{12}^2\right)}{4\left(\bar{F}_{11}^2+\bar{F}_{12}^2+c(\bar{\rho})^2\right)} \quad \Rightarrow \quad \Upsilon_p^{(1)} \cap \Upsilon_r^{(2)} = \Upsilon_r^{(1)} \cap \Upsilon_r^{(2)} = \emptyset.
\end{align*}
This explains why we need to assume the extra conditions on the background velocity in Theorem \ref{variablemain} (ii).

\section{Estimates in Each Case}\label{Estimateineachcase}

Now we want to derive the estimates for each case and eventually obtain the desired estimates from the para-linearized system.  From the discussions in the previous section, we have three isolated cases where we can not apply the Kreiss symmetrizer method. Each relation of $\tau$ and $\eta$ in these cases corresponds to a curve on $\Sigma$ with fixed $(t,x_1)$. In the previous section, we explained that we only focus on the situation where all these curves do not intersect with each other, and hence for each curve, we can construct an open neighborhood around it. Up to shrinking those open neighborhoods, we can guarantee that those neighborhoods do not intersect with each other and  do not contain any point in $\Upsilon_\o$.  We denote, on $\R^3\times\Sigma$, the open neighborhood around $\tau=-\ti v^r\eta$ by $\V_{p_1}^r$, around $\tau=-\ti v^l\eta$ by $\V_{p_1}^l$, around $\tau=-\ti \left(v^{r}+\sqrt{{F_{11}^{r}}^2+{F_{12}^{r}}^2}\right)\eta$ by $\V_{p_2}^1$, around $\tau=-\ti \left(v^{r}-\sqrt{{F_{11}^{r}}^2+{F_{12}^{r}}^2}\right)\eta$ by $\V_{p_2}^2$, around $\tau=-\ti \left(v^{l}+\sqrt{{F_{11}^{l}}^2+{F_{12}^{l}}^2}\right)\eta$ by $\V_{p_2}^3$, around $\tau=-\ti \left(v^{l}-\sqrt{{F_{11}^{l}}^2+{F_{12}^{l}}^2}\right)\eta$ by $\V_{p_2}^4$, around $\tau=\ti V_1\eta$ by $\V_{r}^1$, around $\tau=\ti V_2\eta$ by $\V_{r}^2$ and around $\tau=\ti V_3\eta$ by $\V_{r}^3$. 

In the following argument we always consider $\tau=\ti V_3\eta$ being a root of the Lopatinskii determinant. For the background state \eqref{variablestates} satisfying $\bar{v}^2\leq \bar{F}_{11}^2+\bar{F}_{12}^2$, we can just drop the part of argument about $\tau=\ti V_3\eta$.

\subsection{Case 1 -- points in $\Upsilon_p^{(1)} = \Upsilon_r^{(1)}$}\label{subsec_poles and roots}
In this case, we consider the kind of frequencies which are both poles and  roots of the Lopatinskii determinant. More precisely, we consider the frequencies in $\V_{p_1}^r$ and $\V_{p_1}^l$. Without loss of generality, we just take $\V_{p_1}^r$ as an example. The other cases can be treated exactly the same way.  

We note that $\V_{p_1}^r$ only contains the poles of the equations for $W^+$ in \eqref{paralinearizationsystem}, but does not contain the poles of the equations for $W^-$ in \eqref{paralinearizationsystem}. So the estimates of $W^-$ can be obtained in the exactly same way as the estimates of $W^+$ in the Case 2. In this part we only provide details on how to derive the estimates for $W^+$. 

To isolate the objective frequencies, we introduce the smooth cut-off function $\chi_{p_1}$ whose range is the closed interval $[0,1]$. In particular, on $\R^3_+\times\Sigma$, the support of $\chi_{p_1}$ is contained in $\V_{p_1}^r$ and equals $1$ on a smaller neighborhood of the curve satisfying $\tau=-\ti v^r\eta$. Then we extend $\chi_{p_1}$ by homogeneity of degree $0$ with respect to $(\tau,\eta)$ into the whole space $\R^3_+\times\Pi$. As in \cite{coulombel2004stability}, $\chi_{p_1}$ is in the symbol class $\Gamma_k^0$ for any integer $k$. With this cut-off function, we define
\begin{align*}
W^+_{p_1}:=T^\gamma_{\chi_{p_1}}W^+.
\end{align*}
From \eqref{paralinearizationsystem}, we can obtain 
\begin{align*}
I_2\p_2W^+_{p_1}=I_2T^\gamma_{\p_2\chi_{p_1}}W^++T^{\gamma}_{\chi_{p_1}}F^+-T^\gamma_{\chi_{p_1}}T^\gamma_{\tau A_0^r+\ti\eta A_1^r+A_0^rC^r}W^+.
\end{align*}
By para-differential calculus, we have
\begin{align}\label{para1}
T^\gamma_{\mathcal{A}^r}W^+_{p_1}+T^\gamma_{A_0^rC^r}W^+_{p_1}+T^\gamma_rW^++I_2\p_2W^+_{p_1}=T^\gamma_{\chi_{p_1}}F+R_{-1}W^+,
\end{align}
where   $\mathcal{A}^r:=\tau A_0^r+\ti\eta A_1^r$ and $r$ is a symbol in the class $\Gamma_1^0$ whose support is in the area with $\chi_{p_1}\in (0,1)$. To estimate $W^+_{p_1}$, we need to consider the expression of the differential equations in $\V_{p_1}^r$. 
Moreover, we recall from the definition of $\o^r$ that it fails to be in the symbol class $\Gamma_2^1$.
To get around this difficulty we consider two more cut-off functions $\chi_1$ and $\chi_2$ in $\Gamma_1^0$, such that both of their supports are in 
$$\V_{p_1}^r\cdot\R_+:=\{(t,x_1,x_2,\tau,\eta)\in\O\times\Pi:(t,x_1,x_2,\tau/\sqrt{|\tau|^2+\eta^2},\eta/\sqrt{|\tau|^2+\eta^2})\in\V_{p_1}^r\},$$
with $\chi_1=1$ on the support of $\chi_{p_1}$ and $\chi_2=1$ on the support of $\chi_1$.
Now we multiply \eqref{para1} by the symbol $\chi_2$ to obtain
\begin{align}\label{variablep1}
T^\gamma_{\chi_2\mathcal{A}^r}W^+_{p_1}+T^\gamma_{\chi_2A_0^rC^r}W^+_{p_1}+T^\gamma_rW^++I_2\p_2W^+_{p_1}=R_0F+R_{-1}W^+.
\end{align}
Here the support of the first order symbol $\chi_2\mathcal{A}^r$ is contained in the support of $\chi_2$ which is a subset of $\V_{p_1}^r\cdot\R_+$. So we actually exclude the part of the frequency where $\o^r$ vanishes. Now we can construct the transformation matrix to upper triangularize the first order symbol $\chi_2\mathcal{A}^r$ on the support of $\chi_2$. With the inspiration from the eigenvector $E^r$ of $\mathbb{A}^r$ in the previous subsection, we consider transformation matrix $Q_0^r$  which is homogeneous of degree 0 with respect to $(\tau,\eta)$ and takes the following form on $\V_{p_1}^r$:
\begin{align*}
{Q_0^r}:=\begin{pmatrix}
1 & \widehat{W}_1^r & 0 & 0 & 0 & 0 & 0\\
0 & -\alpha^r(\mu^r+\o^r) & U_2^r & 0 & 0 & 0 & 0\\
0 & -\alpha^rm^r & U_3^r & 0 & 0 & 0 & 0\\
0 & \widehat{W}_4^r & 0 & 1 & 0 & 0 & 0\\
0 & \widehat{W}_5^r & 0 & 0 & 1 & 0 & 0\\
0 & \widehat{W}_6^r & 0 & 0 & 0 & 1 & 0\\
0 & \widehat{W}_7^r & 0 & 0 & 0 & 0 & 1
\end{pmatrix},
\end{align*}
where the second and third rows of the second column are from the eigenvector $E^r$, and the choice of $\widehat{W}_i^r$, for $i=1,4,5,6,7$, is to guarantee that the second column of $\mathcal{A}^r{Q_0^r}$  is zero except for the second and third rows, that is
\begin{align*}
{\small \begin{pmatrix}
\tau+\ti v^r \eta & \frac{\ti\eta {c^r}^2}{\langle\p_1\Phi^r\rangle\rho} & \frac{\ti\eta {c^r}^2}{\langle\p_1\Phi^r\rangle\rho} & -\ti\eta F_{11}^r & 0 & -\ti\eta F_{12}^r & 0\\
-\ti F_{11}^r\eta & 0 & 0 & \tau+\ti v^r\eta & 0 & 0 & 0\\
0 & -\frac{\ti c^r F_{11}^r\eta}{\rho^r} & \frac{\ti c^r F_{11}^r\eta}{\rho^r} & 0 & \tau+\ti v^r\eta & 0 & 0\\
-\ti F_{12}^r\eta & 0 & 0 & 0 & 0 & \tau+\ti v^r\eta & 0\\
0 & -\frac{\ti c^r F_{12}^r\eta}{\rho^r} & \frac{\ti c^r F_{12}^r\eta}{\rho^r} & 0 & 0 & 0 & \tau+\ti v^r\eta
\end{pmatrix}\begin{pmatrix}
\widehat{W}_1^r\\
-\alpha^r(\mu^r+\o^r)\\
-\alpha^rm^r\\
\widehat{W}_4^r\\
\widehat{W}_5^r\\
\widehat{W}_6^r\\
\widehat{W}_7^r\end{pmatrix}=0}.
\end{align*} 
Note that we can solve the above equations for $\widehat{W}_i^r$ with $i=1,4,5,6,7$ at all points in $\R^3_+\times\Pi$. Because all the terms in each equation above have factors of $\tau+\ti v^r \eta$ or $(\tau+\ti v^r\eta)^2+({F_{11}^r}^2+{F_{12}^r}^2)\eta^2$, by cancelling the common factors we can obtain a linear algebraic system that never degenerates at any point in $\R^3_+\times\Pi$.  
However $\widehat{W}_i^r$ defined in this way is not homogeneous of degree $0$, because $\o^r$ degenerates near its zeros. We will use $\chi_1Q_0^r$ in the operations of the para-differential calculus below to exclude the frequencies where $\o^r$ degenerates. 
Moreover $U_2$ and $U_3$ are free to choose, given the matrix $Q_0^r$ invertible. Here to make the   arguments  simpler, we take $U_2=1$ and $U_3=0$ in this case. Thus $\chi_1Q_0^r$ is in the symbol class $\Gamma_2^0$.

To successively upper triangularize the first order operator ${A}^r$ in $\V_{p_1}^r\cdot\R_+$, we construct the symmetrizer $R_0^r$ that is homogeneous of degree $0$ and takes the following form on $\V_{p_1}^r$:
\begin{align*}
R_0^r=\begin{pmatrix}
1 & 0 & 0 & 0 & 0 & 0 & 0\\
0 & 0 & -\frac{1}{\xi} & 0 & 0 & 0 & 0\\
\bar{W}_1^r & \frac{\alpha^rm^r}{\xi} & -\frac{\alpha^r(\mu^r+\o^r)}{\xi} & \bar{W}_4^r & \bar{W}_5^r & \bar{W}_6^r & \bar{W}_7^r\\
0 & 0 & 0 & 1 & 0 & 0 & 0\\
0 & 0 & 0 & 0 & 1 & 0 & 0\\
0 & 0 & 0 & 0 & 0 & 1 & 0\\
0 & 0 & 0 & 0 & 0 & 0 & 1
\end{pmatrix},
\end{align*}
where $\xi$ equals the determinant of $Q_0^r$. $\bar{W}_1^r$, $\bar{W}_4^r$, $\bar{W}_5^r$, $\bar{W}_6^r$ and $\bar{W}_7^r$  are chosen to be homogeneous of degree 0 and make third row of $R_0^r\mathcal{A}^rQ_0^r$ to be zero except the second and third columns in this row, that is
\begin{align*}
{\small\begin{pmatrix}
\bar{W}_1^r\\
\frac{\alpha^rm^r}{\xi}\\
-\frac{\alpha^r(\mu^r+\o^r)}{\xi}\\
\bar{W}_4^r\\
\bar{W}_5^r\\
\bar{W}_6^r\\
\bar{W}_7^r
\end{pmatrix}^\top
\begin{pmatrix}
\tau+\ti v^r \eta  & -\ti\eta F_{11}^r & 0 & -\ti\eta F_{12}^r & 0\\
\frac{\ti\eta\p_2\Phi^r\rho^r}{2c^r\langle\p_1\Phi^r\rangle^2} &  0 & -\frac{\ti\eta\p_2\Phi^r\rho^rF_{11}^r}{2{c^r}^2\langle\p_1\Phi^r\rangle} & 0 & -\frac{\ti\eta\p_2\Phi^r\rho^rF_{12}^r}{2{c^r}^2\langle\p_1\Phi^r\rangle}\\
-\frac{\ti\eta\p_2\Phi^r\rho^r}{2c^r\langle\p_1\Phi^r\rangle^2} &  0 & -\frac{\ti\eta\p_2\Phi^r\rho^rF_{11}^r}{2{c^r}^2\langle\p_1\Phi^r\rangle} & 0 & -\frac{\ti\eta\p_2\Phi^r\rho^rF_{12}^r}{2{c^r}^2\langle\p_1\Phi^r\rangle}\\
-\ti F_{11}^r\eta  & \tau+\ti v^r\eta & 0 & 0 & 0\\
0  & 0 & \tau+\ti v^r\eta & 0 & 0\\
-\ti F_{12}^r\eta  & 0 & 0 & \tau+\ti v^r\eta & 0\\
0 & 0 & 0 & 0 & \tau+\ti v^r\eta
\end{pmatrix}=0}.
\end{align*}
Similarly as in $\chi_1Q_0^r$, we see that $\chi_1R_0^r$ is in the symbol class $\Gamma_2^0$. With the above choice of $Q_0^r$ and $R_0^r$, we can obtain our upper triangularized first order symbol:
\begin{align*}
&\widetilde{A}^r:=R_0^r\mathcal{A}^rQ_0^r=\\
&{\small\begin{pmatrix}
\tau+\ti v^r\eta & 0 & \Theta_1 & -\ti\eta F_{11}^r & 0 & -\ti\eta F_{12}^r & 0\\
\Theta_1 & -\o^r-\frac{\ti\p_2\Phi^r\p_1\Phi^r\eta}{\langle\p_1\Phi^r\rangle^2} & 0 & 0 & \Theta_1 & 0 & \Theta_1\\
0 & 0 & \o^r-\frac{\ti\p_2\Phi^r\p_1\Phi^r\eta}{\langle\p_1\Phi^r\rangle^2} & 0 & 0 & 0 & 0\\
-\ti\eta F_{11}^r & 0 & 0 & \tau+\ti v^r \eta & 0 & 0 & 0\\
0 & 0 & \Theta_1 & 0 & \tau+\ti v^r \eta & 0 & 0\\
-\ti\eta F_{12}^r & 0 & 0 & 0 & 0 & \tau+\ti v^r \eta & 0\\
0 & 0 & \Theta_1 & 0 & 0 & 0 & \tau+\ti v^r \eta
\end{pmatrix}},
\end{align*}
where $\Theta_1$ is some symbol in the class of $\Gamma_2^1$,  whose detailed expression is not important. 

Before we apply the transformation matrix $Q_0^r$ and symmetrizer $R_0^r$ to \eqref{paralinearizationsystem}, it is necessary to construct the matrix symbol $Q_{-1}^r$ and $R_{-1}^r$ in $\Gamma_1^{-1}$ to
decouple the incoming modes and outgoing mode in zero order terms. To be specific, we obtain the following lemma:
\begin{Lemma}\label{diag}
With the appropriate choice of $Q_{-1}^r$ and $R_{-1}^r$ in $\Gamma_1^{-1}$, there is a symbol $D_0=({d}_{i,j})_{7\times7}$ in $\Gamma_1^0$ satisfying $d_{2,3}=d_{3,2}=0$,  such that 
\begin{align*}
R_{-1}^r{R_0^r}^{-1}\widetilde{A}^{r}-\widetilde{A}^{r}Q^r_{-1}Q^r_0+(\p_2{Q_0^r}^{-1}-R^rA_0^rC^r-[R^r,\chi_2\mathcal{A}^r]+[\chi_2\mathcal{A}^r,Q^r])Q_0^r-D_0
\end{align*}
is a symbol that is homogeneous of degree $-1$ and has regularity $1$ on $\chi_2=1$; moreover
\begin{align*}
R_{-1}^rI_2=I_2Q_{-1}^r,
\end{align*}
where $[\cdot,\cdot]$   is defined as 
\begin{align*}
[A,B]:=\frac{1}{i}\left(\frac{\p A}{\p \delta}\frac{\p B}{\p t}+\frac{\p A}{\p \eta}\frac{\p B}{\p x_1}\right),
\end{align*}
for any symbols A and B.
\end{Lemma}
The proof of this lemma is in the same spirit as the Lemma 5.5 in \cite{coulombel2004stability} and hence is omitted here.
\begin{Remark}
Actually the above lemma is also true for the neighborhoods $\V_{p_1}^l$, $\V_r^i$, $\V_{p_2}^j$ for $i=1,2,3$ and $j=1,2,3,4$, if we construct transformation matrices $Q_0^r$, $Q_{-1}^r$ and symmetrizers $R_0^r$, $R_{-1}^r$ in the same way as above.
\end{Remark}

Then we consider the new unknown
\begin{align*}
Z^+=T^\gamma_{\chi_1({Q_0^r}^{-1}+Q_{-1}^r)}W_{p_1}^+,
\end{align*}
and set $Q^r={Q_0^r}^{-1}+Q_{-1}^r$ and $R^r=R_0^r+R_{-1}^r$. We have
\begin{align*}
I_2\p_2Z^+&=I_2T^\gamma_{(\p_2\chi_1)Q^r}W_{p_1}^++I_2T^\gamma_{\chi_1\p_2Q^r}W_{p_1}^++I_2T^\gamma_{\chi_1Q^r}\p_2W_{p_1}^+\\
&=I_2T^\gamma_{(\p_2\chi_1)Q^r}W_{p_1}^++I_2T^\gamma_{\chi_1\p_2Q^r}W_{p_1}^++T^\gamma_{\chi_1R^r}I_2\p_2W_{p_1}^+.
\end{align*}
In the last expression,   $\p_2\chi_1$ is only supported in the place where $\chi_1\in(0,1)$ and by definition this region is disjoint with the support of $\chi_{p_1}$. By the asymptotic expansion of the symbols, it follows that 
\begin{align*}
T^\gamma_{(\p_2\chi_1)Q}W_{p_1}^+=R_{-1}W^+.
\end{align*}

Recalling \eqref{variablep1} and utilizing the property of para-differential calculus, we have
\begin{align*}
&I_2\p_2Z^+=I_2T^\gamma_{\chi_1\p_2{Q_0^r}^{-1}}W_{p_1}^+-T^\gamma_{\chi_1R^r\mathcal{A}^{r}}W_{p_1}^+-T^\gamma_{[\chi_1R^r,\chi_2\mathcal{A}^r]}W_{p_1}^+-T^\gamma_{\chi_1R^rA_0^rC^r}W_{p_1}^+\\
&\qquad\qquad+T^\gamma_rW^++R_0F+R_{-1}W^+\\
&=I_2T^\gamma_{\chi_1\p_2{Q_0^r}^{-1}}W_{p_1}^+-T^\gamma_{\chi_1R^r\mathcal{A}^{r}}W_{p_1}^++T^\gamma_{\chi_2\widetilde{A}^{r}}T^\gamma_{\chi_1Q^r}W^+_{p_1}-T^\gamma_{\chi_2\widetilde{A}^{r}}T^\gamma_{\chi_1Q^r}W^+_{p_1}-T^\gamma_{[\chi_1R^r,\chi_2\mathcal{A}^r]}W_{p_1}^+\\
&\qquad\qquad-T^\gamma_{\chi_1R^rA_0^rC^r}W_{p_1}^++T^\gamma_rW^++R_0F+R_{-1}W^+\\
&=-T^\gamma_{\chi_1(R_{-1}{R_0}^{-1}\widetilde{A}^{r}-\widetilde{A}^{r}Q_{-1}Q_0)}Z^+-T^\gamma_{\chi_2\widetilde{A}^{r}}Z^++I_2T^\gamma_{\chi_1\p_2{Q_0^r}^{-1}}W_{p_1}^+-T^\gamma_{\chi_1R^rA_0^rC^r}W_{p_1}^+\\
&\qquad\qquad-T^\gamma_{[\chi_1R^r,\chi_2\mathcal{A}^r]-[\chi_2\widetilde{A}^r,\chi_1Q^r]}W_{p_1}^++T^\gamma_rW^++R_0F+R_{-1}W^+.
\end{align*}
Then by the asymptotic expansion and the support of each cut-off function, one has
\begin{align*}
I_2\p_2Z^+= & -T^\gamma_{\chi_1(R_{-1}{R_0}^{-1}\widetilde{A}^{r}-\widetilde{A}^{r}Q_{-1}Q_0)}Z^+-T^\gamma_{\chi_2\widetilde{A}^{r}}Z^++I_2T^\gamma_{\chi_1\p_2{Q_0^r}^{-1}Q_0^r}Z^+-T^\gamma_{\chi_1R^rA_0^rC^rQ_0^r}Z^+\\
& -T^\gamma_{\chi_1([R^r,\chi_2\mathcal{A}^r]-[\chi_2\widetilde{A}^r,Q^r])Q_0^r}Z^++T^\gamma_rW^++R_0F+R_{-1}W^+, 
\end{align*}
and Lemma \ref{diag} implies
\begin{align}\label{uppertrianglecase1}
I_2\p_2Z^+=-T^\gamma_{\chi_2\widetilde{A}^{r}}Z^++T^\gamma_{D_0}Z^++T^\gamma_rW^++R_0F+R_{-1}W^+.
\end{align}
By the definition of $Z^+$, it is obvious that the support of the Fourier transform of $Z^+$ is in the support of $\chi_{p_1}$. So the above equation can be rewrite as
\begin{align}\label{de}
I_2\p_2Z^+=-T^\gamma_{\widetilde{D_1}}Z^++T^\gamma_{\widetilde{D_0}}Z^++T^\gamma_rW^++R_0F+R_{-1}W^+,
\end{align}
where $\widetilde{D_1}$ is the same as  $\widetilde{A}^{r}$ except replacing $\o^r$ in each element by $\tilde{\o}^r$, where $\tilde{\o}^r$ is a symbol of degree $1$ with regularity $2$ which equals $\o^r$ on the support of $\chi_2$, and $\widetilde{D_0}$ is an extension of $D_0$ with  $d_{2,3}=d_{3,2}=0$ in the whole space. Moreover, by a direct computation, we see that $\o^r\geq c \Lambda:= c\sqrt{\gamma^2+\delta^2+\eta^2}$ in $\V_{p_1}^r$. This suggests that we can choose the extension $\tilde{\o}^r$ such that $\tilde{\o}^r\geq c \Lambda$  in the whole space. In the following context, for simplicity of the expression,  we will drop the tilde in $\tilde{\o}^r$.

Next we are going to apply the appropriate symmetrizer to  \eqref{de} and obtain the estimate. By the definition of $\widetilde{D_1}$,   
  $\widetilde{D_1}$   degenerates near the pole   in this case, which weakens the estimate of the characteristic part of the unknown, because at the pole the leading term has homogeneity of degree $0$ and we can not neglect the effect of the zeroth order term. To overcome this difficulty, we shall consider different symmetrizers for each equation in the system.

We denote $Z^+=(Z_1,Z_2,...,Z_7)^\top$ and start from the third equation in \eqref{de},  corresponding to the outgoing mode of the system,
\begin{align*}
\p_2Z_3=T^\gamma_{-\o^r+\ti\varpi^r}Z_3+T^\gamma_{\Theta_0}Z_3+\Sigma_{i\neq2,3}T^\gamma_{\Theta_0}Z_i+T^\gamma_rW^++R_0F+R_{-1}W^+,
\end{align*}
where $\varpi^r=\frac{\p_2\Phi^r\p_1\Phi^r}{\langle\p_1\Phi^r\rangle^2} \eta$. For this equation, we consider the following two operators $(T^\gamma_{\sigma})^*T^\gamma_{\Lambda}T^\gamma_{\sigma}$ and $(T^\gamma_{\Lambda})^*T^\gamma_{\Lambda}$ as symmetrizers, where $\sigma$ is defined on $\R^3_+\times\Sigma$ by \eqref{rootsofL},  and extend the domain to the whole frequency space $\R^3_+\times\Pi$ by homogeneity of degree 1. Thus $\sigma$ is in the symbol class $\Gamma_2^1$, and
\begin{equation}\label{variableZ3_1}
\begin{split}
\Re\langle T^\gamma_{\sigma}\p_2Z_3, T^\gamma_{\Lambda}T^\gamma_{\sigma}Z_3\rangle & = \Re\langle T^\gamma_{\Lambda}T^\gamma_{\sigma}Z_3, T^\gamma_{\sigma}\p_2Z_3\rangle \\
& = \Re\langle T^\gamma_{\Lambda}T^\gamma_{\sigma}Z_3, T^\gamma_{\sigma}T^\gamma_{-\o^r+\ti\varpi^r}Z_3\rangle
+\Re\langle T^\gamma_{\Lambda}T^\gamma_{\sigma}Z_3, T^\gamma_{\sigma}T^\gamma_{\Theta_0}Z_3\rangle \\
&\quad +\Sigma_{i\neq2,3}\Re\langle T^\gamma_{\Lambda}T^\gamma_{\sigma}Z_3, T^\gamma_{\sigma}T^\gamma_{\Theta_0}Z_i\rangle +\Re\langle T^\gamma_{\Lambda}T^\gamma_{\sigma}Z_3, T^\gamma_{\sigma}T^\gamma_{r}W^+\rangle \\
&\quad +\Re\langle T^\gamma_{\Lambda}T^\gamma_{\sigma}Z_3, T^\gamma_{\sigma}R_{-1}W^+\rangle+\Re\langle T^\gamma_{\Lambda}T^\gamma_{\sigma}Z_3, T^\gamma_{\sigma}F\rangle.
\end{split}
\end{equation}
For the first term on the right side, by the para-differential calculus, we have
\begin{align*}
\Re\langle T^\gamma_{\Lambda}T^\gamma_{\sigma}Z_3, T^\gamma_{\sigma}T^\gamma_{-\o^r+\ti\varpi^r}Z_3\rangle=\Re\langle T^\gamma_{\Lambda}T^\gamma_{\sigma}Z_3, T^\gamma_{\frac{-\o^r+\ti\varpi^r}{\Lambda}}T^\gamma_{\Lambda}T^\gamma_{\sigma}Z_3\rangle+\Re\langle T^\gamma_{\Lambda}T^\gamma_{\sigma}Z_3, R_1Z_3\rangle.
\end{align*}
From the way of extending  $\o^r$, we have 
\begin{align*}
\Re \frac{-\o^r+\ti\varpi^r}{\Lambda}\geq c
\end{align*}
for some positive $c$   depending only  on the background state. By G\r{a}rding's inequality, we have
\begin{align*}
\Re\langle T^\gamma_{\Lambda}T^\gamma_{\sigma}Z_3, T^\gamma_{\frac{-\o^r+\ti\varpi^r}{\Lambda}}T^\gamma_{\Lambda}T^\gamma_{\sigma}Z_3\rangle\geq c\|T^\gamma_{\Lambda}T^\gamma_{\sigma}Z_3\|_0^2=c\|T^\gamma_{\sigma}Z_3\|_{1,\gamma}^2.
\end{align*}
For the other terms on the right  side of \eqref{variableZ3_1}, we have
\begin{align*}
&\Re\langle T^\gamma_{\Lambda}T^\gamma_{\sigma}Z_3, R_1Z_3\rangle\leq \e\|T^\gamma_{\Lambda}T^\gamma_{\sigma}Z_3\|_0^2+\frac{1}{\e}\|Z_3\|_{1,\gamma}^2,\\
&\Re\langle T^\gamma_{\Lambda}T^\gamma_{\sigma}Z_3, T^\gamma_{\sigma}T^\gamma_{\Theta_0}Z_3\rangle\leq \e\|T^\gamma_{\Lambda}T^\gamma_{\sigma}Z_3\|_0^2+\frac{1}{\e}\|Z_3\|_{1,\gamma}^2,\\
&\Re\langle T^\gamma_{\Lambda}T^\gamma_{\sigma}Z_3, T^\gamma_{\sigma}T^\gamma_{\Theta_0}Z_i\rangle
=\Re\langle T^\gamma_{\Lambda}T^\gamma_{\sigma}Z_3, T^\gamma_{\Theta_0}T^\gamma_{\sigma}Z_i\rangle+\Re\langle T^\gamma_{\Lambda}T^\gamma_{\sigma}Z_3, R_0Z_i\rangle\\
&\hspace{1.4in} \leq\e\|T^\gamma_{\Lambda}T^\gamma_{\sigma}Z_3\|_0^2+\frac{1}{\e}\|T^\gamma_{\sigma}Z_i\|_0^2+\frac{1}{\e}\|Z_i\|_0^2,\\
&\Re\langle T^\gamma_{\Lambda}T^\gamma_{\sigma}Z_3, T^\gamma_{\sigma}T^\gamma_{r}W^+\rangle\leq\e\|T^\gamma_{\Lambda}T^\gamma_{\sigma}Z_3\|_0^2+\frac{1}{\e}\|T^\gamma_{r}W^+\|_{1,\gamma}^2,\\
&\Re\langle T^\gamma_{\Lambda}T^\gamma_{\sigma}Z_3, T^\gamma_{\sigma}R_{-1}W^+\rangle\leq\e\|T^\gamma_{\Lambda}T^\gamma_{\sigma}Z_3\|_0^2+\frac{1}{\e}\|W^+\|_{0}^2,\\
&\Re\langle T^\gamma_{\Lambda}T^\gamma_{\sigma}Z_3, T^\gamma_{\sigma}F\rangle\leq\e\|T^\gamma_{\Lambda}T^\gamma_{\sigma}Z_3\|_0^2+\frac{1}{\e}\|F\|_{1,\gamma}^2,
\end{align*}
where $\e>0$ is sufficiently small. To achieve the desired estimate, we need to deal with the term with $x_2$ derivative. We note
\begin{equation}\label{variableZ3_2}
\begin{split}
\p_2\Re\langle Z_3, (T^\gamma_{\sigma})^*T^\gamma_{\Lambda}T^\gamma_{\sigma}Z_3\rangle & =\Re\p_2\langle T^\gamma_{\sigma}Z_3, T^\gamma_{\Lambda}T^\gamma_{\sigma}Z_3\rangle\\
&=\Re\langle T^\gamma_{\p_2\sigma}Z_3, T^\gamma_{\Lambda}T^\gamma_{\sigma}Z_3\rangle
+\Re\langle T^\gamma_{\sigma}Z_3, T^\gamma_{\Lambda}T^\gamma_{\p_2\sigma}Z_3\rangle\\
&\quad +\Re\langle T^\gamma_{\sigma}Z_3, T^\gamma_{\Lambda}T^\gamma_{\sigma}\p_2Z_3\rangle+\Re\langle T^\gamma_{\sigma}\p_2Z_3, T^\gamma_{\Lambda}T^\gamma_{\sigma}Z_3\rangle.
\end{split}
\end{equation}
The last term on the right side is what we want to estimate. For the first three terms on the right  side of \eqref{variableZ3_2}, we have
\begin{align*}
&\Re\langle T^\gamma_{\p_2\sigma}Z_3, T^\gamma_{\Lambda}T^\gamma_{\sigma}Z_3\rangle\leq\e\|T^\gamma_{\Lambda}T^\gamma_{\sigma}Z_3\|_0^2+\frac{1}{\e}\|Z_3\|_{1,\gamma}^2,\\
&\Re\langle T^\gamma_{\sigma}Z_3, T^\gamma_{\Lambda}T^\gamma_{\p_2\sigma}Z_3\rangle\leq\e\|T^\gamma_{\Lambda}T^\gamma_{\sigma}Z_3\|_0^2+\frac{1}{\e}\|Z_3\|_{1,\gamma}^2,\\
&\Re\langle T^\gamma_{\sigma}Z_3, T^\gamma_{\Lambda}T^\gamma_{\sigma}\p_2Z_3\rangle=\Re\langle T^\gamma_{\Lambda}T^\gamma_{\sigma}Z_3, T^\gamma_{\sigma}\p_2Z_3\rangle+\Re\langle T^\gamma_{\sigma}Z_3, R_0T^\gamma_{\sigma}\p_2Z_3\rangle.
\end{align*}
For the term $\Re\langle T^\gamma_{\sigma}Z_3, R_0T^\gamma_{\sigma}\p_2Z_3\rangle$, we can handle it in a similar way to $\Re\langle T^\gamma_{\Lambda}T^\gamma_{\sigma}Z_3, T^\gamma_{\sigma}\p_2Z_3\rangle$. Adding \eqref{variableZ3_1} and \eqref{variableZ3_2} and integrating with respect to $x_2$, we  obtain
\begin{align*}
&\vertiii{T^\gamma_{\Lambda}T^\gamma_{\sigma}Z_3}_0^2+\Re\langle T^\gamma_{\sigma}Z_3, T^\gamma_{\Lambda}T^\gamma_{\sigma}Z_3\rangle|_{x_2=0}\lesssim \left(C+\frac{1}{\e}\right)\vertiii{Z_3}_{1,\gamma}^2+\Sigma_{i\neq2,3}\frac{1}{\e}\left(\vertiii{T^\gamma_{\sigma}Z_i}_0^2+\vertiii{Z_i}_0^2\right)\\
&\qquad\qquad\qquad\qquad\qquad\qquad\qquad\qquad\qquad+\frac{1}{\e}\left(\vertiii{T_rW^+}_{1,\gamma}^2+\vertiii{W^+}_0^2+\vertiii{F}_{1,\gamma}^2\right).
\end{align*}
From
\begin{align*}
\Re\langle T^\gamma_{\sigma}Z_3, T^\gamma_{\Lambda}T^\gamma_{\sigma}Z_3\rangle|_{x_2=0}=\Re\langle T^\gamma_{\Lambda^\frac{1}{2}}T^\gamma_{\sigma}Z_3, T^\gamma_{\Lambda^\frac{1}{2}}T^\gamma_{\sigma}Z_3\rangle|_{x_2=0}+\langle T^\gamma_{\sigma}Z_3, R_0T^\gamma_{\sigma}Z_3\rangle|_{x_2=0},
\end{align*}
one has
\begin{equation}\label{variableZ3main1}
\begin{split}
\vertiii{T^\gamma_{\Lambda}T^\gamma_{\sigma}Z_3}_0^2+\|T^\gamma_{\Lambda^\frac{1}{2}}T^\gamma_{\sigma}Z_3|_{x_2=0}\|_0^2 & \lesssim \|T^\gamma_{\sigma}Z_3|_{x_2=0}\|_0^2+\left(C+\frac{1}{\e}\right)\vertiii{Z_3}_{1,\gamma}^2\\
&\quad +\Sigma_{i\neq2,3}\frac{1}{\e}\left(\vertiii{T^\gamma_{\sigma}Z_i}_0^2+\vertiii{Z_i}_0^2\right) \\
&\quad +\frac{1}{\e}\left(\vertiii{T_rW^+}_{1,\gamma}^2+\vertiii{W^+}_0^2+\vertiii{F}_{1,\gamma}^2\right).
\end{split}
\end{equation}

For the second symmetrizer $(T^\gamma_{\Lambda})^*T^\gamma_{\Lambda}$, we have
\begin{align*}
&\p_2\Re\langle T^\gamma_{\Lambda}Z_3,T^\gamma_{\Lambda}Z_3\rangle \\
& =2\Re\langle T^\gamma_{\Lambda}Z_3,T^\gamma_{\Lambda}\p_2Z_3\rangle\\
& =2\Re\langle T^\gamma_{\Lambda}Z_3,T^\gamma_{\Lambda}T^\gamma_{-\o^r+\ti\varpi^r}Z_3\rangle+2\Re\langle T^\gamma_{\Lambda}Z_3,T^\gamma_{\Lambda}T^\gamma_{\Theta_0}Z_3\rangle+2\Sigma_{i\neq2,3}\Re\langle T^\gamma_{\Lambda}Z_3,T^\gamma_{\Lambda}T^\gamma_{\Theta_0}Z_i\rangle\\
&\quad +2\Re\langle T^\gamma_{\Lambda}Z_3,T^\gamma_{\Lambda}T^\gamma_{r}W^+\rangle+2\Re\langle T^\gamma_{\Lambda}Z_3,T^\gamma_{\Lambda}R_{-1}W^+\rangle+2\Re\langle T^\gamma_{\Lambda}Z_3,T^\gamma_{\Lambda}F\rangle.
\end{align*}
A similar argument to the above yields the following:
\begin{equation}\label{variableZ3main2}
\begin{split}
\vertiii{Z_3}_{\frac{3}{2},\gamma}^2+\|Z_3|_{x_2=0}\|_{1,\gamma}^2 & \lesssim C\vertiii{Z_3}_{1,\gamma}^2+\frac{1}{\e}\vertiii{Z_3}_{\frac{1}{2},\gamma}^2+\Sigma_{i\neq2,3}(\frac{1}{\e}\vertiii{Z_i}_{\frac{1}{2},\gamma}^2+\frac{1}{\e\gamma}\vertiii{Z_i}_0^2)\\
&\quad+\frac{1}{\e\gamma}\left(\vertiii{T_rW^+}_{1,\gamma}^2+\vertiii{W^+}_0^2+\vertiii{F}_{1,\gamma}^2\right)
\end{split}
\end{equation}

Then we consider the 1st, 4th and 6th equations of the system \eqref{de} which can be written as
\begin{align}\label{para2}
T^\gamma_{\bf a}(Z_1, Z_4, Z_6)^\top+T^\gamma_{{\bf\Theta}_1}Z_3+ T^\gamma_{{\bf\Theta}_0}Z^++T^\gamma_{r}W^++R_{-1}W^+=R_0F,
\end{align}
where ${\bf\Theta}_1$ is a $3\times1$ matrix symbol of degree 1 and regularity 2, ${\bf\Theta}_0$ is a $3\times7$ matrix symbol of degree 0 and regularity 1, and ${\bf a}$ is a $3\times3$ matrix symbol of degree 1 and regularity 2 taking the following form
\begin{align*}
{\bf a}=\begin{pmatrix}
\tau+\ti v^r\eta & -\ti\eta F_{11}^r & -\ti\eta F_{12}^r \\
-\ti\eta F_{11}^r & \tau+\ti v^r\eta & 0\\
-\ti\eta F_{12}^r & 0 &  \tau+\ti v^r\eta
\end{pmatrix}.
\end{align*}
We first apply the symbol $\frac{{\bf a}^*}{\Lambda^2}\in\Gamma^0_2$ to equation \eqref{para2} with ${\bf a}^*$ the adjoint matrix of ${\bf a}$, and obtain
\begin{align*}
T^\gamma_{a}Z_j+T^\gamma_{\Theta_1}Z_3+\Sigma T^\gamma_{\Theta_0}Z_i+T^\gamma_{r}W^++R_{-1}W^+=R_0F,
\end{align*} 
where $j=1,4,6$ and $a=(\tau+\ti v^r\eta)\left((\tau+\ti v^r\eta)^2+({F_{11}^r}^2+{F_{12}^r}^2)\eta^2\right)/\Lambda^2$. 
From the way we construct the cut-off function, $\left((\tau+\ti v^r\eta)^2+({F_{11}^r}^2+{F_{12}^r}^2)\eta^2\right)/\Lambda^2$ is not zero on the support of $\chi_2$.  So we can rewrite
\begin{align*}
a=(1-\chi_2)a+\chi_2(\tau+\ti v^r\eta)\left((\tau+\ti v^r\eta)^2+({F_{11}^r}^2+{F_{12}^r}^2)\eta^2\right)/\Lambda^2.
\end{align*}
Because $Z^+=T^\gamma_{\chi_1(Q_0+Q_{-1})}W_{p_1}^+$, we have
\begin{align*}
T^\gamma_{(1-\chi_2)a}Z_j=T^\gamma_{(1-\chi_2)a}T^\gamma_{\chi_1Q_j^r}T^\gamma_{\chi_{p_2}}W^+=R_{-1}T^\gamma_{\chi_{p_2}}W^+.
\end{align*}
The above is from the fact that the support of $(1-\chi_2)a$ is disjoint with $\chi_1Q_j^r$'s. Since $\sigma=0$ only holds at the points where $\tau+\ti v^r\eta=0$ in the support of $\chi_2$, we have $\chi_2(\tau+\ti v^r\eta)\left((\tau+\ti v^r\eta)^2+({F_{11}^r}^2+{F_{12}^r}^2)\eta^2\right)/\Lambda^2=\chi_2\Theta_0\times(\gamma+\ti\sigma)$, and hence,
\begin{align}\label{variableZj_1}
T^\gamma_{\chi_2\Theta_0\times(\gamma+\ti\sigma)}Z_j+T^\gamma_{\Theta_1}Z_3+\Sigma T^\gamma_{\Theta_0}Z_i+T^\gamma_{r}W^++R_{-1}W^+=R_0F.
\end{align}
For the above equation, we consider the symmetrizer $(T^\gamma_{\sigma})^*T^\gamma_{\sigma}$ and obtain 
\begin{equation*}
\begin{split}
&\Re\langle T^\gamma_{\sigma}Z_j,T^\gamma_{\sigma}T^\gamma_{\chi_2\Theta_0\times(\gamma+\ti\sigma)}Z_j\rangle+\Re\langle T^\gamma_{\sigma}Z_j,T^\gamma_{\sigma}T^\gamma_{\Theta_1}Z_3\rangle+\Sigma\Re\langle T^\gamma_{\sigma}Z_j,T^\gamma_{\sigma}T^\gamma_{\Theta_0}Z_i\rangle\\
&\qquad +\Re\langle T^\gamma_{\sigma}Z_j,T^\gamma_{\sigma}T^\gamma_{r}W^+\rangle+\Re\langle T^\gamma_{\sigma}Z_j,T^\gamma_{\sigma}R_{-1}W^+\rangle =\Re\langle T^\gamma_{\sigma}Z_j,T^\gamma_{\sigma}F\rangle.
\end{split}
\end{equation*}
In the above equality, all the terms except the first one can be easily estimated as follows:
\begin{align*}
&\Re\langle T^\gamma_{\sigma}Z_j,T^\gamma_{\sigma}T^\gamma_{\Theta_1}Z_3\rangle=\Re\langle T^\gamma_{\sigma} Z_j,T^\gamma_{\frac{\Theta_1}{\Lambda}}T^\gamma_{\Lambda}T^\gamma_{\sigma}Z_3\rangle+\Re\langle T^\gamma_{\sigma}Z_j, R_1Z_3\rangle\\
&\qquad\qquad\qquad\qquad \leq \e\gamma\|T^\gamma_{\sigma}Z_j\|_0^2+\frac{1}{\e\gamma}\|T^\gamma_{\Lambda}T^\gamma_{\sigma}Z_3\|_0^2+\frac{1}{\e\gamma}\|Z_3\|_{1,\gamma}^2,\\
&\Re\langle T^\gamma_{\sigma}Z_j,T^\gamma_{\sigma}T^\gamma_{\Theta_0}Z_i\rangle=\Re\langle T^\gamma_{\sigma}Z_j,T^\gamma_{\Theta_0}T^\gamma_{\sigma}Z_i\rangle+\Re\langle T^\gamma_{\sigma}Z_j,R_0Z_i\rangle\\
&\qquad\qquad\qquad\qquad \leq\e\gamma\|T^\gamma_{\sigma}Z_j\|_0^2+\frac{1}{\e\gamma}\|T^\gamma_{\sigma}Z_i\|_0^2+\frac{1}{\e\gamma}\|Z_i\|_0^2,\\
&\Re\langle T^\gamma_{\sigma}Z_j,T^\gamma_{\sigma}T^\gamma_{r}W^+\rangle\leq\e\gamma\|T^\gamma_{\sigma}Z_j\|_0^2+\frac{1}{\e\gamma}\|T^\gamma_{r}W^+\|_{1,\gamma}^2,\\
&\Re\langle T^\gamma_{\sigma}Z_j,T^\gamma_{\sigma}R_{-1}W^+\rangle\leq\e\gamma\|T^\gamma_{\sigma}Z_j\|_0^2+\frac{1}{\e\gamma}\|W^+\|_{0}^2,\\
&\Re\langle T^\gamma_{\sigma}Z_j,T^\gamma_{\sigma}F\rangle\leq\e\gamma\|T^\gamma_{\sigma}Z_j\|_0^2+\frac{1}{\e\gamma}\|F\|_{1,\gamma}^2.
\end{align*}
For the first term, we have
\begin{align*}
&\Re\langle T^\gamma_{\sigma}Z_j,T^\gamma_{\sigma}T^\gamma_{\chi_2\Theta_0\times(\gamma+\ti\sigma)}Z_j\rangle\\
&\qquad\qquad=\Re\langle T^\gamma_{\sigma}Z_j,T^\gamma_{\sigma}T^\gamma_{\chi_2\Theta_0}T^\gamma_{\gamma+\ti\sigma}Z_j\rangle+\Re\langle T^\gamma_{\sigma}Z_j,T^\gamma_{\Theta_0}T^\gamma_{\sigma}Z_j\rangle+\Re\langle T^\gamma_{\sigma}Z_j,R_{0}Z_j\rangle.
\end{align*}
The last two terms in the above can be estimated as
\begin{align*}
&\Re\langle T^\gamma_{\sigma}Z_j,R_{0}Z_j\rangle\leq\e\gamma\|T^\gamma_{\sigma}Z_j\|_0^2+\frac{1}{\e\gamma}\|Z_j\|_0^2,\\
&\Re\langle T^\gamma_{\sigma}Z_j,T^\gamma_{\Theta_0}T^\gamma_{\sigma}Z_j\rangle\leq C\|T^\gamma_{\sigma}Z_j\|_0^2,
\end{align*}
 and the first one can be further written as
\begin{align}
\Re\langle T^\gamma_{\sigma}Z_j,T^\gamma_{\sigma}T^\gamma_{\chi_2\Theta_0}T^\gamma_{\gamma+\ti\sigma}Z_j\rangle=\Re\langle T^\gamma_{\sigma}Z_j,T^\gamma_{\chi_2\Theta_0}T^\gamma_{\sigma}T^\gamma_{\gamma+\ti\sigma}Z_j\rangle+\Re\langle T^\gamma_{\sigma}Z_j,R_0T^\gamma_{\gamma+\ti\sigma}Z_j\rangle.  \label{E510}
\end{align}
In \eqref{E510}, the last term   implies
\begin{align*}
 \Re\langle T^\gamma_{\sigma}Z_j,R_0T^\gamma_{\gamma+\ti\sigma}Z_j\rangle&=\Re\langle T^\gamma_{\sigma}Z_j,R_0T^\gamma_{\gamma}Z_j\rangle+\Re\langle T^\gamma_{\sigma}Z_j,R_0T^\gamma_{\ti\sigma}Z_j\rangle\\
&\leq\e\gamma\|T^\gamma_{\sigma}Z_j\|_0^2+\frac{1}{\e\gamma}\|T^\gamma_\gamma Z_j\|_0^2+C\|T^\gamma_{\sigma}Z_j\|_0^2,
\end{align*}
and the first term leads to
\begin{align}
&\Re\langle T^\gamma_{\sigma}Z_j,T^\gamma_{\chi_2\Theta_0}T^\gamma_{\sigma}T^\gamma_{\gamma+\ti\sigma}Z_j\rangle=\Re\langle T^\gamma_{\sigma}Z_j,T^\gamma_{\chi_2\Theta_0}T^\gamma_{\sigma}T^\gamma_{\gamma}Z_j\rangle+\Re\langle T^\gamma_{\sigma}Z_j,T^\gamma_{\chi_2\Theta_0}T^\gamma_{\sigma}T^\gamma_{\ti\sigma}Z_j\rangle \notag\\
&=\gamma\Re\langle T^\gamma_{\sigma}Z_j,T^\gamma_{\chi_2\Theta_0}R_0Z_j\rangle+\Re\langle T^\gamma_{\sigma}Z_j,T^\gamma_{\chi_2\Theta_0}T^\gamma_{\gamma}T^\gamma_{\sigma}Z_j\rangle+\Re\langle T^\gamma_{\sigma}Z_j,T^\gamma_{\chi_2\Theta_0}T^\gamma_{\ti\sigma}T^\gamma_{\sigma}Z_j\rangle \notag\\
&=\gamma\Re\langle T^\gamma_{\sigma}Z_j,T^\gamma_{\chi_2\Theta_0}R_0Z_j\rangle+\Re\langle T^\gamma_{\sigma}Z_j,T^\gamma_{\chi_2\Theta_0}T^\gamma_{\gamma+\ti\sigma}T^\gamma_{\sigma}Z_j\rangle \notag\\
&=\gamma\Re\langle T^\gamma_{\sigma}Z_j,T^\gamma_{\chi_2\Theta_0}R_0Z_j\rangle+\Re\langle T^\gamma_{\sigma}Z_j,T^\gamma_{\chi_2\Theta_0(\gamma+\ti\sigma)}T^\gamma_{\sigma}Z_j\rangle+\Re\langle T^\gamma_{\sigma}Z_j,R_0T^\gamma_{\sigma}Z_j\rangle.  \label{E511}
\end{align}
In   \eqref{E511}, the first and the third terms can be estimated directly by the Cauchy inequality,
and for the second term, we have
\begin{equation}\label{variableZj_2}
\begin{split}
\Re\langle T^\gamma_{\sigma}Z_j,T^\gamma_{\chi_2\Theta_0(\gamma+\ti\sigma)}T^\gamma_{\sigma}Z_j\rangle&=\Re\langle T^\gamma_{\sigma}Z_j,T^\gamma_{\chi_2a}T^\gamma_{\sigma}Z_j\rangle\\
&=\Re\langle T^\gamma_{\sigma}Z_j,T^\gamma_{\tilde{a}}T^\gamma_{\sigma}Z_j\rangle+\Re\langle T^\gamma_{\sigma}Z_j,T^\gamma_{(\chi_2-1)\tilde{a}}T^\gamma_{\sigma}Z_j\rangle,
\end{split}
\end{equation}
where $\tilde{a}$ is an extension of $\chi_2 a$ to the whole space with $|\Re\tilde{a}|\geq c\gamma$ for some generic positive constant $c$,  and in the last term of \eqref{variableZj_2},
\begin{align*}
T^\gamma_{(\chi_2-1)\tilde{a}}T^\gamma_{\sigma}Z_j & =T^\gamma_{(\chi_2-1)\tilde{a}}T^\gamma_{\sigma}T^\gamma_{\chi_1Q_j^r}T^\gamma_{\chi_{p_2}}W^+\\
&=T^\gamma_{(\chi_2-1)\tilde{a}}T^\gamma_{\chi_1Q_j^r}T^\gamma_{\sigma}T^\gamma_{\chi_{p_2}}W^++T^\gamma_{(\chi_2-1)\tilde{a}}T^\gamma_{O_0}T^\gamma_{\chi_{p_2}}W^+\\
&\quad +T^\gamma_{(\chi_2-1)\tilde{a}}T^\gamma_{O_{-1}}T^\gamma_{\chi_{p_2}}W^++T^\gamma_{(\chi_2-1)\tilde{a}}R_{-2}T^\gamma_{\chi_{p_2}}W^+.
\end{align*}
Note that $O_{0}$ and $O_{-1}$ only support on the support of $\chi_1$ which is disjoint with the support of $(\chi_2-1)\tilde{a}$. 
Thus, we have
\begin{align*}
&T^\gamma_{(\chi_2-1)\tilde{a}}T^\gamma_{O_0}T^\gamma_{\chi_{p_2}}W^+=R_{-1}W^+,\\
&T^\gamma_{(\chi_2-1)\tilde{a}}T^\gamma_{O_{-1}}T^\gamma_{\chi_{p_2}}W^+=R_{-1}W^+,\\
&T^\gamma_{(\chi_2-1)\tilde{a}}T^\gamma_{\chi_1Q^j}T^\gamma_{\sigma}T^\gamma_{\chi_{p_2}}W^+=R_{-1}W^+,
\end{align*}
and hence,
\begin{align}\label{variableZj_2_1}
\Re\langle T^\gamma_{\sigma}Z_j,T^\gamma_{(\chi_2-1)\tilde{a}}T^\gamma_{\sigma}Z_j\rangle\leq\e\gamma\|T^\gamma_{\sigma}Z_j\|_0^2+\frac{1}{\e\gamma}\|W^+\|_0^2.
\end{align}
Moreover, since $|\Re\tilde{a}|\geq c\gamma$, one has
\begin{align}\label{variableZj_2_2}
|\Re\langle T^\gamma_{\sigma}Z_j,T^\gamma_{\tilde{a}}T^\gamma_{\sigma}Z_j\rangle|\geq c\gamma\|T^\gamma_{\sigma}Z_j\|_0^2.
\end{align}
Combining \eqref{variableZj_2}, \eqref{variableZj_2_1} and \eqref{variableZj_2_2} gives us the following:
\begin{equation}\label{variableZjmain1}
\begin{split}
 \gamma\|T^\gamma_{\sigma}Z_j\|_0^2\leq &\frac{1}{\e\gamma}\|Z_j\|_0^2+C\|T^\gamma_{\sigma}Z_j\|_0^2+\frac{1}{\e}\|Z_j\|_0^2+\frac{1}{\e\gamma}\left(\|T^\gamma_{\Lambda}T^\gamma_{\sigma}Z_3\|_0^2+\|Z_3\|_{1,\gamma}^2\right)\\
&+\sum_i\frac{1}{\e\gamma}\left(\|T^\gamma_{\sigma}Z_i\|_0^2+\|Z_i\|_0^2\right)+\frac{1}{\e\gamma}\left(\|T^\gamma_{r}W^+\|_{1,\gamma}^2+\|W^+\|_0^2+\|F\|_{1,\gamma}^2\right),
\end{split}
\end{equation}
for $j=1,4,6$. Similarly to the outgoing mode $Z_3$, we apply another symmetrizer $T^\gamma_{\Lambda}$ to \eqref{variableZj_1} and obtain
\begin{equation}\label{E516}
\begin{split}
&\Re\langle Z_j,T^\gamma_{\Lambda}T^\gamma_{\chi_2\Theta_0\times(\gamma+\ti\sigma)}Z_j\rangle+\Re\langle Z_j,T^\gamma_{\Lambda}T^\gamma_{\Theta_1}Z_3\rangle+\Sigma\Re\langle Z_j,T^\gamma_{\Lambda}T^\gamma_{\Theta_0}Z_i\rangle\\
&\qquad\qquad +\Re\langle Z_j,T^\gamma_{\Lambda}T^\gamma_{r}W^+\rangle+\Re\langle Z_j,T^\gamma_{\Lambda}R_{-1}W^+\rangle =\Re\langle Z_j,T^\gamma_{\Lambda}F\rangle.
\end{split}
\end{equation}
In the above equality \eqref{E516}, the first term can be written as
\begin{align*}
\Re\langle Z_j,T^\gamma_{\Lambda}T^\gamma_{\chi_2\Theta_0\times(\gamma+\ti\sigma)}Z_j\rangle & =\Re\langle T^\gamma_{\Lambda^{\frac{1}{2}}}Z_j,T^\gamma_{\Lambda^{\frac{1}{2}}}T^\gamma_{\chi_2\Theta_0\times(\gamma+\ti\sigma)}Z_j\rangle+\Re\langle Z_j,R_0T^\gamma_{\chi_2\Theta_0\times(\gamma+\ti\sigma)}Z_j\rangle\\
&=\Re\langle T^\gamma_{\Lambda^{\frac{1}{2}}}Z_j,T^\gamma_{\chi_2\Theta_0\times(\gamma+\ti\sigma)}T^\gamma_{\Lambda^{\frac{1}{2}}}Z_j\rangle+\Re\langle T^\gamma_{\Lambda^{\frac{1}{2}}}Z_j,R_{\frac{1}{2}}Z_j\rangle \\
&\quad +\Re\langle Z_j,R_0T^\gamma_{\chi_2\Theta_0\times(\gamma+\ti\sigma)}Z_j\rangle.
\end{align*}
For $\Re\langle T^\gamma_{\Lambda^{\frac{1}{2}}}Z_j,T^\gamma_{\chi_2\Theta_0\times(\gamma+\ti\sigma)}T^\gamma_{\Lambda^{\frac{1}{2}}}Z_j\rangle$, similarly to the above, we have
\begin{align*}
&\Re\langle T^\gamma_{\Lambda^{\frac{1}{2}}}Z_j,T^\gamma_{\chi_2\Theta_0\times(\gamma+\ti\sigma)}T^\gamma_{\Lambda^{\frac{1}{2}}}Z_j\rangle=\Re\langle T^\gamma_{\Lambda^{\frac{1}{2}}}Z_j,T^\gamma_{\tilde{a}}T^\gamma_{\Lambda^{\frac{1}{2}}}Z_j\rangle+\Re\langle T^\gamma_{\Lambda^{\frac{1}{2}}}Z_j,T^\gamma_{(\chi_2-1)\tilde{a}}T^\gamma_{\Lambda^{\frac{1}{2}}}Z_j\rangle,
\end{align*}
and 
\begin{align*}
&\Re\langle T^\gamma_{\Lambda^{\frac{1}{2}}}Z_j,T^\gamma_{(\chi_2-1)\tilde{a}}T^\gamma_{\Lambda^{\frac{1}{2}}}Z_j\rangle\leq C\|Z_j\|_{\frac{1}{2},\gamma}^2+\e\gamma\|Z_j\|_{\frac{1}{2},\gamma}^2+\frac{1}{\e\gamma}\|W^+\|_{-1,\gamma}^2,\\
&|\Re\langle T^\gamma_{\Lambda^{\frac{1}{2}}}Z_j,T^\gamma_{\tilde{a}}T^\gamma_{\Lambda^{\frac{1}{2}}}Z_j\rangle|\geq c\gamma\|Z_j\|_{\frac{1}{2},\gamma}^2.
\end{align*}
The estimates of the other terms in  \eqref{E516} are straightforward as the following:
\begin{align*}
&\Re\langle Z_j,R_0T^\gamma_{\chi_2\Theta_0\times(\gamma+\ti\sigma)}Z_j\rangle\leq C\|Z_j\|_{\frac{1}{2},\gamma}^2,\\
&\Re\langle Z_j,T^\gamma_{\Lambda}T^\gamma_{\Theta_1}Z_3\rangle\leq\e\gamma\|Z_j\|_{\frac{1}{2},\gamma}^2+\frac{1}{\e\gamma}\|Z_3\|_{\frac{3}{2},\gamma}^2,\\
&\Re\langle Z_j,T^\gamma_{\Lambda}T^\gamma_{\Theta_0}Z_i\rangle\leq\e\gamma\|Z_j\|_{\frac{1}{2},\gamma}^2+\frac{1}{\e\gamma}\|Z_i\|_{\frac{1}{2},\gamma}^2,\\
&\Re\langle Z_j,T^\gamma_{\Lambda}T^\gamma_{r}W^+\rangle\leq\e\gamma^2\|Z_j\|_0^2+\frac{1}{\e\gamma^2}\|T^\gamma_{r}W^+\|_{1,\gamma}^2,\\
&\Re\langle Z_j,T^\gamma_{\Lambda}R_{-1}W^+\rangle\leq\e\gamma^2\|Z_j\|_0^2+\frac{1}{\e\gamma^2}\|W^+\|_0^2,\\
&\Re\langle Z_j,T^\gamma_{\Lambda}F\rangle\leq\e\gamma^2\|Z_j\|_0^2+\frac{1}{\e\gamma^2}\|F\|_{1,\gamma}^2.
\end{align*}
Thus these estimates together lead to
\begin{equation}\label{variableZjmain2}
\begin{split}
\gamma\|Z_j\|_{\frac{1}{2},\gamma}^2 & \leq C\|Z_j\|_{\frac{1}{2},\gamma}^2+\frac{1}{\e\gamma}\|Z_3\|_{\frac{3}{2},\gamma}^2\\
&\quad+\frac{1}{\e\gamma}\|Z_i\|_{\frac{1}{2},\gamma}^2+\frac{1}{\e\gamma^2}(\|T^\gamma_{r}W^+\|_{1,\gamma}^2+\|W^+\|_0^2+\|F\|_{1,\gamma}^2),
\end{split}
\end{equation}
for $j=1,4,6$, by taking $\e$ small enough.

For the fifth and seventh equations in \eqref{de}, we have
\begin{align*}
T^\gamma_{\tau+\ti v^r\eta}Z_j+T^\gamma_{\Theta_1}Z_3+\Sigma T^\gamma_{\Theta_0}Z_i+T^\gamma_{r}W^++R_{-1}W^+=R_0F,
\end{align*} 
where $j=5,7$. We just follow the similar argument for $Z_j$ with $j=1,4,6$ to obtain
\begin{align}
\gamma\|T^\gamma_{\sigma}Z_j\|_0^2 & \leq\frac{1}{\e\gamma}\|Z_j\|_0^2+C\|T^\gamma_{\sigma}Z_j\|_0^2+\frac{1}{\e}\|Z_j\|_0^2+\frac{1}{\e\gamma}\left(\|T^\gamma_{\Lambda}T^\gamma_{\sigma}Z_3\|_0^2+\|Z_3\|_{1,\gamma}^2\right) \nonumber \\
&\quad + \sum_i\frac{1}{\e\gamma}\left(\|T^\gamma_{\sigma}Z_i\|_0^2+\|Z_i\|_0^2\right)+\frac{1}{\e\gamma}\left(\|T^\gamma_{r}W^+\|_{1,\gamma}^2+\|W^+\|_0^2+\|F\|_{1,\gamma}^2\right), \label{variableZjmain3}\\
\gamma\|Z_j\|_{\frac{1}{2},\gamma}^2 & \leq C\|Z_j\|_{\frac{1}{2},\gamma}^2+\frac{1}{\e\gamma}\|Z_3\|_{\frac{3}{2},\gamma}^2 \nonumber\\
&\quad+\frac{1}{\e\gamma}\|Z_i\|_{\frac{1}{2},\gamma}^2+\frac{1}{\e\gamma^2}(\|T^\gamma_{r}W^+\|_{1,\gamma}^2+\|W^+\|_0^2+\|F\|_{1,\gamma}^2), \label{variableZjmain4}
\end{align}
for $j=5,7$. Note that the estimates for $Z_j$, $j=1,4,5,6,7$, are the same. 

Now we consider the equation of incoming mode $Z_2$ in \eqref{de}
\begin{equation}\label{variableZ2}
\begin{split}
\p_2Z_2 & =T^\gamma_{\o^r}Z_2+T^\gamma_{\Theta_1}Z_1+T^\gamma_{\Theta_1}Z_5\\
&\quad+T^\gamma_{\Theta_1}Z_7+T^\gamma_{\Theta_0}Z_2+\Sigma_{i\neq2,3}T^\gamma_{\Theta_0}Z_i+T^\gamma_{r}W^++R_{-1}W^++F.
\end{split}
\end{equation}
We apply the two symmetrizers $(T^\gamma_{\sigma})^*T^\gamma_{\frac{1}{\Lambda}}T^\gamma_{\sigma}$ and   $1$ to \eqref{variableZ2}.
For $(T^\gamma_{\sigma})^*T^\gamma_{\frac{1}{\Lambda}}T^\gamma_{\sigma}$, we have
\begin{align*}
\Re\langle T^\gamma_{\sigma}Z_2,T^\gamma_{\frac{1}{\Lambda}}T^\gamma_{\sigma}\p_2Z_2\rangle & = \Re\langle T^\gamma_{\sigma}Z_2,T^\gamma_{\frac{1}{\Lambda}}T^\gamma_{\sigma}T^\gamma_{\o^r+\ti\varpi^r}Z_2\rangle+\Re\langle T^\gamma_{\sigma}Z_2,T^\gamma_{\frac{1}{\Lambda}}T^\gamma_{\sigma}T^\gamma_{\Theta_1}Z_j\rangle\\
&\quad +\Re\langle T^\gamma_{\sigma}Z_2,T^\gamma_{\frac{1}{\Lambda}}T^\gamma_{\sigma}T^\gamma_{\Theta_0}Z_2\rangle+\Sigma_{i\neq2,3}\Re\langle T^\gamma_{\sigma}Z_2,T^\gamma_{\frac{1}{\Lambda}}T^\gamma_{\sigma}T^\gamma_{\Theta_0}Z_i\rangle \\
&\quad +\Re\langle T^\gamma_{\sigma}Z_2,T^\gamma_{\frac{1}{\Lambda}}T^\gamma_{\sigma}T^\gamma_{r}W^+\rangle +\Re\langle T^\gamma_{\sigma}Z_2,T^\gamma_{\frac{1}{\Lambda}}T^\gamma_{\sigma}R_{-1}W^+\rangle \\
&\quad +\Re\langle T^\gamma_{\sigma}Z_2,T^\gamma_{\frac{1}{\Lambda}}T^\gamma_{\sigma}F\rangle.
\end{align*}
Similarly for the outgoing modes, we have
\begin{align*}
 \p_2\Re\langle T^\gamma_{\sigma}Z_2,T^\gamma_{\frac{1}{\Lambda}}T^\gamma_{\sigma}Z_2\rangle & =\Re\langle T^\gamma_{\p_2\sigma}Z_2,T^\gamma_{\frac{1}{\Lambda}}T^\gamma_{\sigma}Z_2\rangle+\Re\langle T^\gamma_{\sigma}Z_2,T^\gamma_{\frac{1}{\Lambda}}T^\gamma_{\p_2\sigma}Z_2\rangle\\
& +\Re\langle T^\gamma_{\sigma}\p_2Z_2,T^\gamma_{\frac{1}{\Lambda}}T^\gamma_{\sigma}Z_2\rangle+\Re\langle T^\gamma_{\sigma}Z_2,T^\gamma_{\frac{1}{\Lambda}}T^\gamma_{\sigma}\p_2Z_2\rangle,\\
 \Re\langle T^\gamma_{\sigma}Z_2,T^\gamma_{\frac{1}{\Lambda}}T^\gamma_{\sigma}T^\gamma_{\o^r+\ti\varpi^r}Z_2\rangle & =\Re\langle T^\gamma_{\sigma}Z_2,T^\gamma_{\frac{\o^r+\ti\varpi^r}{\Lambda}}T^\gamma_{\sigma}Z_2\rangle+\Re\langle T^\gamma_{\sigma}Z_2,R_0Z_2\rangle,
\end{align*}
with the following estimates:
\begin{align*}
&\Re\langle T^\gamma_{\p_2\sigma}Z_2,T^\gamma_{\frac{1}{\Lambda}}T^\gamma_{\sigma}Z_2\rangle\leq \e\|T^\gamma_{\sigma}Z_2\|_0^2+\frac{1}{\e}\|Z_2\|_0^2,\\
&\Re\langle T^\gamma_{\sigma}Z_2,T^\gamma_{\frac{1}{\Lambda}}T^\gamma_{\p_2\sigma}Z_2\rangle\leq \e\|T^\gamma_{\sigma}Z_2\|_0^2+\frac{1}{\e}\|Z_2\|_0^2,\\
&\Re\langle T^\gamma_{\sigma}\p_2Z_2,T^\gamma_{\frac{1}{\Lambda}}T^\gamma_{\sigma}Z_2\rangle=\Re\langle T^\gamma_{\frac{1}{\Lambda}}T^\gamma_{\sigma}\p_2Z_2,T^\gamma_{\sigma}Z_2\rangle+\Re\langle T^\gamma_{\sigma}\p_2Z_2,R_{-2}T^\gamma_{\sigma}Z_2\rangle,\\
&\Re\langle T^\gamma_{\sigma}Z_2,T^\gamma_{\frac{\o^r+\ti\varpi^r}{\Lambda}}T^\gamma_{\sigma}Z_2\rangle\leq -c\|T^\gamma_{\sigma}Z_2\|_0^2,\\
&\Re\langle T^\gamma_{\sigma}Z_2,R_0Z_2\rangle\leq \e\|T^\gamma_{\sigma}Z_2\|_0^2+\frac{1}{\e}\|Z_2\|_0^2,\\
&\Re\langle T^\gamma_{\sigma}Z_2,T^\gamma_{\frac{1}{\Lambda}}T^\gamma_{\sigma}T^\gamma_{\Theta_1}Z_j\rangle=\Re\langle T^\gamma_{\sigma}Z_2,T^\gamma_{\frac{1}{\Lambda}}T^\gamma_{\Theta_1}T^\gamma_{\sigma}Z_j\rangle+\Re\langle T^\gamma_{\sigma}Z_2,R_0Z_j\rangle\\
&\qquad\qquad \qquad\qquad \leq \e\|T^\gamma_{\sigma}Z_2\|_0^2+\frac{1}{\e}\|T^\gamma_{\sigma}Z_j\|_0^2+\frac{1}{\e}\|Z_j\|_0^2,\\
&\Re\langle T^\gamma_{\sigma}Z_2,T^\gamma_{\frac{1}{\Lambda}}T^\gamma_{\sigma}T^\gamma_{\Theta_0}Z_2\rangle=\Re\langle T^\gamma_{\sigma}Z_2,T^\gamma_{\frac{1}{\Lambda}}T^\gamma_{\Theta_0}T^\gamma_{\sigma}Z_2\rangle+\Re\langle T^\gamma_{\sigma}Z_2,R_{-1}Z_2\rangle\\
&\qquad\qquad \qquad\qquad  \leq \e\|T^\gamma_{\sigma}Z_2\|_0^2+\frac{1}{\e}\|T^\gamma_{\sigma}Z_2\|_{-1,\gamma}^2+\frac{1}{\e}\|Z_2\|_{-1,\gamma}^2,\\
&\Re\langle T^\gamma_{\sigma}Z_2,T^\gamma_{\frac{1}{\Lambda}}T^\gamma_{\sigma}T^\gamma_{\Theta_0}Z_i\rangle\leq \e\|T^\gamma_{\sigma}Z_2\|_0^2+\frac{1}{\e}\|T^\gamma_{\sigma}Z_i\|_{-1,\gamma}^2+\frac{1}{\e}\|Z_i\|_{-1,\gamma}^2,\\
&\Re\langle T^\gamma_{\sigma}Z_2,T^\gamma_{\frac{1}{\Lambda}}T^\gamma_{\sigma}T^\gamma_{r}W^+\rangle\leq \e\|T^\gamma_{\sigma}Z_2\|_0^2+\frac{1}{\e}\|T^\gamma_{r}W^+\|_{0}^2,\\
&\Re\langle T^\gamma_{\sigma}Z_2,T^\gamma_{\frac{1}{\Lambda}}T^\gamma_{\sigma}R_{-1}W^+\rangle\leq \e\|T^\gamma_{\sigma}Z_2\|_0^2+\frac{1}{\e}\|W^+\|_{-1,\gamma}^2,\\
&\Re\langle T^\gamma_{\sigma}Z_2,T^\gamma_{\frac{1}{\Lambda}}T^\gamma_{\sigma}F\rangle\leq \e\|T^\gamma_{\sigma}Z_2\|_0^2+\frac{1}{\e}\|F\|_{0}^2.
\end{align*}
 These estimates imply
\begin{equation}\label{variableZ2main1}
\begin{split}
 \vertiii{T^\gamma_{\sigma}Z_2}_0^2 & \leq\Re\langle T^\gamma_{\sigma}Z_2,T^\gamma_{\frac{1}{\Lambda}}T^\gamma_{\sigma}Z_2\rangle|_{x_2=0}+\frac{1}{\e}\vertiii{Z_2}_0^2+\frac{1}{\e}\left(\vertiii{T^\gamma_{\sigma}Z_j}_0^2+\vertiii{Z_j}_0^2\right)\\
&\quad +\frac{1}{\e}\left(\vertiii{T^\gamma_{\sigma}Z_2}_{-1,\gamma}^2+\vertiii{Z_2}_{-1,\gamma}^2\right)+\sum_{i\neq2,3}\frac{1}{\e}\left(\vertiii{T^\gamma_{\sigma}Z_i}_{-1,\gamma}^2+\vertiii{Z_i}_{-1,\gamma}^2\right)\\
&\quad +\frac{1}{\e}\left(\vertiii{T^\gamma_{r}W^+}_{0}^2+\vertiii{W^+}_{-1,\gamma}^2+\vertiii{F}_{0}^2\right).
\end{split}
\end{equation}
Next applying the symmetrizer $1$ to obtain
\begin{align*}
\p_2\Re\langle Z_2, Z_2\rangle & =2\Re\langle Z_2, \p_2Z_2\rangle\\
& =2\Re\langle Z_2, T^\gamma_{\o^r+\ti\varpi^r}Z_2\rangle+2\Re\langle Z_2, T^\gamma_{\Theta_1}Z_j\rangle+2\Re\langle Z_2, T^\gamma_{\Theta_0}Z_2\rangle\\
&\quad +\Sigma_{i\neq2,3}2\Re\langle Z_2, T^\gamma_{\Theta_0}Z_i\rangle+2\Re\langle Z_2, T^\gamma_{r}W^+\rangle+2\Re\langle Z_2, R_{-1}W^+\rangle+2\Re\langle Z_2, F\rangle.
\end{align*}
Noting that
\begin{align*}
2\Re\langle Z_2, T^\gamma_{\o^r+\ti\varpi^r}Z_2\rangle & =2\Re\langle Z_2,(T^\gamma_{\Lambda^{\frac{1}{2}}})^* T^\gamma_{\frac{\o^r+\ti\varpi^r}{\Lambda^{\frac{1}{2}}}}Z_2\rangle+2\Re\langle Z_2, R_0Z_2\rangle\\
&=2\Re\langle T^\gamma_{\Lambda^{\frac{1}{2}}}Z_2, T^\gamma_{\frac{\o^r+\ti\varpi^r}{\Lambda}}T^\gamma_{\Lambda^{\frac{1}{2}}}Z_2\rangle+2\Re\langle T^\gamma_{\Lambda^{\frac{1}{2}}}Z_2, R_{-\frac{1}{2}}Z_2\rangle+2\Re\langle Z_2, R_0Z_2\rangle,
\end{align*}
and from the   estimates
\begin{align*}
&\Re\langle T^\gamma_{\Lambda^{\frac{1}{2}}}Z_2, T^\gamma_{\frac{\o^r+\ti\varpi^r}{\Lambda}}T^\gamma_{\Lambda^{\frac{1}{2}}}Z_2\rangle\leq-c\|Z_2\|_{\frac{1}{2},\gamma}^2,\\
&\Re\langle T^\gamma_{\Lambda^{\frac{1}{2}}}Z_2, R_{-\frac{1}{2}}Z_2\rangle\leq \e\|Z_2\|_{\frac{1}{2},\gamma}^2+\frac{1}{\e}\|Z_2\|_{-\frac{1}{2},\gamma}^2,\\
&\Re\langle Z_2, R_0Z_2\rangle\leq C\|Z_2\|_0^2,
\end{align*}
as well as the straightforward estimates of the other terms,  one has
\begin{equation}\label{variableZ2main2}
\begin{split}
 \vertiii{Z_2}_{\frac{1}{2},\gamma}^2\leq &\|Z_2|_{x_2=0}\|_0^2+(C+\frac{1}{\e})\vertiii{Z_2}_0^2+\frac{1}{\e}\vertiii{Z_j}_{\frac{1}{2},\gamma}^2+\Sigma_{i\neq2,3}\frac{1}{\e}\vertiii{Z_i}_{-\frac{1}{2},\gamma}^2\\
& +\frac{1}{\e}\left(\vertiii{T^\gamma_{r}W^+}_{-\frac{1}{2},\gamma}^2+\vertiii{W^+}_{-\frac{3}{2},\gamma}^2+\vertiii{F}_{-\frac{1}{2},\gamma}^2\right).
\end{split}
\end{equation}

We consider \eqref{variableZ3main1}, \eqref{variableZ3main2}, \eqref{variableZjmain1}, \eqref{variableZjmain2}, \eqref{variableZjmain3}, \eqref{variableZjmain4}, \eqref{variableZ2main1} and \eqref{variableZ2main2} and divide them by appropriate powers of $\gamma$ to obtain
\begin{align*}
\frac{1}{\gamma}\vertiii{T^\gamma_{\Lambda}T^\gamma_{\sigma}Z_3}_0^2 & +\frac{1}{\gamma}\|T^\gamma_{\Lambda^\frac{1}{2}}T^\gamma_{\sigma}Z_3|_{x_2=0}\|_0^2 \\
& \lesssim \frac{1}{\gamma}\|T^\gamma_{\sigma}Z_3|_{x_2=0}\|_0^2+\left(C+\frac{1}{\e}\right)\frac{1}{\gamma}\vertiii{Z_3}_{1,\gamma}^2\\
&\quad+\sum_{i\neq2,3}\frac{1}{\e\gamma}\left(\vertiii{T^\gamma_{\sigma}Z_i}_0^2+\vertiii{Z_i}_0^2\right) \\
&\quad +\frac{1}{\e\gamma}\left(\vertiii{T_rW^+}_{1,\gamma}^2+\vertiii{W^+}_0^2+\vertiii{F}_{1,\gamma}^2\right),\\
\vertiii{Z_3}_{\frac{3}{2},\gamma}^2+\|Z_3|_{x_2=0}\|_{1,\gamma}^2 & \lesssim C\vertiii{Z_3}_{1,\gamma}^2+\frac{1}{\e}\vertiii{Z_3}_{\frac{1}{2},\gamma}^2+\sum_{i\neq2,3}(\frac{1}{\e}\vertiii{Z_i}_{\frac{1}{2},\gamma}^2+\frac{1}{\e\gamma}\vertiii{Z_i}_0^2)\\
&\quad+\frac{1}{\e\gamma}\left(\vertiii{T_rW^+}_{1,\gamma}^2+\vertiii{W^+}_0^2+\vertiii{F}_{1,\gamma}^2\right),\\
\gamma\vertiii{T^\gamma_{\sigma}Z_j}_0^2 & \leq\frac{1}{\e\gamma}\vertiii{Z_j}_0^2+C\vertiii{T^\gamma_{\sigma}Z_j}_0^2+\frac{1}{\e}\vertiii{Z_j}_0^2 \\
&\quad +\frac{1}{\e\gamma}\left(\vertiii{T^\gamma_{\Lambda}T^\gamma_{\sigma}Z_3}_0^2+\vertiii{Z_3}_{1,\gamma}^2\right)+\sum_{i\ne 2,3}\frac{1}{\e\gamma}\left(\vertiii{T^\gamma_{\sigma}Z_i}_0^2+\vertiii{Z_i}_0^2\right) \\
& \quad +\frac{1}{\e\gamma}\left(\vertiii{T^\gamma_{r}W^+}_{1,\gamma}^2+\vertiii{W^+}_0^2+\vertiii{F}_{1,\gamma}^2\right),\\
\gamma^2\vertiii{Z_j}_{\frac{1}{2},\gamma}^2 & \leq C\gamma\vertiii{Z_j}_{\frac{1}{2},\gamma}^2+\frac{1}{\e}\vertiii{Z_3}_{\frac{3}{2},\gamma}^2+\frac{1}{\e}\vertiii{Z_i}_{\frac{1}{2},\gamma}^2\\
&\quad+\frac{1}{\e\gamma}\left(\vertiii{T^\gamma_{r}W^+}_{1,\gamma}^2+\vertiii{W^+}_0^2+\vertiii{F}_{1,\gamma}^2\right),\\
\gamma\vertiii{T^\gamma_{\sigma}Z_2}_0^2 & \leq \gamma\Re\langle T^\gamma_{\sigma}Z_2,T^\gamma_{\frac{1}{\Lambda}}T^\gamma_{\sigma}Z_2\rangle|_{x_2=0}+\frac{\gamma}{\e}\vertiii{Z_2}_0^2+\frac{\gamma}{\e}\left(\vertiii{T^\gamma_{\sigma}Z_j}_0^2+\vertiii{Z_j}_0^2\right)\\
&\quad+\frac{\gamma}{\e}\left(\vertiii{T^\gamma_{\sigma}Z_2}_{-1,\gamma}^2+\vertiii{Z_2}_{-1,\gamma}^2\right)+\sum_{i\neq2,3}\frac{\gamma}{\e}\left(\vertiii{T^\gamma_{\sigma}Z_i}_{-1,\gamma}^2+\vertiii{Z_i}_{-1,\gamma}^2\right)\\
&\quad+\frac{\gamma}{\e}\left(\vertiii{T^\gamma_{r}W^+}_{0}^2+\vertiii{W^+}_{-1,\gamma}^2+\vertiii{F}_{0}^2\right),\\
\gamma^2\vertiii{Z_2}_{\frac{1}{2},\gamma}^2 & \leq\gamma^2\|Z_2|_{x_2=0}\|_0^2+(C+\frac{1}{\e})\gamma^2\vertiii{Z_2}_0^2+\frac{\gamma^2}{\e}\vertiii{Z_j}_{\frac{1}{2},\gamma}^2+\sum_{i\neq2,3}\frac{\gamma^2}{\e}\vertiii{Z_i}_{-\frac{1}{2},\gamma}^2\\
&\quad+\frac{\gamma^2}{\e}\left(\vertiii{T^\gamma_{r}W^+}_{-\frac{1}{2},\gamma}^2+\vertiii{W^+}_{-\frac{3}{2},\gamma}^2+\vertiii{F}_{-\frac{1}{2},\gamma}^2\right).
\end{align*}
Then adding all the  above inequalities together and taking $\gamma$ large enough lead to
\begin{equation}\label{polerootW+}
{\small
\begin{split}
\hspace{-.1in}\frac{1}{\gamma} &\vertiii{T^\gamma_{\Lambda}T^\gamma_{\sigma}Z_3}_0^2 +\vertiii{Z_3}_{\frac{3}{2},\gamma}^2+\gamma\vertiii{T^\gamma_{\sigma}Z_j}_0^2+\gamma^2\vertiii{Z_j}_{\frac{1}{2},\gamma}^2+\gamma\vertiii{T^\gamma_{\sigma}Z_2}_0^2+\gamma^2\vertiii{Z_2}_{\frac{1}{2},\gamma}^2\\
& \qquad +\frac{1}{\gamma} \|T^\gamma_{\Lambda^\frac{1}{2}}T^\gamma_{\sigma}Z_3|_{x_2=0}\|_0^2+\|Z_3|_{x_2=0}\|_{1,\gamma}^2 \\
& \lesssim\gamma\Re\langle T^\gamma_{\sigma}Z_2,T^\gamma_{\frac{1}{\Lambda}}T^\gamma_{\sigma}Z_2\rangle|_{x_2=0}+\gamma^2\|Z_2|_{x_2=0}\|_0^2+\frac{1}{\gamma}\vertiii{Z_3}_{1,\gamma}^2 \\
& \quad+\Sigma\left(\vertiii{Z_i}_{\frac{1}{2},\gamma}^2+\frac{1}{\gamma}\vertiii{T^\gamma_{\sigma}Z_i}_0^2\right)+\frac{1}{\gamma}\left(\vertiii{T_rW^+}_{1,\gamma}^2+\vertiii{W^+}_0^2+\vertiii{F}_{1,\gamma}^2\right)\\
&\lesssim\gamma\Re\langle T^\gamma_{\sigma}Z_2,T^\gamma_{\frac{1}{\Lambda}}T^\gamma_{\sigma}Z_2\rangle|_{x_2=0}+\gamma^2\|Z_2|_{x_2=0}\|_0^2 +\frac{1}{\gamma}\left(\vertiii{T_rW^+}_{1,\gamma}^2+\vertiii{W^+}_0^2+\vertiii{F}_{1,\gamma}^2\right).
\end{split}}
\end{equation}

In our setting, since $v^r\neq v^l$, then the point where $\tau=-\ti v^r\eta$ is not a pole of the differential equation for $W^-$ in \eqref{paralinearizationsystem} and thus it can be treated in the same way as in Case 2, with  the following estimate:
\begin{equation}\label{polerootW-}
{\small
\begin{split}
&\frac{1}{\gamma}\vertiii{T^\gamma_{\Lambda}T^\gamma_{\sigma}Z_{10}}_0^2+\vertiii{Z_{10}}_{\frac{3}{2},\gamma}^2+\vertiii{T^\gamma_{\sigma}Z_j}_{\frac{1}{2},\gamma}^2+\gamma\vertiii{Z_j}_{1,\gamma}^2+\gamma\vertiii{T^\gamma_{\sigma}Z_9}_0^2+\gamma^2\vertiii{Z_9}_{\frac{1}{2},\gamma}^2\\
&\quad +\frac{1}{\gamma}\|T^\gamma_{\Lambda^\frac{1}{2}}T^\gamma_{\sigma}Z_{10}|_{x_2=0}\|_0^2+\|Z_{10}|_{x_2=0}\|_{1,\gamma}\\
&\lesssim\gamma\Re\langle T^\gamma_{\sigma}Z_9,T^\gamma_{\frac{1}{\Lambda}}T^\gamma_{\sigma}Z_9\rangle|_{x_2=0}+\gamma^2\|Z_9|_{x_2=0}\|_0^2+\frac{1}{\gamma}\left(\vertiii{T_rW^+}_{1,\gamma}^2+\vertiii{W^+}_0^2+\vertiii{F}_{1,\gamma}^2\right),
\end{split}}
\end{equation}
where $Z^-=(Z_9,Z_{10},...,Z_{14})^\top:=T^\gamma_{\chi_1Q^l}T^\gamma_{\chi_{p_1}}W^-$ and $Q^l$ is the transformation matrix for $W^-$ which is defined in a similar way to $Q^r$. In the above, $Z_3$ and $Z_{10}$ are the outgoing modes of the system and $Z_2$ and $Z_9$ are the incoming modes, thus we denote
\begin{align*}
&Z_{in}=(Z_2,Z_9)^\top,\text{ and }Z_{out}=(Z_3,Z_{10})^\top.
\end{align*}
Finally, we need to use the boundary conditions in \eqref{paralinearizationsystem} to estimate $\|Z_2|_{x_2=0}\|_0^2$, $\|Z_9|_{x_2=0}\|_0^2$, $\gamma\Re\langle T^\gamma_{\sigma}Z_2,T^\gamma_{\frac{1}{\Lambda}}T^\gamma_{\sigma}Z_2\rangle|_{x_2=0}$ and $\gamma\Re\langle T^\gamma_{\sigma}Z_9,T^\gamma_{\frac{1}{\Lambda}}T^\gamma_{\sigma}Z_9\rangle|_{x_2=0}$. However, it is obvious that 
\begin{align*}
&\gamma\Re\langle T^\gamma_{\sigma}Z_2,T^\gamma_{\frac{1}{\Lambda}}T^\gamma_{\sigma}Z_2\rangle|_{x_2=0}\lesssim \|T^\gamma_{\sigma}Z_2|_{x_2=0}\|_0^2,\\
&\gamma\Re\langle T^\gamma_{\sigma}Z_9,T^\gamma_{\frac{1}{\Lambda}}T^\gamma_{\sigma}Z_9\rangle|_{x_2=0}\lesssim \|T^\gamma_{\sigma}Z_9|_{x_2=0}\|_0^2,
\end{align*}
thus we only need to estimate $T^\gamma_{\sigma}Z_{in}|_{x_2=0}$ and $Z_{in}|_{x_2=0}$, which is the same situation as in \cite{coulombel2004stability}. So by the same method in \cite{coulombel2004stability} and utilizing the boundary conditions in \eqref{paralinearizationsystem}, we obtain
\begin{align}\label{polerootZin}
\gamma^2\|Z_{in}|_{x_2=0}\|_0^2+\|T^\gamma_{\sigma}Z_{in}|_{x_2=0}\|_0^2\lesssim \|G\|_{1,\gamma}^2+\|Z_{out}|_{x_2=0}\|_{1,\gamma}^2+\|W^\nc|_{x_2=0}\|_0^2.
\end{align}
Moreover, from \eqref{variableZ3main2}, we have
\begin{align}\label{polerootZout}
\|Z_{out}|_{x_2=0}\|_{1,\gamma}^2\lesssim\Sigma_{i\neq2,3}\vertiii{Z_i}_{\frac{1}{2},\gamma}^2+\frac{1}{\gamma}\left(\vertiii{T_rW^+}_{1,\gamma}^2+\vertiii{W^+}_0^2+\vertiii{F}_{1,\gamma}^2\right).
\end{align}
Adding \eqref{polerootW+}, \eqref{polerootW-}, \eqref{polerootZin} and \eqref{polerootZout} together yields the following estimate:
\begin{equation}\label{polerootestimate}
\begin{split}
&\frac{1}{\gamma}\vertiii{T^\gamma_{\Lambda}T^\gamma_{\sigma}Z_{out}}_0^2+\vertiii{Z_{out}}_{\frac{3}{2},\gamma}^2+\gamma\vertiii{T^\gamma_{\sigma}Z_c}_0^2+\gamma^2\vertiii{Z_c}_{\frac{1}{2},\gamma}^2+\gamma\vertiii{T^\gamma_{\sigma}Z_{in}}_0^2+\gamma^2\vertiii{Z_{in}}_{\frac{1}{2},\gamma}^2\\
&\qquad +\frac{1}{\gamma}\|T^\gamma_{\Lambda^\frac{1}{2}}T^\gamma_{\sigma}Z_{out}|_{x_2=0}\|_0^2+\|Z_{out}|_{x_2=0}\|_{1,\gamma}+\gamma^2\|Z_{in}|_{x_2=0}\|_0^2+\|T^\gamma_{\sigma}Z_{in}|_{x_2=0}\|_0^2\\
&\lesssim\|G\|_{1,\gamma}^2+\|W^\nc|_{x_2=0}\|_0^2+\frac{1}{\gamma}\left(\vertiii{T_rW}_{1,\gamma}^2+\vertiii{W}_0^2+\vertiii{F}_{1,\gamma}^2\right),
\end{split}
\end{equation}
where $Z_c=(Z_1,Z_4,Z_5,Z_6,Z_7,Z_8,Z_{11},Z_{12},Z_{13},Z_{14})^\top$.

\subsection{Case 2 -- points in $\Upsilon_r^{(2)}$}\label{subsec_roots} 

We want to estimate the part of $W^\pm$ corresponding to $\V_r^1$, $\V_r^2$ and $\V_r^3$. Here we will just show the details for the differential equations of $W^+$ in $\V_r^1$. For the other neighborhoods $\V_r^2$ and $\V_r^3$ and unkown $W^-$,  we can obtain the estimates in exactly the same way. Now we consider the cut-off function $\chi_{rt}$ in $\Gamma_k^0$ for any integer $k$, whose support on $\R^3_+\times\Sigma$ is contained in the neighborhood $\V_{r}^1$,  and equal to $1$ in a smaller neighborhood of the curve where $\tau=\ti V_1\eta$.  Denote
\begin{align*}
W^+_{rt}:=T^\gamma_{\chi_{rt}}W^+.
\end{align*}
Similarly, as in the previous case, we have
\begin{align*}
T^\gamma_{\tau A_0^r+\ti\eta A_1^r}W^+_{rt}+T^\gamma_{A_0^rC^r}W^+_{rt}+T^\gamma_rW^++I_2\p_2W^+_{rt}=T^\gamma_{\chi_{rt}}F+R_{-1}W^+,
\end{align*}
where $r$ is in the class of $\Gamma^0_1$, bounded and supported only in the set where $\chi_{rt}\in(0,1)$. Then we take two cut-off functions $\chi_1$ and $\chi_2$ in the class $\Gamma^0_k$ for any integer $k$ and both with support in $\V_{rt}$, and $\chi_1=1$ on the support of $\chi_{rt}$ and $\chi_2=1$ on the support of $\chi_1$. Similarly, as in the previous case, after applying the cut-off symbol to the differential equation, we can find the transform matrices $Q_0^r$ and $Q_{-1}^r$, and symmetrizers $R_0^r$ and $R_{-1}^r$ such that
\begin{align*}
I_2\p_2Z^+=-T^\gamma_{\chi_2\widetilde{A}^{r}}Z^++T^\gamma_{D_0}Z^++T^\gamma_rW^++R_0F+R_{-1}W^+,
\end{align*}
where $\widetilde{A}^{r}$ is the same as \eqref{uppertrianglecase1} in the Case 1, $\chi_1Q_0^r$ and $\chi_1R_0^r$ are invertible symbols in $\Gamma^0_{2}$,  $Q_{-1}^r$ and $R_{-1}^r$ are symbols in $\Gamma^{-1}_{1}$,  and $Z^+$ is defined by
\begin{align*}
Z^+=T^\gamma_{\chi_1({Q_0^r}^{-1}+Q_{-1}^r)}W_{rt}^+,
\end{align*}
just as in the previous case. Then we do the same extension as the previous case for $\chi_2\widetilde{A}^{r}$ and $D_0$, and obtain
\begin{align}\label{der}
I_2\p_2Z^+=-T^\gamma_{\widetilde{D_1}}Z^++T^\gamma_{\widetilde{D_0}}Z^++T^\gamma_rW^++R_0F+R_{-1}W^+,
\end{align}
where, in $\widetilde{D_1}$, $\o^r$ is a symbol in $\Gamma^1_2$ and $\o^r\geq c \Lambda^{1,\gamma}$, and in $\widetilde{D_0}$, $d_{2,3}=d_{3,2}=0$.
We need to consider each differential equation and construct appropriate symmetrizers to obtain the desired estimate. Still we denote $Z^+=(Z_1,Z_2,...,Z_7)^\top$. For the third equation in \eqref{der}, we have
\begin{align*}
\p_2Z_3=T^\gamma_{-\o^r+\ti\varpi^r}Z_3+T^\gamma_{\Theta_0}Z_3+\Sigma_{i\neq2,3}T^\gamma_{\Theta_0}Z_i+T^\gamma_rW^++R_0F+R_{-1}W^+,
\end{align*}
which is the same equation as the previous case. So we can still take $(T^\gamma_{\sigma})^*T^\gamma_{\Lambda}T^\gamma_{\sigma}$ and $(T^\gamma_{\Lambda})^*T^\gamma_{\Lambda}$ as symmetrizers. After the same argument, we   have
\begin{equation*}
\begin{split}
&\vertiii{T^\gamma_{\Lambda}T^\gamma_{\sigma}Z_3}_0^2+\|T^\gamma_{\Lambda^\frac{1}{2}}T^\gamma_{\sigma}Z_3|_{x_2=0}\|_0^2\lesssim \|T^\gamma_{\sigma}Z_3|_{x_2=0}\|_0^2+(C+\frac{1}{\e})\vertiii{Z_3}_{1,\gamma}^2\\
&\qquad+\Sigma_{i\neq2,3}\frac{1}{\e}\left(\vertiii{T^\gamma_{\sigma}Z_i}_0^2+\vertiii{Z_i}_0^2\right)+\frac{1}{\e}\left(\vertiii{T_rW^+}_{1,\gamma}^2+\vertiii{W^+}_0^2+\vertiii{F}_{1,\gamma}^2\right),
\end{split}
\end{equation*}
and 
\begin{equation*}
\begin{split}
&\vertiii{Z_3}_{\frac{3}{2},\gamma}^2+\|Z_3|_{x_2=0}\|_{1,\gamma}^2\lesssim C\vertiii{Z_3}_{1,\gamma}^2+\frac{1}{\e}\vertiii{Z_3}_{\frac{1}{2},\gamma}^2\\
&\qquad\qquad +\Sigma_{i\neq2,3}(\frac{1}{\e}\vertiii{Z_i}_{\frac{1}{2},\gamma}^2+\frac{1}{\e\gamma}\vertiii{Z_i}_0^2)+\frac{1}{\e\gamma}\left(\vertiii{T_rW^+}_{1,\gamma}^2+\vertiii{W^+}_0^2+\vertiii{F}_{1,\gamma}^2\right).
\end{split}
\end{equation*}
Since in this case the root of the Lopatinskii determinant does not coincide with any pole of the differential equation, we can estimate $Z_j$ for $j=1,4,5,6,7$, in the same strategy. 
In particular, we multiply the 1st, 4th, 5th, 6th and 7th equations in \eqref{der} by some appropriate matrix symbol in $\Gamma_1^0$ and obtain
\begin{align}
T^\gamma_{a}Z_j+T^\gamma_{\Theta_1}Z_3+\Sigma T^\gamma_{\Theta_0}Z_i+T^\gamma_{r}W^++R_{-1}W^+=R_0F,
\end{align}  
where $|\Re a|\geq c\Lambda$ on the support of $\chi_2$. Therefore we can do the same extension argument of $a$ as in Case 1 and obtain a new symbol $\tilde{a}$ such that $|\Re \tilde{a}|\geq c\Lambda$. Then in this case, we consider the symmetrizers $(T^\gamma_{\sigma})^*T^\gamma_{\sigma}$ and $T^\gamma_{\Lambda}$, and by a similar argument, we have
\begin{equation*}
\begin{split}
\|T^\gamma_{\sigma}Z_j\|_{\frac{1}{2},\gamma}^2& \leq  C\|Z_j\|_{1,\gamma}^2
+\frac{1}{\e}\|Z_3\|_{\frac{3}{2},\gamma}^2+\Sigma\frac{1}{\e}\|Z_i\|_{\frac{1}{2},\gamma}^2\\
&\quad +\frac{1}{\e}\left(\|T^\gamma_{r}W^+\|_{\frac{1}{2},\gamma}^2+\|W^+\|_{-\frac{1}{2},\gamma}^2+\|F\|_{\frac{1}{2},\gamma}^2\right),
\end{split}
\end{equation*}
and
\begin{align*}
\|Z_j\|_{1,\gamma}^2\leq \frac{1}{\e}\|Z_j\|_0^2+\frac{1}{\e}\|Z_3\|_{1,\gamma}^2+\Sigma\frac{1}{\e}\|Z_i\|_0^2+\frac{1}{\e}\left(\|T^\gamma_{r}W^+\|_0^2+\|W^+\|_{-1,\gamma}^2+\|F\|_0^2\right).
\end{align*}
For the second equation in \eqref{der}, corresponding to the incoming mode of the system, we can still follow the previous argument to obtain
\begin{equation*}
\begin{split}
\vertiii{T^\gamma_{\sigma}Z_2}_0^2& \leq \Re\langle T^\gamma_{\sigma}Z_2,T^\gamma_{\frac{1}{\Lambda}}T^\gamma_{\sigma}Z_2\rangle|_{x_2=0}+\frac{1}{\e}\vertiii{Z_2}_0^2+\frac{1}{\e}\left(\vertiii{T^\gamma_{\sigma}Z_j}_0^2+\vertiii{Z_j}_0^2\right)\\
&\quad+\frac{1}{\e}\left(\vertiii{T^\gamma_{\sigma}Z_2}_{-1,\gamma}^2+\vertiii{Z_2}_{-1,\gamma}^2\right)+\Sigma_{i\neq2,3}\frac{1}{\e}\left(\vertiii{T^\gamma_{\sigma}Z_i}_{-1,\gamma}^2+\vertiii{Z_i}_{-1,\gamma}^2\right)\\
&\quad +\frac{1}{\e}\left(\vertiii{T^\gamma_{r}W^+}_{0}^2+\vertiii{W^+}_{-1,\gamma}^2+\vertiii{F}_{0}^2\right),
\end{split}
\end{equation*}
and 
\begin{equation*}
\begin{split}
\vertiii{Z_2}_{\frac{1}{2},\gamma}^2 & \leq \|Z_2|_{x_2=0}\|_0^2+(C+\frac{1}{\e})\vertiii{Z_2}_0^2+\frac{1}{\e}\vertiii{Z_j}_{\frac{1}{2},\gamma}^2+\Sigma_{i\neq2,3}\frac{1}{\e}\vertiii{Z_i}_{-\frac{1}{2},\gamma}^2\\
&\quad +\frac{1}{\e}\left(\vertiii{T^\gamma_{r}W^+}_{-\frac{1}{2},\gamma}^2+\vertiii{W^+}_{-\frac{3}{2},\gamma}^2+\vertiii{F}_{-\frac{1}{2},\gamma}^2\right).
\end{split}
\end{equation*}
Again, combining all the above estimates, dividing them by suitable powers of $\gamma$ and taking $\gamma$ large enough, we have
\begin{equation}\label{rootW+}
\begin{split}
&\frac{1}{\gamma}\vertiii{T^\gamma_{\Lambda}T^\gamma_{\sigma}Z_3}_0^2+\vertiii{Z_3}_{\frac{3}{2},\gamma}^2+\vertiii{T^\gamma_{\sigma}Z_j}_{\frac{1}{2},\gamma}^2+\gamma\vertiii{Z_j}_{1,\gamma}^2+\gamma\vertiii{T^\gamma_{\sigma}Z_2}_0^2+\gamma^2\vertiii{Z_2}_{\frac{1}{2},\gamma}^2\\
&\qquad +\frac{1}{\gamma}\|T^\gamma_{\Lambda^\frac{1}{2}}T^\gamma_{\sigma}Z_3|_{x_2=0}\|_0^2+\|Z_3|_{x_2=0}\|_{1,\gamma}^2\\
&\lesssim \gamma\Re\langle T^\gamma_{\sigma}Z_2,T^\gamma_{\frac{1}{\Lambda}}T^\gamma_{\sigma}Z_2\rangle|_{x_2=0}
+\gamma^2\|Z_2|_{x_2=0}\|_0^2+\frac{1}{\gamma}\left(\vertiii{T^\gamma_rW^+}_{1,\gamma}^2+\vertiii{W^+}_0^2+\vertiii{F}_{1,\gamma}^2\right).
\end{split}
\end{equation}
Similarly, for $Z^-=(Z_8, Z_9,...,Z_{14})^\top:=T^\gamma_{\chi_1(Q^l_0+Q^l_{-1})}T^\gamma_{\chi_{rt}}W^-$, we have
\begin{equation}\label{rootW-}
\begin{split}
&\frac{1}{\gamma}\vertiii{T^\gamma_{\Lambda}T^\gamma_{\sigma}Z_{10}}_0^2+\vertiii{Z_{10}}_{\frac{3}{2},\gamma}^2+\vertiii{T^\gamma_{\sigma}Z_j}_{\frac{1}{2},\gamma}^2+\gamma\vertiii{Z_j}_{1,\gamma}^2+\gamma\vertiii{T^\gamma_{\sigma}Z_9}_0^2+\gamma^2\vertiii{Z_9}_{\frac{1}{2},\gamma}^2\\
&\quad +\frac{1}{\gamma}\|T^\gamma_{\Lambda^\frac{1}{2}}T^\gamma_{\sigma}Z_{10}|_{x_2=0}\|_0^2+\|Z_{10}|_{x_2=0}\|_{1,\gamma}^2\\
&\lesssim \gamma\Re\langle T^\gamma_{\sigma}Z_9,T^\gamma_{\frac{1}{\Lambda}}T^\gamma_{\sigma}Z_9\rangle|_{x_2=0}
 +\gamma^2\|Z_9|_{x_2=0}\|_0^2+\frac{1}{\gamma}\left(\vertiii{T^\gamma_rW^+}_{1,\gamma}^2+\vertiii{W^+}_0^2+\vertiii{F}_{1,\gamma}^2\right).
\end{split}
\end{equation}

Note that   as in the Case 1, the boundary terms in \eqref{paralinearizationsystem} can be used to estimate $\gamma\Re\langle T^\gamma_{\sigma}Z_2,T^\gamma_{\frac{1}{\Lambda}}T^\gamma_{\sigma}Z_2\rangle|_{x_2=0}$, $\|Z_2|_{x_2=0}\|_0^2$, $\gamma\Re\langle T^\gamma_{\sigma}Z_9,T^\gamma_{\frac{1}{\Lambda}}T^\gamma_{\sigma}Z_9\rangle|_{x_2=0}$ and $\|Z_9|_{x_2=0}\|_0^2$. By the same method in estimating $T^\gamma_{\sigma}Z_{in}|_{x_2=0}$ and $Z_{in}|_{x_2=0}$ and combining \eqref{rootW+} and \eqref{rootW-}, we have
\begin{equation}\label{rootestimate}
\begin{split}
&\frac{1}{\gamma}\vertiii{T^\gamma_{\Lambda}T^\gamma_{\sigma}Z_{out}}_0^2+\vertiii{Z_{out}}_{\frac{3}{2},\gamma}^2+\vertiii{T^\gamma_{\sigma}Z_c}_{\frac{1}{2},\gamma}^2+\gamma\vertiii{Z_c}_{1,\gamma}^2+\gamma\vertiii{T^\gamma_{\sigma}Z_{in}}_0^2+\gamma^2\vertiii{Z_{in}}_{\frac{1}{2},\gamma}^2\\
&\quad +\frac{1}{\gamma}\|T^\gamma_{\Lambda^\frac{1}{2}}T^\gamma_{\sigma}Z_{out}|_{x_2=0}\|_0^2+\|Z_{out}|_{x_2=0}\|_{1,\gamma}+\gamma^2\|Z_{in}|_{x_2=0}\|_0^2+\|T^\gamma_{\sigma}Z_{in}|_{x_2=0}\|_0^2\\
&\lesssim\|G\|_{1,\gamma}^2+\|W^\nc|_{x_2=0}\|_0^2+\frac{1}{\gamma}\left(\vertiii{T_rW}_{1,\gamma}^2+\vertiii{W}_0^2+\vertiii{F}_{1,\gamma}^2\right),
\end{split}
\end{equation}
where $Z_c=(Z_1,Z_4,Z_5,Z_6,Z_7,Z_8,Z_{11},Z_{12},Z_{13},Z_{14})^\top$. 

\subsection{Case 3 -- points in $\Upsilon_p^{(2)}$} 

Here we will discuss the case of a pole which is not a root of the Lopatinskii determinant, that is, the part of $W^\pm$ whose frequencies are in $\V_{p_2}^1$, $\V_{p_2}^2$, $\V_{p_2}^3$ and $\V_{p_2}^4$. Since those frequencies are always from the roots of the Lopatinskii determinant, the boundary conditions in \eqref{paralinearizationsystem} provide a stronger control on the incoming modes on $\V_{p_2}^1$, $\V_{p_2}^2$, $\V_{p_2}^3$ and $\V_{p_2}^4$ than on the neighborhoods in the previous two cases. In this case, we only need to construct one symmetrizer for each equation in the system. We take the neighborhood $\V_{p_2}^1$ as an example. 
This neighborhood contains the curve in the frequency space such that $\tau=-\ti \left(v^{r}+\sqrt{{F_{11}^{r}}^2+{F_{12}^{r}}^2}\right)\eta$ which is a pole of the differential equations for $W^+$ but not for $W^-$.  In fact, the equations of $W^-$ can be reduced to the non-characteristic case. We want to point out that the strategy we discuss here works for both equations of $W^+$ and $W^-$, provided that the neighborhood does not contain any points at which $\o^r=0$ or $\o^l=0$. In the following we will only give an argument for $W^+$. 

First we take cut-off functions $\chi_{p_2}$, $\chi_1$ and $\chi_2$ as in the last pole case, pick the transform matrices $Q_0^r$ and $Q_{-1}^r$ as well as the  symmetrizers $R_0^r$ and $R_{-1}^r$,  and make the appropriate adjustment of $\chi_2\widetilde{A}^{r}$ and $D_0$ to obtain
\begin{align}\label{dep}
I_2\p_2Z^+=-T^\gamma_{\widetilde{D_1}}Z^++T^\gamma_{\widetilde{D_0}}Z^++T^\gamma_rW^++R_0F+R_{-1}W^+,
\end{align}
where
\begin{align*}
Z^+=T^\gamma_{\chi_1({Q_0^r}^{-1}+Q_{-1}^r)}T^\gamma_{\chi_{p_2}}W^+,
\end{align*}
and the symbols in the equations are the same as in the previous two cases. For the third equation:
\begin{align*}
\p_2Z_3=T^\gamma_{-\o^r+\ti\varpi^r}Z_3+T^\gamma_{\Theta_0}Z_3+\Sigma_{i\neq2,3}T^\gamma_{\Theta_0}Z_i+T^\gamma_rW^++R_0F+R_{-1}W^+,
\end{align*}
we consider the symmetrizer $(T^\gamma_{\Lambda})^*T^\gamma_{\Lambda}T^\gamma_{\Lambda}$ and obtain
\begin{align*}
\Re\langle T^\gamma_{\Lambda}T^\gamma_{\Lambda}Z_3, T^\gamma_{\Lambda}\p_2Z_3\rangle & =\Re\langle T^\gamma_{\Lambda}T^\gamma_{\Lambda}Z_3, T^\gamma_{\Lambda}T^\gamma_{-\o^r+\ti\varpi^r}Z_3\rangle+\Re\langle T^\gamma_{\Lambda}T^\gamma_{\Lambda}Z_3, T^\gamma_{\Lambda}T^\gamma_{\Theta_0}Z_3\rangle\\
&\quad +\Re\langle T^\gamma_{\Lambda}T^\gamma_{\Lambda}Z_3, T^\gamma_{\Lambda}T^\gamma_{\Theta_0}Z_i\rangle+\Re\langle T^\gamma_{\Lambda}T^\gamma_{\Lambda}Z_3, T^\gamma_{\Lambda}T^\gamma_rW^+\rangle\\
&\quad +\Re\langle T^\gamma_{\Lambda}T^\gamma_{\Lambda}Z_3, T^\gamma_{\Lambda}R_{-1}W^+\rangle+\Re\langle T^\gamma_{\Lambda}T^\gamma_{\Lambda}Z_3, T^\gamma_{\Lambda}R_0F\rangle.
\end{align*}
Using the similar argument and taking $\e$ small enough, we have
\begin{align*}
\vertiii{Z_3}_{2,\gamma}^2+\|Z_3|_{x_2=0}\|_{\frac{3}{2},\gamma}^2 & \lesssim \|Z_3|_{x_2=0}\|_{1,\gamma}^2+(\frac{1}{\e}+C)\vertiii{Z_3}_{1,\gamma}^2+\Sigma_{i\neq2,3}\frac{1}{\e}\vertiii{Z_i}_{1,\gamma}^2\\
&\quad +\frac{1}{\e}\left(\vertiii{T^\gamma_rW^+}_{1,\gamma}^2+\vertiii{W^+}_0^2+\vertiii{F}_{1,\gamma}^2\right).
\end{align*}
For the 1st, 4th, 5th, 6th and 7th equations in \eqref{dep},  we consider the symmetrizer $(T^\gamma_{\Lambda})^*T^\gamma_{\Lambda}$, follow the similar  argument   to the previous pole case, and obtain
\begin{align*}
\gamma\|Z_j\|_{1,\gamma}^2\lesssim C\|Z_j\|_{1,\gamma}^2+\frac{1}{\gamma}\|Z_3\|_{2,\gamma}^2+\Sigma\frac{1}{\gamma}\|Z_i\|_{1,\gamma}^2+\frac{1}{\gamma}\left( \|T^\gamma_{r}W^+\|_{1,\gamma}^2+\|W^+\|_0^2+\|F\|_{1,\gamma}^2\right),
\end{align*}
for $j=1,4,5,6,7$.

For the differential equation of the incoming mode $Z_2$, we take the symmetrizer $T^\gamma_{\Lambda}$, and use the similar strategy to obtain
\begin{align*}
\vertiii{Z_2}_{1,\gamma}^2 & \lesssim\Re\langle Z_2, T^\gamma_{\Lambda}Z_2\rangle|_{x_2=0}+\frac{1}{\e}\vertiii{Z_2}_0^2+\frac{1}{\e}\vertiii{Z_j}_{1,\gamma}^2+\Sigma_{i\neq2,3}\frac{1}{\e}\vertiii{Z_i}_0^2\\
&\quad +\frac{1}{\e}\left(\vertiii{T^\gamma_{r}W^+}_0^2+\vertiii{W^+}_{-1,\gamma}^2+\vertiii{F}_0^2\right).
\end{align*}
Combining the above estimates together leads to the estimate:
\begin{align*}
\frac{1}{\gamma}&\vertiii{Z_3}_{2,\gamma}^2+\gamma\vertiii{Z_j}_{1,\gamma}^2+\gamma\vertiii{Z_2}_{1,\gamma}^2+\frac{1}{\gamma}\|Z_3|_{x_2=0}\|_{\frac{3}{2},\gamma}^2\\
&\lesssim \gamma\Re\langle Z_2, T^\gamma_{\Lambda}Z_2\rangle|_{x_2=0}+\frac{1}{\gamma}\left(\vertiii{T^\gamma_{r}W^+}_{1,\gamma}^2+\vertiii{W^+}_0^2+\vertiii{F}_{1,\gamma}^2\right).
\end{align*}
Similarly, for $Z^-=(Z_8,Z_9,...,Z_14):=T^\gamma_{\chi_1Q^l}T^\gamma_{\chi_{p_1}}W^-$, we have
\begin{align*}
\frac{1}{\gamma}&\vertiii{Z_{10}}_{2,\gamma}^2+\gamma\vertiii{Z_j}_{1,\gamma}^2+\gamma\vertiii{Z_9}_{1,\gamma}^2+\frac{1}{\gamma}\|Z_{10}|_{x_2=0}\|_{\frac{3}{2},\gamma}^2\\
&\lesssim \gamma\Re\langle Z_9, T^\gamma_{\Lambda}Z_9\rangle|_{x_2=0}+\frac{1}{\gamma}\left(\vertiii{T^\gamma_{r}W^+}_{1,\gamma}^2+\vertiii{W^+}_0^2+\vertiii{F}_{1,\gamma}^2\right).
\end{align*}
Therefore, we only need to estimate $\gamma\Re\langle Z_2, T^\gamma_{\Lambda}Z_2\rangle|_{x_2=0}$ and $\gamma\Re\langle Z_9, T^\gamma_{\Lambda}Z_9\rangle|_{x_2=0}$ which can be controlled by $\|Z_{in}|_{x_2=0}\|_{1,\gamma}^2$. The estimate of $\|Z_{in}|_{x_2=0}\|_{1,\gamma}^2$ can be obtained through the boundary condition, using the fact that the Lopatinskii determinant has a positive lower bound in the open neighborhood $\V_{p_1}^1$, that is,
\begin{align*}
\|Z_{in}|_{x_2=0}\|_{1,\gamma}^2\lesssim \|G\|_{1,\gamma}^2+\|Z_{out}|_{x_2=0}\|_{1,\gamma}^2+\|W^\nc|_{x_2=0}\|_0^2.
\end{align*}
Combining all these estimates together, we have
\begin{equation}\label{poleestimate}
\begin{split}
&\frac{1}{\gamma}\vertiii{Z_{out}}_{2,\gamma}^2+\gamma\vertiii{Z_c}_{1,\gamma}^2+\gamma\vertiii{Z_{in}}_{1,\gamma}^2+\frac{1}{\gamma}\|Z_{out}|_{x_2=0}\|_{\frac{3}{2},\gamma}^2+\|Z_{in}|_{x_2=0}\|_{1,\gamma}^2\\
&\lesssim \|G\|_{1,\gamma}^2+\|W^\nc|_{x_2=0}\|_0^2+\frac{1}{\gamma}\left(\vertiii{T^\gamma_{r}W^+}_{1,\gamma}^2+\vertiii{W^+}_0^2+\vertiii{F}_{1,\gamma}^2\right),
\end{split}
\end{equation}
where $Z_c=(Z_1,Z_4,Z_5,Z_6,Z_7,Z_8,Z_{11},Z_{12},Z_{13},Z_{14})^\top$.

\subsection{Case 4 -- others} 

At all of the remaining points, the Lopatinskii determinant does not vanish and the system can be reduced to the non-characteristic form. This is exactly the good frequency case in \cite{coulombel2004stability}. By the same treatment as in \cite{coulombel2004stability}, we can construct the Kreiss symmetrizer. In particular, starting from the  original system \eqref{paralinearizationsystem}, we consider the cut-off symbol $\chi_{re}=1-\bar{\chi}_{p_1}-\bar{\chi}_{p_2}-\bar{\chi}_{rt}$ in $\Gamma^0_k$ for any integer $k$, where $\bar{\chi}_{p_1}$ is the sum of all four cut-off functions $\chi_{p_1}$ for the four neighborhood $\V_{p_1}^i, i=1,2,3,4$, $\bar{\chi}_{p_2}$ is the sum of the two cut-off functions $\chi_{p_2}$ for the two neighborhood $\V_{p_2}^1$ and $\V_{p_2}^2$, and $\bar{\chi}_{rt}$ is the sum of three cut-off functions $\chi_{rt}$ for the three neighborhood $\V_{rt}^i, i=1,2,3$. Since all the neighborhoods above do not overlap,   $\chi_{re}$ is also a cut-off function that equals $0$ near the roots of the Lopatinskii determinant $\Upsilon_r$ and the poles of the system $\Upsilon_p$. Then we can construct an open neighborhood $\V_{re}$ that contains the support of $\chi_{re}$ but does not contain a small neighborhood of $\Upsilon_r$ and $\Upsilon_p$.  Denote
\begin{align*}
&W^\pm_{re}:=T^\gamma_{\chi_{re}}W^\pm, \text{ and }W_{re}:=(W_{re}^+,W_{re}^-)^\top.
\end{align*}
By the same idea in \cite{coulombel2002weak,coulombel2004stability}, we can eliminate all the components of $W^\pm_{re}$ in the kernel of $I_2$ and obtain a differential equation in $W^\nc_{re}:=T^\gamma_{\chi_{re}}W^\nc$ of the following form:
\begin{align}
\p_2W^\nc_{re} = T^\gamma_{\chi_2\mathbb{A}}W^\nc+T^\gamma_{\mathbb{E}}W^\nc+T^\gamma_{r}W+R_0F+R_{-1}W,
\end{align}
where $\mathbb{A}=\text{diag}\{\mathbb{A}^r, \mathbb{A}^l\}$ and $\mathbb{E}$ is a symbol in $\Gamma^0_1$. Same as before, $r$ is a symbol in $\Gamma^0_1$ which supports only in the place where $\chi_{re}\in (0,1)$. Combining with the boundary conditions in \eqref{paralinearizationsystem}, we can construct the Kreiss symmetrizer as in  \cite[Proposition 3.1]{coulombel2002weak}. By a standard argument in \cite{coulombel2002weak, coulombel2004stability} for a boundary value problem on the half space, we can obtain the following estimate
\begin{equation}\label{restestimate}
\begin{split}
&\gamma\vertiii{W_{re}}_{1,\gamma}^2+\|W^\nc_{re}|_{x_2=0}\|_{1,\gamma}^2\\
& \lesssim \|G\|_{1,\gamma}^2 +\|W^\nc|_{x_2=0}\|_0^2+\frac{1}{\gamma}\left(\vertiii{F}_{1,\gamma}^2+\vertiii{W}_0^2+\vertiii{T^\gamma_r W}_{1,\gamma}^2\right).
\end{split}
\end{equation}

\section{Proof of the Main Theorem}\label{ProofofMain}
Now we combine all the estimates in the above four cases. First we point out that the left   sides of \eqref{polerootestimate}, \eqref{rootestimate}, \eqref{poleestimate} and \eqref{restestimate} dominate $\gamma^3\vertiii{W_{\bar{\chi}_{p_1}}}_{0}^2+\gamma^2\|W^\nc_{\bar{\chi}_{p_1}}|_{x_2=0}\|_{0}^2$, $\gamma^3\vertiii{W_{\bar{\chi}_{rt}}}_{0}^2+\gamma^2\|W^\nc_{\bar{\chi}_{rt}}|_{x_2=0}\|_{0}^2$, $\gamma^3\vertiii{W_{\bar{\chi}_{p_2}}}_{0}^2+\gamma^2\|W^\nc_{\bar{\chi}_{p_2}}|_{x_2=0}\|_{0}^2$,  and $\gamma^3\vertiii{W_{\bar{\chi}_{re}}}_{0}^2+\gamma^2\|W^\nc_{\bar{\chi}_{re}}|_{x_2=0}\|_{0}^2$ respectively. We note that
\begin{align*}
&\gamma^3\vertiii{W_{\bar{\chi}_{p_1}}}_{0}^2+\gamma^2\|W^\nc_{\bar{\chi}_{p_1}}|_{x_2=0}\|_{0}^2+\gamma^3\vertiii{W_{\bar{\chi}_{rt}}}_{0}^2+\gamma^2\|W^\nc_{\bar{\chi}_{rt}}|_{x_2=0}\|_{0}^2+\gamma^3\vertiii{W_{\bar{\chi}_{p_2}}}_{0}^2\\
&\quad +\gamma^2\|W^\nc_{\bar{\chi}_{p_2}}|_{x_2=0}\|_{0}^2+\gamma^3\vertiii{W_{\bar{\chi}_{re}}}_{0}^2+\gamma^2\|W^\nc_{\bar{\chi}_{re}}|_{x_2=0}\|_{0}^2\\
&=\gamma^3\vertiii{W}_0^2+\gamma^2\|W^\nc|_{x_2=0}\|_0^2.
\end{align*}
Summing up all these equations and taking $\gamma$ large enough, we can absorb $\|W^\nc|_{x_2=0}\|_0^2$ and $\frac{1}{\gamma}\vertiii{W}_0^2$ into the terms on the left   side of the sum of equations \eqref{polerootestimate}, \eqref{rootestimate}, \eqref{poleestimate} and \eqref{restestimate}. The only thing left   is to absorb $\frac{1}{\gamma}\vertiii{T^\gamma_r W}_{1,\gamma}^2$ into the terms on the left   side of the sum of equations \eqref{polerootestimate}, \eqref{rootestimate}, \eqref{poleestimate} and \eqref{restestimate}. Following the same argument in \cite{coulombel2004stability,coulombel2002weak}, we observe that the support of $r$ in each case is contained in the following set:
 \begin{align*}
 \left\{(t,x_1,x_2,\delta,\eta)\in\R^3_+\times\Pi; \bar{\chi}_{p_1}\in(0,1) \text{ or } \bar{\chi}_{p_2}\in(0,1) \text{ or } \bar{\chi}_{rt}\in(0,1) \text{ or } \chi_{re}\in(0,1)\right\},
 \end{align*}
 where we recall that the symbol $r$ is introduced in \eqref{para1}.
Since  $\bar{\chi}_{p_1}+\bar{\chi}_{p_2}+\bar{\chi}_{rt}+\chi_{re}=1$, then $r=0$ if one of $\bar{\chi}_{p_1}$, $\bar{\chi}_{p_2}$, $\bar{\chi}_{rt}$ and $\chi_{re}$ equals 1. Moreover,  $\sigma$ only vanishes at some points in the set where $\bar{\chi}_{p_1}=1$ and $\bar{\chi}_{rt}=1$. Thus $\sigma$ has a positive lower bound on the support of $r$ and $r$ can be decomposed into
\begin{align*}
r=a_{p_2}\bar{\chi}_{p_2}+a_{re}\chi_{re}+a_{p_1}\sigma\chi_1^{p_1}
\begin{pmatrix}
Q^r & 0\\
0 & Q^l
\end{pmatrix}\bar{\chi}_{p_1}+a_{rt}\sigma\chi_1^{rt}
\begin{pmatrix}
Q^r & 0\\
0 & Q^l
\end{pmatrix}\bar{\chi}_{rt},
\end{align*}
where $a_{p_2}$, $a_{re}$, $a_{p_1}$ and $a_{rt}$ are some block diagonal symbol in the class $\Gamma^0_1$; and $\chi_1^{p_1}$ and $\chi_1^{rt}$ denote the corresponding cut-off functions $\chi_1$ in the Case 1 and Case 2 respectively. This suggests that $\frac{1}{\gamma}\vertiii{T^\gamma_r W}_{1,\gamma}^2$ can be absorbed by 
\begin{align*}
\gamma\vertiii{T^\gamma_\sigma T^\gamma_{\bar{\chi}_{p_1}}W}_0^2+\gamma\vertiii{T^\gamma_\sigma T^\gamma_{\bar{\chi}_{rt}}W}_0^2+\gamma\vertiii{ T^\gamma_{\bar{\chi}_{p_2}}W}_{1,\gamma}^2+\gamma\vertiii{ T^\gamma_{\chi_{re}}W}_{1,\gamma}^2,
\end{align*}
with $\gamma$ large enough. The above sum can also be dominated by the sum of the left   sides of equations \eqref{polerootestimate}, \eqref{rootestimate}, \eqref{poleestimate} and \eqref{restestimate} with $\gamma$ large enough. Therefore by adding the equations \eqref{polerootestimate}, \eqref{rootestimate}, \eqref{poleestimate} and \eqref{restestimate}, we obtain
\begin{align*}
\|W^{\nc}|_{x_2=0}\|_0^2\leq C_0\left(\frac{1}{\gamma^3}\vertiii{\tilde{F}}_{1,\gamma}^2+\frac{1}{\gamma^2}\|\tilde{G}\|_{1,\gamma}^2\right),
\end{align*}
and thus complete the proof of Theorem \ref{variablemain}.

\bigskip
\section*{Acknowledgements}
R. M. Chen's research was supported in part   by the Simons Foundation under Grant 354996 and the NSF grant DMS-1613375.
J. Hu's research was supported in part by  the NSF Foundation under Grant DMS-1312800.
D. Wang was supported in part by the NSF grants DMS-1312800 and DMS-1613213.

The authors would like to thank the anonymous referees for valuable comments and suggestions.


\bigskip\bigskip
\begin{thebibliography}{10}

\bibitem{alinhac1989existence}
{ S.~Alinhac}, {\em Existence d'ondes de rarefaction pour des systemes quasi-lin\'eaires hyperboliques multidimensionnels}. Comm. Partial Differential Equations 14 (1989), no. 2, 173--230.

\bibitem{benzoni2007multi}
{ S.~Benzoni-Gavage and D.~Serre}, {\em Multidimensional Hyperbolic Partial
 DifferentialEquations: First-order Systems and Applications}. Oxford Mathematical Monographs. The Clarendon Press, Oxford University Press, Oxford, 2007.

\bibitem{bony1981calcul}
{ J. M.~Bony}, {\em Calcul symbolique et propagation des singularit\'{e}s pour les \'{e}quations aux d\'{e}riv\'{e}es partielles non lin\'{e}aires}. Ann. Sci. \'{E}cole Norm. Sup. (4) 14 (1981), no. 2, 209--246.   

\bibitem{chazarain1982introduction}
{ J.~Chazarain and A.~Piriou}, {\em Introduction to the Theory of Linear Partial Differential Equations}. Studies in Mathematics and its Applications, 14. North-Holland Publishing Co., Amsterdam-New York, 1982. 

\bibitem{chen2008existence}
{ G.-Q.~Chen and Y.-G.~Wang}, {\em Existence and stability of compressible
 current-vortex sheets in three-dimensional magnetohydrodynamics}.  Arch. Ration. Mech. Anal. 187 (2008), no. 3, 369--408.

\bibitem{chen2017relativistic}
{G.-Q.~Chen, P.~Secchi, and T.~Wang}, {\em Nonlinear stability of relativistic vortex sheets in three-dimensional Minkowski spacetime}. arXiv:1707.02672 [math.AP],  2017.
  
 \bibitem{chen2015linear}
{ R. M.~Chen, J.~Hu, and D.~Wang}, {\em Linear stability of compressible vortex sheets in two-dimensional elastodynamics}. Adv. Math. 311 (2017), 18--60. 
  

\bibitem{chen2018nonlinear}
{ R. M.~Chen, J.~Hu, and D.~Wang}, {\em Nonlinear stability and existence of compressible vortex sheets in two-dimensional elastodynamics}. Preprint, 2018. 

\bibitem{coulombel2002weak}
{ J.-F.~Coulombel}, {\em Weak stability of nonuniformly stable
  multidimensional shocks}. SIAM J. Math. Anal. 34 (2002), no. 1, 142--172 (electronic).

\bibitem{coulombel2004weakly}
{ J.-F.~Coulombel}, {\em Weakly stable
  multidimensional shocks}. Ann. Inst. H. Poincar\'{e} Anal. Non Lin\'{e}aire 21 (2004), no. 4, 401--443.

\bibitem{coulombel2004stability}
{ J.-F.~Coulombel and P.~Secchi}, {\em The stability of compressible vortex
  sheets in two space dimensions}. Indiana Univ. Math. J. 53 (2004), no. 4, 941--1012.

\bibitem{dafermos2010hyperbolic}
{ C. M. ~Dafermos}, {\em Hyperbolic Conservation Laws in Continuum Physics}. 4th edition. Grundlehren der Mathematischen Wissenschaften [Fundamental Principles of Mathematical Sciences], 325. Springer-Verlag, Berlin, 2016.

\bibitem{francheteau2000existence}
{ J.~Francheteau and G.~M{\'e}tivier}, {\em Existence de chocs faibles pour
  des systemes quasi-lin{\'e}aires hyperboliques multidimensionnels}. Ast\'{e}risque No. 268 (2000).
  
 \bibitem{GuWang}
 X. Gu, F. Wang,
 {\em  Well-posedness of the free boundary problem in incompressible elastodynamics under the mixed type stability condition}.
 arXiv:1804.00972 [math.AP], 2018.


\bibitem{gurtin1981introduction}
{M. E.~Gurtin}, {\em An Introduction to Continuum Mechanics}. Mathematics in Science and Engineering, 158. Academic Press, Inc. [Harcourt Brace Jovanovich, Publishers], New York-London, 1981.
  
 
 \bibitem{hu2011global}
{X.~Hu and D.~Wang}, {\em Global existence for the multi-dimensional compressible viscoelastic flows}. 
J. Differential Equations 250 (2011), no. 2, 1200--1231.
 
 \bibitem{joseph1990fluid}
{D.~Joseph}, {\em Fluid Dynamics of Viscoelastic Liquids}. Applied Mathematical Sciences, 84. Springer-Verlag, New York, 1990.
 
 \bibitem {LiWangZhang}
 H. Li, W. Wang, Z. Zhang,
 {\em Well-posedness of the free boundary problem in incompressible elastodynamics}. 
arXiv:1802.08819 [math.AP], 2018.
 
 \bibitem{meyer1981remarques}
{ Y.~Meyer}, {\em Remarques sur un th\'{e}or\`{e}me de J.-M. Bony}. 
Rend. Circ. Mat. Palermo (2) 1981, suppl. 1, 1--20.
 
  \bibitem{morando2008two}
{ A.~Morando and P.~Trebeschi}, {\em Two-dimensional vortex sheets for the
  nonisentropic Euler equations: linear stability}. J. Hyperbolic Differ. Equ. 5 (2008), no. 3, 487--518. 

\bibitem{ruanrectilinear}
{ L.~Ruan, D.~Wang, S.~Weng, and C.~Zhu}, {\em Rectilinear vortex sheets of
  inviscid liquid-gas two-phase flow: linear stability}. Commun. Math. Sci.  14 (2016), 735--776.
  
  \bibitem{Serre2000}
{ D.~Serre}, {\em Systems of Conservation Laws  2: Geometric Structures, Oscillations, and Initial-Boundary Value Problems}.
 Cambridge University Press, Cambridge, 2000.
  
 \bibitem{Stevens}
 B. Stevens,  
{\em Short-time structural stability of compressible vortex sheets with surface tension}.   
Arch. Ration. Mech. Anal. 222 (2016), no. 2, 603--730. 

 \bibitem{trakhinin2005existencevairable}
{ Y.~Trakhinin}, {\em Existence of compressible current-vortex sheets:
  variable coefficients linear analysis}. Arch. Ration. Mech. Anal. 177 (2005), 331--366. 
  
  \bibitem{trakhinin2005existencelinearized}
{ Y.~Trakhinin}, {\em On the existence of
  incompressible current--vortex sheets: study of a linearized free boundary
  value problem}. Math. Methods Appl. Sci. 28 (2005), no. 8, 917--945.
  
  \bibitem{trakhinin2009existence}
{ Y.~Trakhinin}, {\em The existence of
  current-vortex sheets in ideal compressible magnetohydrodynamics}. Arch. Ration. Mech. Anal. 191 (2009), no. 2, 245--310. 

\bibitem{wang2013stability}
{ Y.-G.~Wang and F.~Yu}, {\em Stability of contact discontinuities in
  three-dimensional compressible steady flows}. J. Differential Equations 255 (2013),  1278--1356. 
  
\bibitem{wang2015stability}
{ Y.-G.~Wang and H.~Yuan}, {\em Weak stability of transonic contact discontinuities in three-dimensional steady non-isentropic compressible Euler flows}. Z. Angew. Math. Phys., 66 (2015),
341--388.



\end{thebibliography}
\end{document}